\documentclass[a4paper,11pt]{article}
\usepackage[centertags]{amsmath}
\usepackage{amscd}
\usepackage{amsthm}
\usepackage{amssymb}
\usepackage{array}
\usepackage{enumerate}
\usepackage{hyperref}
\usepackage{cleveref}
\usepackage{multicol}
\usepackage{centernot}
\usepackage{mathtools}
\usepackage{tikz-cd}
\usepackage{dynkin-diagrams}
\usepackage{stmaryrd}
\usepackage{pifont}
\usepackage{tikz}
\usetikzlibrary{cd}
\usepackage{tikz-cd}
\usetikzlibrary{matrix}
\usepackage[sc]{mathpazo}

\usepackage[english]{babel}
\usepackage{color}
\usepackage{tikz}
\usepackage{stackengine}
\usepackage{scalerel}
\usepackage{indentfirst} 

\usepackage[T1]{fontenc}

\usepackage[utf8]{inputenc}

\usepackage{geometry}
\pgfkeys{/Dynkin diagram,
	edge length=1.5cm,
	fold radius=1.5cm,
	root radius=0.1cm
}

\newtheorem{defn0}{Definition}[section]
\newtheorem{prop0}[defn0]{Proposition}
\newtheorem{thm0}[defn0]{Theorem}
\newtheorem{lemma0}[defn0]{Lemma}
\newtheorem{corollary0}[defn0]{Corollary}
\newtheorem{example0}[defn0]{Example}
\newtheorem{remark0}[defn0]{Remark}
\newtheorem{conjecture0}[defn0]{Conjecture}
\newtheorem{ques0}[defn0]{Question}

\newenvironment{definition}{ \begin{defn0}}{\end{defn0}}
\newenvironment{proposition}{\bigskip \begin{prop0}}{\end{prop0}}
\newenvironment{theorem}{\bigskip \begin{thm0}}{\end{thm0}}
\newenvironment{lemma}{\bigskip \begin{lemma0}}{\end{lemma0}}
\newenvironment{corollary}{\bigskip \begin{corollary0}}{\end{corollary0}}
\newenvironment{example}{ \begin{example0}\rm}{\end{example0}}
\newenvironment{remark}{ \begin{remark0}\rm}{\end{remark0}}

\usepackage{fancyhdr}
\newcommand{\Z}{\mathbb{Z}}
\newcommand{\R}{\mathbb{R}}

\newcommand{\RP}{\mathbb{R}P}
\newcommand{\C}{\mathbb{C}}
\newcommand{\D}{\mathbb{D}}
\newcommand{\CP}{\mathbb{C}P}
\newcommand{\im}{\operatorname{Im}}
\renewcommand{\ker}{\operatorname{Ker}}

\newcommand{\Ext}{\operatorname{Ext}}

\newcommand{\Hom}{\operatorname{Hom}}
\newcommand{\Isom}{\operatorname{Isom}}

\newcommand{\Ogroup}{\operatorname{O}}
\newcommand{\SO}{\operatorname{SO}}

\newcommand{\Spin}{\operatorname{Spin}}
\newcommand{\U}{\operatorname{U}}
\newcommand{\SU}{\operatorname{SU}}
\newcommand{\Sp}{\operatorname{Sp}}

\newcommand{\Ad}{\operatorname{Ad}}

\newcommand{\Aut}{\operatorname{Aut}}

\newcommand{\acts}{\curvearrowright}
\newcommand{\diag}{\operatorname{diag}}
\newcommand{\id}{\mathrm{id}}

\title{Spin structures on cohomogeneity one manifolds}

\author{
\begin{tabular}{cc}
Jordi Daura Serrano   & Marie-Am\'elie Lawn\\
\small  University of Cyprus &\small Imperial College London (UK) \\
\small daura-serrano.jordi@ucy.ac.cy   &\small  m.lawn@imperial.ac.uk
\end{tabular}
}
\date{\today}

\begin{document}





\selectlanguage{english}
\setcounter{page}{1}




\maketitle
\begin{abstract} 
We study the existence of spin and invariant spin structures on cohomogeneity one manifolds, obtaining explicit criteria in terms of their isotropy data and fundamental groups. We treat all four types of cohomogeneity one actions without exceptional orbits, under suitable connectedness assumptions. 
We illustrate our results on several families of manifolds, notably showing that the natural cohomogeneity one actions on certain Brieskorn homotopy spheres, including exotic spheres, do not lift to their spin structures.
\end{abstract}

\section{Introduction}

Spin structures have played an important role in differential geometry since the work of Lichnerowicz on Dirac operators and the Atiyah--Singer index theorem. Besides their topological applications, spinorial methods have proved particularly useful in the study of geometric structures, notably through the relation between special spinors, holonomy and curvature \cite{Wan89,Bar93,Agr06}. It is well-known that an oriented manifold admits a spin structure if and only if its second Stiefel--Whitney class vanishes. Although this gives a complete characterisation, the obstruction can be difficult to compute. When the manifold admits isometric, orientation preserving group actions of compact Lie groups, one can often exploit these symmetries to obtain more explicit, algebraic criteria. In this situation further question arises. Namely whether the action lifts to an action on the total space of the spin bundle. This is stronger than merely preserving the cohomology class representing the isomorphism class of a spin structure, which is automatic when $G$ is connected. If a spin structure admits a lift of the group $G$ we say that the spin structure is $G$-invariant.

Both questions can be studied particularly well for homogeneous spaces $G/H$. The existence of a spin structure can be studied using the isotropy representation and the topology of the groups $G$ and $H$ \cite{alekseevsky2019spin}, while a $G$-invariant spin structure exists if and only if the isotropy representation lifts to the spin group \cite{alekseevsky2019spin,daura2022g}. In particular, even when $G$ is connected and the manifold has a unique spin structure, the action need not lift \cite{daura2022g}. When it does, one can exploit the group
action in the study of spinors and spinorial equations (see, for example, \cite{agricola2023invariant,artacho2025invariant}).

A natural setting in which to investigate these questions beyond homogeneous spaces is provided by cohomogeneity one manifolds. Although the group action is no longer transitive, the orbit structure remains sufficiently simple to admit a concrete description in terms of Lie groups and their isotropy subgroups. 
By Mostert's description \cite{mostert1957compact} (see also \cite[Chapter IV. \S 8]{bredon1972introduction}),
compact cohomogeneity one manifolds have orbit
space either a circle or an interval. In the first case the manifold is the total space of a fiber bundle over the circle, while in the second case corresponds to two disk bundles glued along a principal orbit. Similarly, non-compact cohomogeneity one manifolds have orbit either a the real line or the nonnegative real half-line. In the first case, the manifold is the product of a homogeneous space and the real line, and the second case is a disk bundle.

This structure has made cohomogeneity one manifolds a particularly interesting class of examples in Riemannian geometry, notably in the study of positive and non-negative curvature \cite{grove2008positively,grove2006nonnegative}. Moreover, invariant geometric equations can often be reduced to ordinary differential equations along a geodesic orthogonal to the orbits, an approach which has been particularly fruitful in the construction of metrics with special holonomy (see for instance
\cite{hitchin2001stable,cleyton2002cohomogeneity}).

It is therefore natural to ask to what extent the existence and invariance of spin structures, can be read directly from the groups describing the action. 
For cohomogeneity one manifolds, however, the
situation is more delicate. The existence of
spin structures depends not only on the isotropy
representations, but also on the topology of
the non-principal orbits and their normal bundles.
In particular, when the orbit space is an
interval, one must determine whether spin
structures on the two disk bundles agree along
their common boundary. 
This gluing problem is especially subtle 
when the non-principal orbits have codimension two
or the principal isotropy group is disconnected.

To our knowledge, no systematic treatment of spin structures and invariant spin structures on cohomogeneity one manifolds has appeared so far. The aim of this paper is to give criteria for their existence directly in terms of the groups describing the action, their isotropy representations and the induced maps on fundamental groups.

In Sections~2 and~3, we recall the topology of
cohomogeneity one manifolds and establish some preliminary results on spin structures. These include an equivariant version of the two-out-of-three lemma and a study of spin structures on quotients by free group actions. We also revisit the criteria for spin and invariant spin structures on homogeneous spaces, which will be used throughout the paper.

We study all four types of cohomogeneity without exceptional orbits, classified according to their orbit spaces:
Type (1), when $M/G\cong S^1$; Type (2), when
$M/G\cong[-1,1]$; Type (3), when $M/G\cong\R$;
and Type (4), when $M/G\cong[0,\infty)$.
We consider these cases in increasing order of
difficulty, beginning with manifolds without
singular orbits and ending with the gluing
problem for manifolds with two singular orbits.

In Section~4, we first treat Types (3) and (1).
For Type (3), the manifold is a product
$G/H\times\R$, and the existence of spin and
$G$-invariant spin structures reduces to the
corresponding questions on $G/H$, see \cref{main_theorem (3) section}.
For Type (1), the manifold is the total space of a fiber bundle over the circle with fiber $G/H$. We prove that, provided $M$ is orientable
and the principal isotropy group $H$ is connected,
$M$ is spin (respectively, $G$-invariantly spin)
if and only if its principal orbit $G/H$ is, see \cref{main_theorem (1) section}. We also show that the connectedness assumption
is essential by constructing a homogeneous
space $G/H$ with exactly two spin structures and an
orientation-preserving $G$-equivariant
diffeomorphism which exchanges them (\cref{disconnected isotropy counterexample}). The corresponding mapping tori admits a cohomogeneity one action of $G$ and it is not spin.

In Section~5, we consider Type (4), where
$M=G\times_K\R^{l+1}$ has a single non-principal
orbit $G/K$. Assuming that $G$ and $K$ are
connected, we obtain criteria for both spin and
$G$-invariant spin structures (\cref{main_theorem noncompact section} and \cref{main_theorem noncompact invariant section}), with the following
consequences:
\begin{itemize}
\item If the non-principal orbit has codimension at
least three, $M$ is spin (respectively,
$G$-invariantly spin) if and only if its principal
orbit $G/H$ is.
\item If the non-principal orbit has codimension two
and $H$ has an even number of connected components,
the same equivalence holds with the non-principal
orbit $G/K$ instead.
\end{itemize}
When $H$ has an odd number of connected components,
the situation is different: in particular, $M$
and $G/K$ cannot both admit $G$-invariant spin
structures. We give precise criteria in this case
in terms of the isotropy and normal representations.

In Section~6, we finally address Type (2), which is the most delicate case, where $M$ is obtained by gluing two disk bundles $M_+$ and $M_-$ along their common principal orbit $G/H$. Assuming that
the non-principal isotropy groups $K_+$ and $K_-$ are connected, we can study the existence and invariance of spin structures on $M_+$ and $M_-$ using the results obtained for cohomogeneity one manifolds of Type (4). Then, we prove the following:
\begin{itemize}
\item If at least one non-principal orbit has
codimension at least three, $M$ is spin
(respectively, $G$-invariantly spin) if and
only if both $M_+$ and $M_-$ are.
\item If both non-principal orbits have codimension two and $H$ has an odd number of connected
components, then $G$-invariant spin structures on $M_+$ and $M_-$ always glue.
\end{itemize}

The general gluing problem is considerably more
subtle: even when both disk bundles are spin, their spin structures need not agree on the common boundary. We explicitly determine which
spin structures on the principal orbit extend
to each disk bundle, obtaining necessary and
sufficient conditions for their compatibility.
When both non-principal orbits have codimension two,
the obstructions take different forms depending
on the parity of the number of connected components of $H$. We also establish an equivariant gluing lemma and obtain corresponding criteria for $G$-invariant spin  structures. These results are proved through Section~6, in particular, we refer to \cref{gluing spin structures a even}, \cref{gluing spin structures a odd}, \cref{gluing spin structures a even}, \cref{gluing G-inv spin structures} and \cref{gluing G-invariant spin structures a even}.
We show that these obstructions are genuine by constructing examples where both disk bundles admit spin structures, even $G$-invariant ones, but their union is not spin.

We conclude the paper in Section~7 by illustrating our results on several families of cohomogeneity one manifolds. We recover the classical spin conditions for spheres, complex projective spaces, and real, complex and
quaternionic Grassmannians, and investigate the invariance of their spin structures under the
corresponding group actions. We also apply our criteria in the opposite direction, deducing that the homogeneous space $\Spin(7)/(S^1\cdot\SU(3))$ is spin, and show that
an entire family of Bazaikin spaces, including
the Berger space, is not spin. We show that the natural $S^1\times\SO(n)$-action on the Brieskorn
manifolds $W_d^{2n-1}$ does not admit an invariant spin structure when $n$ and $d$ are odd. 

\paragraph{Acknowledgments.}
We would like to thank Wilderich Tuschmann for bringing to our attention several connections between cohomogeneity one actions, invariant spin structures and curvature, and in particular for suggesting that we investigate the natural actions on Brieskorn manifolds.

\textbf{Funding.} The first author was supported by the research grant EXCELLENCE/0925-Hopf problem and monotone invariants, under the program “THALIA 2021-2027”, co-funded by the European Union through the Research and Innovation Foundation of Cyprus, and
by the grant PID2023-147642NB-I00 from the Spanish Ministry of Science, Innovation and Universities.

The authors thank the Heilbronn Institute for Mathematical Research (HIMR) for supporting, through its Small Grants Scheme, the first author's research visit to Imperial College London, during which much of this work was carried out.

\paragraph{AI use.} All mathematical results and their original proofs were written without the help of generative AI. AI tools were subsequently used to assist with reviewing and checking proofs, and improving the exposition. The authors take full responsibility for the content of the manuscript.

\section{Preliminaries on the topology of cohomogeneity one manifolds}\label{sec: preliminaris cohomogeneity one topology}

The aim of this section is to review some results on the topology of cohomogeneity one manifolds, focusing on computing their fundamental group.

\subsection{Principal orbits, isotropy types and cohomogeneity}

Let $M$ be an $n$-dimensional connected smooth manifold, and let $G$ be a compact connected Lie group acting smoothly on $M$. Since $G$ is compact, there exists a Riemannian metric on $M$ such that $G$ acts by isometries. For completeness, we briefly recall some standard facts about the cohomogeneity of a group action; see \cite[Chapters~I--II]{bredon1972introduction} and \cite[Chapters~3--4]{kawakubo1991theory} for details.

For $x\in M$, denote by $G_x:=\{g\in G \mid gx=x\}$ the isotropy subgroup of $G$ at $x$ and by $\mathcal{O}_x:=\{gx\mid g\in G\}$ the orbit of $x$. The orbit $\mathcal{O}_x$ is an embedded Riemannian submanifold of $M$. The map $G/G_x\longrightarrow \mathcal{O}_x$ defined by $gG_x\mapsto gx$ is a $G$-equivariant diffeomorphism. Moreover, if $y\in\mathcal{O}_x$, then $G_x$ and $G_y$ are conjugate subgroups of $G$.

We denote by $[H]$ the conjugacy class of a subgroup $H$ of $G$ and by $\mathcal{G}$ the set of conjugacy classes of subgroups of $G$. We define a partial order on $\mathcal{G}$ by setting $[H]\leq[K]$ if and only if $H$ is conjugate to a subgroup of $K$. The conjugacy class $[G_x]$ is called the isotropy type of $x$, or equivalently of the orbit $\mathcal{O}_x$. Note that two points having the same isotropy type do not necessarily belong to the same orbit, although their orbits are $G$-equivariantly diffeomorphic.

The principal orbit theorem (see \cite[IV.3.1]{bredon1972introduction}) asserts that if $G$ is a compact Lie group acting smoothly on a connected smooth manifold $M$, then there exists a unique minimal isotropy type with respect to the partial order introduced above. The orbits with this minimal isotropy type are called principal or regular, whereas other orbits are called exceptional if they have the same dimension as principal orbits, and singular otherwise. We denote by $M_0$ the set of all points of $M$ whose orbit is principal. The set $M_0$ is an open dense submanifold of $M$.

The cohomogeneity of an action is the codimension of its principal orbits. In particular, a connected manifold admits a cohomogeneity zero action if and only if it is a homogeneous space. Cohomogeneity one manifolds can therefore be viewed as the simplest non-homogeneous generalisation of homogeneous spaces, while retaining strong group-theoretic control over their geometry.

The following well-known result of Mostert (see \cite{mostert1957compact} and \cite[Chapter IV.8]{bredon1972introduction}) describes the structure of cohomogeneity one manifolds:

\begin{theorem}\label{cohomogenity one structure}
\cite{mostert1957compact}
Let $M$ be an $n$-dimensional connected smooth manifold without boundary, with a cohomogeneity one action of a compact Lie group $G$. Let $G/H$ be a principal orbit. If $M$ is compact, then we have one of the following two cases:
\begin{itemize}
    \item[(1)] The quotient is $M/G\cong S^1$. Then $M$ is equivariantly diffeomorphic to the total space of a fiber bundle $\pi:M\longrightarrow S^1$ with fiber $G/H$ and structure group $N_G(H)/H$. The action has no non-principal orbits.

    \item[(2)] Then there exist subgroups $K_+$ and $K_-$ of $G$ containing $H$ and satisfying $K_\pm/H\cong S^{l_\pm}$ for some $l_{+}, l_{-} \geq 0$. The actions of $K_\pm$ on $S^{l_\pm}$ are linear and thus extend to linear actions of $K_\pm$ on $\mathbb{D}^{l_\pm+1}$. The non-principal orbits of the $G$-action on $M$ are $G / K_\pm$ and correspond to the preimages of $\pm 1 \in [-1,1]$, respectively. Moreover, $M$ is $G$-equivariantly diffeomorphic to the union of disk bundles over $G/K_\pm$ glued along the common boundary $G/H$, 
		$$M \cong \left( G\times_{K_+}\mathbb{D}^{l_++1} \right) \cup_{G/H} \left( G\times_{K_-}\mathbb{D}^{l_-+1} \right).$$

\end{itemize}
If $M$ is not compact, then we have one of the following two cases:
\begin{itemize}
    \item[(3)] The quotient is $M/G\cong\R$. Then $M$ is equivariantly diffeomorphic to $G/H\times\R$ with the trivial action on $\R$.

    \item[(4)] The quotient is $M/G\cong[0,\infty)$. Then $M$ is equivariantly diffeomorphic to a tube,
    \[
    M\cong G\times_K V,
    \]
    where $V\cong\R^{l+1}$ is a Euclidean space on which $K$ acts orthogonally and transitively on the unit sphere $S^l\subseteq\R^{l+1}$ with isotropy subgroup conjugate to $H$ in $K$ (hence $K/H\cong S^l$). Here $l\geq0$, with $l=0$ corresponding to an exceptional orbit and $l\geq1$ to a singular orbit. The non-principal orbit $G/K$ corresponds to the preimage of $0\in[0,\infty)$.
\end{itemize}
\end{theorem}

The above theorem gives the general classification of cohomogeneity one manifolds, including actions with exceptional orbits. In the remainder of the paper, we restrict our study of spin structures to actions without exceptional orbits, so that $l_\pm\geq 1$ in case (2) and $l\geq 1$ in case (4). Note that these restrictions are automatically satisfied when the corresponding isotropy groups $K_\pm$ or $K$ are connected, as will be assumed in our main results.

\begin{remark}\label{type (2) commutative diagram}
Note that a cohomogeneity one action of $G$ on $M$ of type (2) induces a commutative diagram of inclusions of Lie groups 
\[\begin{tikzcd}
    & G& \\
    K_+ \ar{ru}{j_+} & & K_-\ar{lu}[xshift=3ex,yshift=2.5ex]{j_-}\\
    & H\ar{lu}{i_+}\ar{ru}[xshift=3ex,yshift=-3ex]{i_-} &
\end{tikzcd}\]
satisfying $K_\pm/H\cong S^{l_\pm}$. Conversely, compact Lie groups $G$, $K_+$, $K_-$ and $H$ fitting in the above diagram and satisfying $K_\pm/H\cong S^{l_\pm}$ determine a manifold with a cohomogeneity one action of $G$ of type (2) given by $G\times_{K_+}\D^{l_++1}\cup_{G/H} G\times_{K_-}\D^{l_-+1}$. Hence a cohomogeneity one action of type (2) is fully determined by the above commutative diagram.

Similarly, a cohomogeneity one action of $G$ on $M$ of type (4) is fully determined by the inclusions of 3 compact Lie groups 
\[\begin{tikzcd}
H \ar{r}{i} & K \ar{r}{j} & G
\end{tikzcd}\]
where $K/H\cong S^l$. 

We refer to $(H,K,G,i,j,l)$ as a cohomogeneity one
tuple of type (4).
\end{remark}
\begin{remark}
Note that a manifold can admit cohomogeneity one actions of different types. For example, let $M=S^n\times S^1$, with $n\geq2$. The action of $\SO(n+1)$ on $M$ which acts transitively on the $S^n$ and trivially on the $S^1$ is a cohomogeneity one action of type (1). On the other hand, the action of $\SO(n)\times S^1$ on $M$, where the $\SO(n)$-factor acts by cohomogeneity one on $S^n$ and the $S^1$-factor acts by left multiplication on $S^1$, is a cohomogeneity one action of type (2) on $M$.
\end{remark}

\subsection{Fundamental groups of cohomogeneity one manifolds}
Our next goal is to recall some facts on the fundamental group of manifolds admitting cohomogeneity one actions. 

Assume that \(M\) is a manifold with a cohomogeneity one action of a compact connected Lie group \(G\) of type (1). Then \(M\) is a fiber bundle over \(S^1\) with fiber \(G/H\). The long exact sequence of homotopy gives a short exact sequence
\[
1\longrightarrow \pi_1(G/H)\longrightarrow \pi_1(M)\longrightarrow \mathbb Z\longrightarrow 1.
\]
Since \(\mathbb Z\) is free, this sequence splits. For a suitable choice of splitting, we have
\[
\pi_1(M)\cong \pi_1(G/H)\rtimes_{\phi}\mathbb Z,
\]
where \(\phi:\mathbb Z\to\operatorname{Aut}(\pi_1(G/H))\) is the action induced by the \(G\)-equivariant diffeomorphism of \(G/H\) defining the bundle. In particular, \(\pi_1(M)\) surjects onto \(\mathbb Z\), so a simply connected manifold cannot admit a cohomogeneity one action of type (1).

Assume that \(M\) is a manifold with a cohomogeneity one action of a compact connected Lie group \(G\) of type (3). Then \(M\cong \mathbb{R}\times G/H\), and hence \(\pi_1(M)\cong \pi_1(G/H)\).

Assume that \(M\) is a manifold with a cohomogeneity one action of a compact connected Lie group \(G\) of type (4). Then \(M\cong G\times_K\mathbb{R}^{l+1}\), and there is a strong deformation retraction of \(G\times_K\mathbb{R}^{l+1}\) onto \(G\times_K\{0\}\cong G/K\). Consequently, \(M\) and \(G/K\) are homotopy equivalent, and \(\pi_1(M)\cong\pi_1(G/K)\).

The remaining case, when $M$ admits a cohomogeneity one action of a compact connected Lie group $G$ of type (2), is significantly more involved. To study this case, we need to introduce some notation. We set
$$M_\pm:=G\times_{K_\pm}\D^{l_\pm+1}$$
and
$$M_\pm':=G\times_{K_\pm}\partial(\D^{l_\pm+1})=G\times_{K_\pm}S^{l_\pm}.$$

We also define the inclusions
$$\iota_\pm:M'_\pm\longrightarrow M_\pm$$
and the retractions
$$r_\pm:M_\pm\longrightarrow G\times_{K_\pm}\{0\}=G/K_\pm.$$

Finally, since $S^{l_\pm}$ is $K_\pm$-equivariantly diffeomorphic to $K_\pm/H$, we can define $G$-equivariant diffeomorphisms
$$\zeta_\pm:G/H\longrightarrow G\times_{K_\pm}(K_\pm/H)\cong M_\pm'$$
satisfying $\zeta_\pm(gH)=[g,eH]$, where $e$ denotes the identity element of $G$. The inverse maps are given by $\zeta_\pm^{-1}([g,k_\pm H])=gk_\pm H$.

On the other hand, the inclusions of Lie groups $i_\pm:H\longrightarrow K_\pm$ and $j_\pm:K_\pm\longrightarrow G$ induce homogeneous fibrations
\[
\begin{tikzcd}
K_\pm/H\ar{r}{\overline{j_\pm}} & G/H \ar{r}{\overline{i_\pm}} & G/K_\pm.
\end{tikzcd}
\]

The next lemma provides a relation between the above homogeneous fibration and the inclusions $\iota_\pm$.

\begin{lemma}\cite[Proposition 1.8]{hoelscher2010classification}\label{iota=i}
The maps $\overline{i_\pm}:G/H\longrightarrow G/K_\pm$ and
$r_\pm\circ\iota_\pm\circ\zeta_\pm:G/H\longrightarrow G/K_\pm$
are equal. In particular, $\iota_{\pm*}$ is surjective, and it is an isomorphism if $l_\pm\geq 2$.
\end{lemma}

Since we can write $M$ as the union $M_+\cup_{G/H}M_-$, the Seifert--van Kampen theorem yields the following result.

\begin{lemma}\cite[Proposition 1.8]{hoelscher2010classification}\label{cohomogenity one fundamental group}
Let $M$ be a manifold with a cohomogeneity one action of a compact connected Lie group $G$ of type (2). Then
$$\pi_1(M)\cong \pi_1(G/H)/(\ker(\overline{i_+}_*)\cdot\ker(\overline{i_-}_*)).$$
\end{lemma}

\begin{corollary}\cite[Corollary 1.12]{hoelscher2010classification}\cite[Lemma 1.6]{grove2008positively}\label{K connected}
Let $M$ be a manifold with a cohomogeneity one action of a compact connected Lie group $G$ of type (2).
\begin{itemize}
    \item[1.] If $l_+>1$ and $l_-\geq 1$, then $\pi_1(M)\cong\pi_1(G/K_-)$. If $K_+$ is connected, then $H$ is connected.
    \item[2.] If $l_\pm>1$, then $\pi_1(M)\cong\pi_1(G/H)\cong\pi_1(G/K_\pm)$. If $M$ is simply connected, then $H$, $K_+$ and $K_-$ are connected.
\end{itemize}
\end{corollary}

\subsection{Some additional results for codimension-two singular orbits}

It remains to study the case where $l_\pm=1$. In this subsection, we assume that $K_\pm$ are connected. We need to introduce some extra notation that will also be used when studying (invariant) spin structures in this setting. Recall that the transitive actions of $K_\pm$ on $S^1$ are given by surjective Lie group homomorphisms $q_\pm:K_\pm\longrightarrow S^1$. Define the positive integer $a$ and the nonnegative integers $b_+$ and $b_-$ satisfying
$$q_{\pm*}(\pi_1(K_\pm))=a\Z$$
and
$$q_{\pm*}(\ker j_{\pm*})=b_\pm\Z.$$

Note that $a$ is the number of connected components of $H$. We use the convention $\Z/0=\Z$.

\begin{lemma}\label{kernels}
We have $\ker\iota_{\pm*}\cong\Z/b_\pm$. The group morphisms obtained by composing
\[
\begin{tikzcd}
\Z/b_\pm\cong\ker\iota_{\pm*}\ar{r}{} & \pi_1(M_\pm')\ar{r}{\zeta_{\pm*}^{-1}} & \pi_1(G/H)\ar{r}{} & \pi_0(H)\cong\Z/a
\end{tikzcd}
\]
are surjective. In particular, $b_\pm$ are divisible by $a$.
\end{lemma}

\begin{proof}
To ease the notation, we will omit the subscript $\pm$ in the proof of this lemma. We consider the following commutative diagram:
\[
\begin{tikzcd}
K\ar{r}{\mathrm{id}_K}\ar{d}{(j,-q)} & K\ar{d}{j}\\
G\times S^1\ar{r}{}\ar{d}{\rho} & G\times\D^2\ar{d}{}\\
M'\ar{r}{\iota} & M
\end{tikzcd}
\]
where both columns represent principal $K$-bundles and the middle row is the map obtained by the pullback of $G\times\D^2\longrightarrow M$ by $\iota$. The long exact sequence of homotopy induces a commutative diagram on homotopy groups where all the rows and columns are exact:
\[
\begin{tikzcd}
& & \pi_1(K)\ar{r}{\mathrm{id}_{\pi_1(K)}}\ar{d}{(j_*,-q_*)} & \pi_1(K)\ar{d}{j_*} & \\
0\ar{r}{} & \Z\ar{r}{}\ar[dashed]{d}{\kappa} & \pi_1(G)\times\Z\ar{d}{\rho_*}\ar{r}{} & \pi_1(G)\ar{r}{}\ar{d}{} & 0\\
0\ar{r}{} & \ker\iota_*\ar{r}{} & \pi_1(M')\ar{r}{\iota_*}\ar{d}{} & \pi_1(M)\ar{r}{}\ar{d}{} & 0\\
& & 0 & 0 &
\end{tikzcd}
\]
where $\kappa:\Z\longrightarrow\ker\iota_*$ is obtained by restricting $\rho_*$ to the subgroup generated by $(0,1)\in\pi_1(G)\times\Z$. By a standard diagram chase argument, the morphism $\kappa$ is well-defined and surjective. Moreover, note that $z\in\ker\kappa$ if and only if $(0,z)\in\ker\rho_*=\im(j_*,-q_*)$. Hence, $z\in\ker\kappa$ if and only if $z\in q_*(\ker j_*)=b\Z$, and consequently $\ker\iota_*\cong\Z/b$.

To prove the second part of the lemma, we consider the following commutative diagram:
\[
\begin{tikzcd}
H\ar{r}{i}\ar{d}{j\circ i} & K\ar{d}{(j,-q)}\\
G\ar[hook]{r}{}\ar{d}{} & G\times S^1\ar{d}{\rho}\\
G/H\ar{r}{\zeta} & M'
\end{tikzcd}
\]
where both columns are principal bundles. The long exact sequence of homotopy induces a commutative diagram on homotopy groups where all the columns are exact:
\[
\begin{tikzcd}
\pi_1(H)\ar{d}{j_*\circ i_*}\ar{r}{i_*} & \pi_1(K)\ar{d}{(j_*,-q_*)}\\
\pi_1(G)\ar{d}{}\ar{r}{} & \pi_1(G)\times\Z\ar{d}{\rho_*}\\
\pi_1(G/H)\ar{r}{\zeta_*}\ar{d}{\partial} & \pi_1(M')\\
\Z/a\ar{d}{} & \\
0 &
\end{tikzcd}
\]
Let $x$ denote a generator of $\Z/a$. By the surjectivity of $\partial:\pi_1(G/H)\longrightarrow\Z/a$, we can choose $y\in\pi_1(G/H)$ such that $\partial(y)=x$. We consider now $\zeta_*(y)\in\pi_1(M')$. By the surjectivity of $\rho_*:\pi_1(G)\times\Z\longrightarrow\pi_1(M')$, there exists $(\gamma,z)\in\pi_1(G)\times\Z$ satisfying $\rho_*(\gamma,z)=\zeta_*(y)$. Since $(\gamma,z)+(0,-z)=(\gamma,0)$, we have
$$(\zeta_*^{-1}\circ\rho_*)(\gamma,0)=y\cdot(\zeta_*^{-1}\circ\rho_*)(0,-z)\in\ker\partial,$$
which implies
$$(\partial\circ\zeta_*^{-1}\circ\rho_*)(0,z)=x.$$
Since $\rho_*(0,z)\in\ker\iota_*$, we conclude that
$$\partial\circ\zeta_{*|\ker\iota_*}^{-1}:\Z/b\longrightarrow\Z/a$$
is surjective, as claimed.
\end{proof}

\begin{lemma}\label{l=1 homogeneous space}
Consider the inclusion of compact connected Lie groups $K_\pm\longrightarrow G\times S^1$ given by $k_\pm\mapsto(j_\pm(k_\pm),q_\pm(k_\pm)^{-1})$. We have a $G$-equivariant diffeomorphism $(G\times S^1)/K_\pm\cong M_\pm'$. Moreover, the isotropy representations satisfy
$$\sigma'_\pm\cong\sigma_{G/K_\pm}\oplus 1_{K_\pm},$$
where $1_{K_\pm}$ denotes the trivial one-dimensional representation. In particular, $\sigma'_{\pm*}=\sigma_{G/K_\pm*}$.
\end{lemma}

\begin{proof}
As in the preceding lemmas, we omit the subscript $\pm$ to simplify the notation. Recall that $M'=G\times_K S^1$. The action of $K$ on $G$ is given by right multiplication via the inclusion $j:K\longrightarrow G$, while the action of $K$ on $S^1$ is transitive and by isometries, inducing a map $q:K\longrightarrow S^1$. Since $\Isom_0(S^1)\cong S^1$ as Lie groups and $K$ is connected, the action of $K$ on $S^1$ takes values in $\Isom_0(S^1)$. Identifying $\Isom_0(S^1)$ with $S^1$, we may therefore assume that $q$ is a Lie group morphism. Consequently, the map $K\longrightarrow G\times S^1$ given by $k\longmapsto(j(k),q(k)^{-1})$ is an injective Lie group homomorphism. Moreover, the map
$$
M'\longrightarrow(G\times S^1)/K,\qquad
[g,z]\longmapsto(g,z)K,
$$
is a $G$-equivariant diffeomorphism.

Since $K/H\cong S^1$, the subgroup $H$ is normal in $K$ and $K/H$ is abelian. Using an $\Ad(K)$-invariant inner product on $\mathfrak{k}$, we can choose an $\Ad(K)$-invariant decomposition
$$\mathfrak{k}=\mathfrak{h}\oplus\R u.$$
Since $K/H$ is abelian, the action of $\Ad(K)$ on the summand $\R u$ is trivial. The morphism $q:K\longrightarrow S^1$ induces a surjective Lie algebra homomorphism $dq:\mathfrak{k}\longrightarrow\R$, and we choose $u$ so that $dq(u)=1$.

Similarly, the inclusion $j:K\longrightarrow G$ induces a Lie algebra morphism $dj:\mathfrak{k}\longrightarrow\mathfrak g$. Since $G$ is compact, we may choose an $\Ad(K)$-invariant reductive decomposition
$$\mathfrak g=dj(\mathfrak k)\oplus\mathfrak m,$$
where $T_{eK}(G/K)\cong\mathfrak m$.

The Lie algebra of $G\times S^1$ is $\mathfrak g\oplus\R$, and the differential of the inclusion $K\longrightarrow G\times S^1$ is given by
$$v+tu\longmapsto(dj(v+tu),-t),$$
where $v\in\mathfrak h$ and $t\in\R$.

Using the identification $\mathfrak g\oplus\R=dj(\mathfrak h)\oplus\R dj(u)\oplus\mathfrak m\oplus\R$, define
$$\mathfrak m'=\{(t\,dj(u)+x,t):x\in\mathfrak m,\ t\in\R\}.$$
Then $\mathfrak g\oplus\R=d(j,-q)(\mathfrak{k})\oplus\mathfrak m'$, and this decomposition is $\Ad(K)$-invariant. Therefore, the isotropy representation of $(G\times S^1)/K$ is $\sigma'=\Ad(K)_{|\mathfrak m'}$.

Finally, the isomorphism
$$\mathfrak m'\cong\mathfrak m\oplus\R,\qquad
(t\,dj(u)+x,t)\longmapsto(x,t),$$
is $\Ad(K)$-equivariant. 
Since both homogeneous spaces are reductive, their isotropy representations are equivalent to the restrictions of the corresponding adjoint representations to the reductive complements. As the above isomorphism is $\Ad(K)$-equivariant and $\Ad(K)$ acts trivially on the factor $\R$, we obtain
$$
\sigma'\cong\sigma_{G/K}\oplus 1_K,
$$
where $1_K$ denotes the trivial one-dimensional representation of $K$. Consequently, after composing the induced maps on fundamental groups with the corresponding homomorphisms $\delta_n$, we obtain
$$\sigma'_*=\sigma_{G/K*},$$
as claimed.
\end{proof}

\section{Invariant spin structures and free group actions}\label{sec: preliminaries spin structures}

In this section, we recall the definition of an invariant spin structure and its main properties. We then establish the technical results needed to determine which manifolds admitting a cohomogeneity one action admit a spin structure or an invariant spin structure. The first main result is an equivariant version of the two-out-of-three lemma. We use this equivariant version to study the relation between invariant spin structures on connected manifolds $M$ admitting free actions of compact connected Lie groups $G$ and spin structures on the orbit spaces $M/G$. As an application, we recover the criteria for a homogeneous space $G/H$ to admit a spin structure or a $G$-invariant spin structure, where $G$ is a compact Lie group and $H$ is a closed connected subgroup. We conclude this section by studying spin structures and invariant spin structures on manifolds admitting free proper actions of discrete groups.

For $n\ge 2$, we consider the unique connected $2$--covering
\[
\lambda_n:\mathrm{Spin}(n)\longrightarrow \mathrm{SO}(n).
\]
If $n=1$, then $\mathrm{SO}(1)$ is trivial. In this case, $\mathrm{Spin}(1)\cong \Z/2$ and
$\lambda_1:\mathrm{Spin}(1)\to \mathrm{SO}(1)$ is the trivial group morphism.

Let $X$ be a connected manifold and let $p:E\to X$ be a real vector bundle of rank $n$,
endowed with a Euclidean bundle metric. Assume moreover that $E$ is orientable, i.e.\ $w_1(E)=0$.
Associated to $E$ we have the principal $\mathrm{SO}(n)$-bundle of oriented orthonormal frames
\[
q:FE\longrightarrow X.
\]
A spin structure on $q:FE\to X$ is a pair $(P,\Lambda)$ where $P\to X$ is a principal $\mathrm{Spin}(n)$-bundle and
$\Lambda:P\to FE$ is a bundle map covering $\mathrm{id}_X$ such that
\[
\Lambda(u\cdot a)=\Lambda(u)\cdot \lambda_n(a)\qquad\text{for all }u\in P,\ a\in \mathrm{Spin}(n),
\]
and such that $\Lambda$ restricts on each fibre to a $2$--fold covering.

Like orientability, the existence of spin structures on a principal $\mathrm{SO}(n)$-bundle is topologically
obstructed. Indeed, it is well known that $q:FE\to X$ admits a spin structure if and only if $w_2(E)=0$
(see \cite{lawson2016spin}). An orientable manifold is said to be spin if the oriented frame bundle associated to
its tangent bundle admits a spin structure.

We will use the following description of spin structures:

\begin{lemma}\label{Classification via FG}\cite{friedrich2000dirac}
Let $X$ be a connected manifold, let $E\to X$ be an orientable real vector bundle of rank
$n\ge 1$, and let $q:FE\to X$ be its oriented orthonormal frame bundle. Fix a fibre inclusion
$\iota:\mathrm{SO}(n)\hookrightarrow FE$.

Spin structures on $q:FE\to X$ are in one-to-one correspondence with group homomorphisms
\[
f:\pi_1(FE)\longrightarrow \Z/2
\]
such that the restriction to the fibre equals the homomorphism induced by the spin covering
$\lambda_n:\mathrm{Spin}(n)\to \mathrm{SO}(n)$, i.e.
\[
f\circ \iota_*=\delta_n,
\]
where $\delta_n:\pi_1(\mathrm{SO}(n))\to \Z/2$ is given as follows:
\[
\delta_n=
\begin{cases}
0 & (n=1),\\[2pt]
\text{reduction mod $2$: }\Z\cong \pi_1(\mathrm{SO}(2))\to \Z/2 & (n=2),\\[2pt]
\mathrm{id}_{\Z/2}:\ \Z/2\cong \pi_1(\mathrm{SO}(n))\to \Z/2 & (n\ge 3).
\end{cases}
\]
\end{lemma}

\begin{remark}
Equivalently, for a connected closed manifold $Y$ one has a natural identification
$H^{1}(Y,\Z/2)\cong \Hom(\pi_{1}(Y),\Z/2)$.
Thus spin structures correspond to cohomology classes $x\in H^{1}(FE,\Z/2)$ whose restriction to a fibre
$\mathrm{SO}(n)$ equals the class determined by the spin covering $\lambda_n$; explicitly,
\[
\iota^{*}x=
\begin{cases}
0 & (n=1),\\[2pt]
\text{the nonzero class in }H^{1}(\mathrm{SO}(2),\Z/2)\cong \Z/2 & (n=2),\\[2pt]
\text{the generator of }H^{1}(\mathrm{SO}(n),\Z/2)\cong \Z/2 & (n\ge 3).
\end{cases}
\]
\end{remark}

\subsection{Invariant spin structures and fundamental groups}
\label{sec:Invariant spin structures and fundamental groups}

Let $G$ be a compact Lie group. A $G$-vector bundle is a vector bundle $p:E\longrightarrow X$ such that:
\begin{itemize}
    \item[1.] we have actions of $G$ on $E$ and $X$ such that $p$ is $G$-equivariant;
    \item[2.] the action of $G$ on $E$ is by isometries. That is, given $x\in X$ and $g\in G$ the map
    $g:E_x\longrightarrow E_{gx}$ obtained by restricting the homeomorphism $g:E\longrightarrow E$ to the fibre
    over $x$ is a linear isometry of inner product spaces.
\end{itemize}

If $p:E\longrightarrow X$ is an orientable $G$-vector bundle of rank $n$ and the action of $G$ preserves the orientation then this action induces an action of $G$ on $FE$ which commutes with the action of $\SO(n)$ on $FE$.

\begin{example}
Let $G$ be a Lie group acting by orientation-preserving isometries on a connected oriented Riemannian manifold $M$. We can define an action of $G$ on $TM$ given by $g(x,v)=(gx,d_xg(v))$, which makes the tangent bundle $TM\longrightarrow M$ into a $G$-vector bundle. The action of $G$ on the oriented orthonormal frame bundle $FM$ is the following: given $(x,\{v_1,...,v_n\})\in FM$, we set
\[
g(x,\{v_1,...,v_n\})=(gx,\{d_xg(v_1),...,d_xg(v_n)\}).
\]
\end{example}

\begin{definition}[$G$-invariant spin structure]
Let $E\to X$ be a spin $G$-vector bundle and let $(P,\Lambda)$ be a spin structure on $q:FE\to X$.
We say that $(P,\Lambda)$ is \emph{$G$-invariant} if the induced action of $G$ on $FE$ lifts to an action of $G$ on
$P$ commuting with the right $\mathrm{Spin}(n)$-action, such that $\Lambda$ is $G$-equivariant.
\end{definition}

We want to determine which spin structures on an orientable $G$-vector bundle are $G$-invariant when $G$ is a compact connected Lie group. Thus, from now we will assume that $G$ is a compact connected Lie group until the end of the section. To formulate the criterion, let $x\in FE$. Define the evaluation map $ev_x:G\longrightarrow FE$ by $ev_x(g)=g x$.
The induced map
\[
(ev_x)_*:\pi_1(G,e)\longrightarrow \pi_1(FE,x)
\]
does not depend on the choice of $x$ (see \cite[Theorem 2.4.2]{lee2010seifert}). This is because the image of $(ev_x)_*$ is a central subgroup of $\pi_1(FE,x)$ and changing the basepoint change the fundamental group via a conjugation. Thus, we may suppress basepoints and denote the induced map simply by
\[
ev_*:\pi_1(G)\longrightarrow \pi_1(FE).
\]

\begin{remark}
If the action of $G$ on $FE$ is free, then $FE\longrightarrow FE/G$ is a principal $G$-bundle and $ev_x$ identifies
$G$ with the fibre (orbit) $G\cdot x\subset FE$ through $x$. For example, this holds for the induced action on $FM$
of any effective action of $G$ by orientation-preserving
isometries on a connected Riemannian manifold $M$.
\end{remark}

The next lemma is a direct consequence of \cite[Theorem 2.5.1]{lee2010seifert}.

\begin{lemma0}\label{invariant spin structure criteria}
Let $G$ be a compact connected Lie group. A spin structure $f:\pi_1(FE)\longrightarrow\Z/2$ is $G$-invariant if and only if
$f\circ ev_*:\pi_1(G)\longrightarrow \Z/2$ is trivial.
\end{lemma0}

\begin{remark}
Let $p:E\longrightarrow X$ be a spin $G$-vector bundle. Assume that a spin structure
$f:\pi_1(FE)\longrightarrow \Z/2$ is not $G$-invariant and let $\Tilde{G}\longrightarrow G$ be the connected
$2$-fold covering of Lie groups corresponding to the index-$2$ subgroup
$\ker(f\circ ev_*)\leq \pi_1(G,e)$.
Then we can define actions of $\Tilde{G}$ on $E$ and $X$ by composing $\Tilde{G}\to G$
with the given $G$-actions, making $p:E\longrightarrow X$ into a $\Tilde{G}$-vector bundle. The spin structure
$f:\pi_1(FE)\longrightarrow \Z/2$ is $\Tilde{G}$-invariant.
\end{remark}
\begin{lemma}\label{2 out of 3 homotopy}
Let $G$ be a compact connected Lie group and $E_i\longrightarrow X$ with $i=1,2$ be $G$-vector bundles, and let $E=E_1\oplus E_2\longrightarrow X$ denote the Whitney sum of $E_1$ and $E_2$, which is a $G$-vector bundle with the diagonal action. Assume that all three vector bundles
$E_1$, $E_2$, and $E$ are spin.

Let $f_i:\pi_1(FE_i)\to \Z/2$ be spin structures on $E_i$ $(i=1,2)$, and let
$f:\pi_1(FE)\to \Z/2$ be the spin structure on $E$ obtained from $f_1$ and $f_2$
as in the proof below. Then
\[
f\circ ev_* \;=\; (f_1\circ ev_{1*})+(f_2\circ ev_{2*})
\qquad\text{in }\Hom(\pi_1(G),\Z/2).
\]
If any two of them are not $G$-invariant but are invariant under the same connected $2$-fold covering $\Tilde{G}\longrightarrow G$ then the third one is $G$-invariant. In particular, if any two of the three spin structures $f_1,f_2,f$ are $G$-invariant, then so is the third. 

\end{lemma}
\begin{proof}
    Let $q_i:FE_i\longrightarrow X$ be the principal
$\SO(n_i)$-bundles associated to the oriented
vector bundles $E_i$, and write $n=n_1+n_2$.
Then $FE_1\times_X FE_2\longrightarrow X$
is a principal $\SO(n_1)\times\SO(n_2)$-bundle,
and the oriented frame bundle of $E$ is the
associated bundle
\[
q:(FE_1\times_X FE_2)
\times_{\SO(n_1)\times\SO(n_2)}
\SO(n)\longrightarrow X,
\]
constructed using the block-diagonal inclusion
\[
j:\SO(n_1)\times\SO(n_2)\longrightarrow\SO(n),
\qquad
j(A,B)=
\begin{pmatrix}
A&0\\
0&B
\end{pmatrix}.
\]
Let $p_i:FE_1\times_X FE_2\to FE_i$
denote the natural projections, and let
\[
\pi:(FE_1\times_X FE_2)\times\SO(n)
\longrightarrow FE
\]
be the quotient map.
Recall from \cref{Classification via FG}
that the homomorphisms
$\delta_k:\pi_1(\SO(k))\to\Z/2$
are given by
\[
\delta_k=
\begin{cases}
0 & k=1,\\
\text{reduction mod }2 & k=2,\\
\mathrm{id}_{\Z/2} & k\ge3.
\end{cases}
\]
For the block-diagonal inclusion $j$,
they satisfy
\[
\delta_n\circ j_*=
\delta_{n_1}+\delta_{n_2}.
\] We first construct $f$ from $f_1$ and $f_2$.
Consider the group homomorphism
\[
\overline f:
\pi_1(FE_1\times_X FE_2)
\times\pi_1(\SO(n))
\longrightarrow\Z/2
\]
defined by
\[
\overline f(a,b)=
f_1(p_{1*}(a))+f_2(p_{2*}(a))
+\delta_n(b).
\]
The map $\pi$ is a principal
$\SO(n_1)\times\SO(n_2)$-bundle.
The identity
$\delta_n\circ j_*=
\delta_{n_1}+\delta_{n_2}$
and the fibre conditions for $f_1$ and $f_2$
show that $\overline f$ vanishes on the image
of the fundamental group of its fibre.
By the homotopy exact sequence,
$\overline f$ therefore descends to a
homomorphism
\[
f:\pi_1(FE)\longrightarrow\Z/2.
\]
Moreover, since $\pi$ restricts on the
$\SO(n)$-factor to the usual fibre inclusion
in $FE$, we have
\[
f\circ i_*=\delta_n.
\]
Hence $f$ is a spin structure on $E$.

The action of $G$ on $FE_1$ and $FE_2$
induces an action on
$(FE_1\times_X FE_2)\times\SO(n)$
given by
\[
g(e_1,e_2,A)=(ge_1,ge_2,A).
\]
This action descends to the action of $G$
on $FE$. Consequently, for any
$\alpha\in\pi_1(G)$ we have
\[
\begin{aligned}
(f\circ ev_*)(\alpha)
&=(\overline f\circ\overline{ev}_*)(\alpha)\\
&=f_1(ev_{1*}(\alpha))
  +f_2(ev_{2*}(\alpha)),
\end{aligned}
\]
since the $\SO(n)$-coordinate is fixed
by $G$. Thus
\[
f\circ ev_*=
(f_1\circ ev_{1*})+(f_2\circ ev_{2*}).
\]
Conversely, start with spin structures
$f:\pi_1(FE)\to\Z/2$ and
$f_1:\pi_1(FE_1)\to\Z/2$.
Consider the homomorphism
\[
\overline f:
\pi_1(FE_1\times_X FE_2)
\times\pi_1(\SO(n))
\longrightarrow\Z/2
\]
defined by
\[
\overline f(a,b)=
f_1(p_{1*}(a))
+f(\pi_*(a,b))
+\delta_n(b).
\]

Consider the natural projection
\[
r:(FE_1\times_X FE_2)\times\SO(n)
\longrightarrow FE_2.
\]
This is a fibre bundle with fibre
$\SO(n_1)\times\SO(n)$.
Using the fibre conditions for $f$ and $f_1$,
a direct computation shows that
$\overline f$ vanishes on the image of
the fundamental group of this fibre.
Indeed, on
$\pi_1(\SO(n_1))\times\pi_1(\SO(n))$
its value is
\[
\delta_{n_1}(a_1)
+\delta_n(j_*(a_1,0)+b)
+\delta_n(b)=0.
\]
Since the fibre is connected, the homotopy
exact sequence shows that $r_*$ is
surjective and that its kernel is the
image of the fundamental group of the
fibre. Thus $\overline f$ descends to a
homomorphism
\[
f_2:\pi_1(FE_2)\longrightarrow\Z/2.
\]

The fibre condition for $f_2$ follows
from that of $f$ and the identity
$\delta_n\circ j_*=
\delta_{n_1}+\delta_{n_2}$.
Consequently, $f_2$ is a spin structure
on $E_2$.

Moreover, the same evaluation map
computation gives
\[
f_2\circ ev_{2*}
=(f_1\circ ev_{1*})+(f\circ ev_*).
\]
The remaining case, starting from
$f$ and $f_2$, is analogous.

Finally, the invariance conclusions follow
from these identities and
\cref{invariant spin structure criteria}.
Indeed, if any two of the three homomorphisms
$f\circ ev_*$, $f_1\circ ev_{1*}$, and
$f_2\circ ev_{2*}$ vanish, then so does
the third.
If two spin structures are not $G$-invariant
but are invariant under the same connected
$2$-fold covering $\Tilde G\to G$,
then their obstruction characters coincide
with the nontrivial character corresponding
to this covering. Their sum therefore
vanishes, so the third spin structure
is $G$-invariant.
\end{proof}

\subsection{Invariant spin structures and free actions of connected groups on manifolds}
\label{sec:Invariant spin structures and free actions on manifolds}

Recall that if a compact Lie group $G$ acts freely and by orientation-preserving isometries on a
connected manifold $M$, then the orbit space $M/G$ is a connected smooth manifold of dimension
$\dim M-\dim G$ and the orbit map $\pi:M\longrightarrow M/G$ is a principal $G$-bundle.

If $G$ is connected, then the action of $G$ on an oriented $M$ preserves the orientation, hence $M/G$ is
orientable. Let $\tau_G:=\ker(d\pi)\subset TM$ be the bundle along the fibres. Then there is a short exact
sequence of $G$-vector bundles
\[
0\longrightarrow \tau_G \longrightarrow TM \xrightarrow{d\pi} \pi^*T(M/G)\longrightarrow 0.
\]
Choosing a $G$-invariant connection (for instance the orthogonal complement of $\tau_G$ since the action is
isometric) yields a $G$-equivariant splitting and hence a decomposition
\[
TM\cong \pi^*T(M/G)\oplus \tau_G.
\]
Moreover, $\tau_G$ is (non-canonically) isomorphic to the trivial bundle $M\times\mathfrak{g}\to M$ via
fundamental vector fields, and the $G$-action on $M\times\mathfrak{g}$ is given by
\[
(x,v)g=(xg,\Ad_{g^{-1}}(v)).
\]
For more information on the topic we refer to any standard reference on principal bundles, such as \cite{bredon1972introduction} or \cite[Appendix B]{friedrich2000dirac}.

The goal of this section is to study the relation between spin structures on $M/G$ and invariant spin
structures on $M$.

\begin{lemma}\label{invariant spin structure pullback}
The tangent bundle $T(M/G)\longrightarrow M/G$ admits a spin structure if and only if the pullback bundle
$\pi^*T(M/G)\longrightarrow M$ admits a $G$-invariant spin structure.
\end{lemma}

\begin{proof}
Let $n:=\dim(M/G)$. Write $F(M/G)\to M/G$ for the oriented orthonormal frame bundle of $T(M/G)$.
The pullback bundle $\pi^*T(M/G)\to M$ has oriented frame bundle canonically identified with the pullback
principal $\SO(n)$-bundle $\pi^*F(M/G)\to M$.

The projection
\[
\rho:\pi^*F(M/G)\longrightarrow F(M/G),\qquad (x,u)\longmapsto u,
\]
is a principal $G$-bundle, where $G$ acts on $\pi^*F(M/G)$ by $(x,u)\cdot g=(xg,u)$ (so $\rho$ is $G$-invariant).
In particular, for any $x\in \pi^*F(M/G)$ the composition $\rho\circ ev_x:G\to F(M/G)$ is constant, hence
$\rho_*\circ ev_*=0$ on $\pi_1$.

\smallskip

Assume first that $T(M/G)$ admits a spin structure, encoded, as given in Lemma \ref{Classification via FG}, by a
homomorphism $f:\pi_1(F(M/G))\to \Z/2$ satisfying the fibre condition $f\circ \iota_*=\delta_n$.
Then
\[
f':=f\circ \rho_*:\pi_1(\pi^*F(M/G))\longrightarrow \Z/2
\]
is a spin structure on $\pi^*T(M/G)$, since $\rho$ restricts to the identity on each $\SO(n)$-fibre, so the same
fibre condition holds (with the same $\delta_n$, including the cases $n=1,2$).
Moreover,
\[
f'\circ ev_* = f\circ (\rho_*\circ ev_*) = 0,
\]
so $f'$ is $G$-invariant by Lemma~\ref{invariant spin structure criteria}.

\smallskip

Conversely, let $f':\pi_1(\pi^*F(M/G))\to \Z/2$ be a $G$-invariant spin structure on $\pi^*T(M/G)$.
By Lemma~\ref{invariant spin structure criteria} we have $f'\circ ev_*=0$.
Since $\rho:\pi^*F(M/G)\to F(M/G)$ is a principal $G$-bundle, the homotopy long exact sequence gives
\[
\ker(\rho_*)=\im(ev_*).
\]
Hence $f'$ vanishes on $\ker(\rho_*)$ and therefore descends uniquely to a homomorphism
$\bar f:\pi_1(F(M/G))\to \Z/2$ such that $f'=\bar f\circ \rho_*$.
Finally, the fibre condition for $\bar f$ follows from that of $f'$. If $\iota:\SO(n)\hookrightarrow F(M/G)$ is a
fibre inclusion and $\iota'$ is the corresponding fibre inclusion in $\pi^*F(M/G)$ with $\rho\circ \iota'=\iota$,
then
\[
\bar f\circ \iota_* \;=\; f'\circ \iota'_* \;=\; \delta_n.
\]

Thus $\bar f$ is a spin structure on $T(M/G)$.
\end{proof}
\begin{corollary}\label{spin structures orbit spaces}
The manifold $M/G$ admits a spin structure if and only if there exist spin structures on
$TM\longrightarrow M$ and $\tau_G$, denoted respectively by $f:\pi_1(FM)\longrightarrow\Z/2$ and
$f_\tau:\pi_1(M)\times\pi_1(\SO(\dim G))
\longrightarrow \Z/2$
(using the identification
$F(\tau_G)\cong M\times\SO(\mathfrak g)$), satisfying
\[
f\circ ev_* \;+\; f_\tau\circ\bigl(ev^{G\acts M}_*,\Ad_*\bigr)=0
\qquad\text{in }\Hom(\pi_1(G),\Z/2).
\]
Here $ev_*$ denotes the evaluation map for the $G$-action on $FM$, and
\[
\Ad_*:\pi_1(G)\longrightarrow \pi_1(SO(\dim G))
\] is the map induced by the adjoint representation 
$\Ad:G\to \SO(\mathfrak g)$ 
(with respect to an $\Ad$-invariant inner product).

In particular, if $G$ is either simply connected or abelian, then $M/G$ is spin if and only if $M$ has a
$G$-invariant spin structure.
\end{corollary}

\begin{proof}
The first part of the corollary is obtained by combining \cref{2 out of 3 homotopy} and
\cref{invariant spin structure pullback}. Concretely, one applies \cref{2 out of 3 homotopy} to the splitting
$TM\cong \pi^*T(M/G)\oplus\tau_G$ and uses \cref{invariant spin structure pullback} to relate spin structures on
$T(M/G)$ to $G$-invariant spin structures on $\pi^*T(M/G)$; the evaluation map for $\tau_G$ corresponds to
$(ev^{G\acts M}_*,\Ad_*)$ under the identification $F(\tau_G)\cong M\times \SO(\mathfrak g)$.

For the second part of the statement, note that if $G$ is either simply connected or abelian, then
$\Ad_*:\pi_1(G)\longrightarrow\pi_1(\SO(\dim G))$ is trivial. Thus, the spin structure
\[
f_\tau:\pi_1(M)\times\pi_1(\SO(\dim G))
\longrightarrow\Z/2,\qquad
f_\tau(a,b)=\delta_{\dim G}(b),
\]
is $G$-invariant.
Consequently, the second part of the corollary follows from \cref{2 out of 3 homotopy} and
\cref{invariant spin structure pullback}.
\end{proof}

Assume now that $M$ is an orientable manifold and that $G$ and $H$ are connected compact Lie groups. Assume that we have a free right action of $G$ on $M$ and a left action of $H$ on $M$ which commute with the action of $G$, that is $h(xg)=(hx)g$ for all $x\in M$, $g\in G$ and $h\in H$. Then $H$ acts on $M/G$ and the projection map $\pi:M\longrightarrow M/G$ is $H$-equivariant. Then the $G$-vector bundles $\pi^*T(M/G)$ and $\tau_G$ are also $H$-vector bundles.

\begin{lemma}\label{invariant spin structure pullback invariant}
The tangent bundle $T(M/G)\longrightarrow M/G$ admits a $H$-invariant spin structure if and only if
$\pi^*T(M/G)\longrightarrow M$ admits a spin structure which is $G$-invariant and $H$-invariant.    
\end{lemma}

\begin{proof}
We use the same notation as in \cref{invariant spin structure pullback}. By hypothesis, the map
$$\rho:\pi^*F(M/G)\longrightarrow F(M/G)$$ is $H$-equivariant, and hence
\[
\rho_* \circ ev_*^{H\acts \pi^*F(M/G)} = ev_*^{H\acts F(M/G)}.
\]

If $f:\pi_1(F(M/G))\longrightarrow \Z/2$ is a $H$-invariant spin structure on $M/G$, then
$f\circ ev_*^{H\acts F(M/G)}=0$. Consequently, the $G$-invariant spin structure
$f\circ \rho_*:\pi_1(\pi^*F(M/G))\longrightarrow \Z/2$ satisfies
\[
(f\circ \rho_*) \circ ev_*^{H\acts \pi^*F(M/G)}
= f\circ (\rho_* \circ ev_*^{H\acts \pi^*F(M/G)})
= f\circ ev_*^{H\acts F(M/G)}=0,
\]
which implies that $f\circ\rho_*$ is also $H$-invariant.

Conversely, if $f:\pi_1(\pi^*F(M/G))\longrightarrow \Z/2$ is a $G$-invariant and $H$-invariant spin structure on
$\pi^*F(M/G)\longrightarrow M$, then the induced spin structure
$\overline{f}:\pi_1(F(M/G))\longrightarrow \Z/2$ on $M/G$ satisfies that $\overline{f}\circ\rho_*=f$.
Therefore
\[
\overline{f}\circ ev_*^{H\acts F(M/G)}
=\overline{f}\circ(\rho_* \circ ev_*^{H\acts \pi^*F(M/G)})
=(\overline{f}\circ\rho_*)\circ ev_*^{H\acts \pi^*F(M/G)}=0
\]
and the spin structure $\overline{f}$ is $H$-invariant.
\end{proof}

The following corollary is the invariant spin structure counterpart to \cref{spin structures orbit spaces}.

\begin{corollary}\label{invariant spin structures orbit spaces}
The manifold $M/G$ admits a $H$-invariant spin structure if and only if there exist spin structures on
$p:TM\longrightarrow M$ and $\tau_G$, encoded by
$f:\pi_1(FM)\longrightarrow \Z/2$ and
$f_\tau:\pi_1(M)\times\pi_1(\SO(\dim G))
\longrightarrow \Z/2$
(using the identification
$F(\tau_G)\cong M\times\SO(\mathfrak g)$),
satisfying that
\[
f\circ ev_*^{G\acts FM}+f_\tau\circ(ev^{G\acts M}_*,\Ad_*)=0
\qquad\text{and}\qquad
f\circ ev_*^{H\acts FM}=f_\tau\circ ev_*^{H\acts F(\tau_G)}.
\]
\end{corollary}

Using \cref{spin structures orbit spaces}, \cref{invariant spin structure criteria},
\cref{invariant spin structure pullback invariant} and \cref{invariant spin structures orbit spaces}, we can
reprove the following criterion for the existence of (invariant) spin structures on homogeneous spaces.
Let $G/H$ be an orientable Riemannian homogeneous space, where $G$ and $H$ are connected compact Lie groups, and
let $\mathfrak{g}$ and $\mathfrak{h}$ denote their respective Lie algebras. Then $G/H$ is reductive, therefore
there exists a $\Ad(H)$-invariant subspace $\mathfrak{m}\subseteq \mathfrak{g}$ such that
$\mathfrak{g}=\mathfrak{h}\oplus\mathfrak{m}$. Recall that $\mathfrak{m}\cong T_{eH}G/H$ and that the isotropy
representation $\sigma:H\longrightarrow \SO(T_{eH}G/H)$ defined as
\[
\sigma(h)=d_{eH}h:T_{eH}G/H\longrightarrow T_{eH}G/H
\]
satisfies that $\sigma(h)=\Ad(h)_{|\mathfrak{m}}$ for all $h\in H$.

\begin{proposition}\label{spin structures homogeneous spaces}
Let $G$ be a compact connected Lie group, let $H$ be a closed connected subgroup and let $G/H$ be a closed
Riemannian $n$-dimensional homogeneous space. Let $i:H\longrightarrow G$ denote the inclusion, let
$\pi:G\longrightarrow G/H$ denote the projection map and let $\sigma:H\longrightarrow \SO(n)$ denote the isotropy
representation. We write
\[
\sigma_*:=\delta_{n}\circ (\sigma)_*:\pi_1(H)\to \Z/2,
\]
where $\delta_n$ is as in Lemma \ref{Classification via FG}. Then:
\begin{itemize}
	\item[(1)] If $G/H$ is spin then $\sigma_*(\ker i_*)=0$. Conversely, if $\sigma_*(\ker i_*)=0$ and
	$\pi^*:\Ext^1_\Z(\pi_1(G/H),\Z/2)\longrightarrow \Ext^1_\Z(\pi_1(G),\Z/2)$ is injective, then $G/H$ is spin.
	\item[(2)] $G/H$ has a $G$-invariant spin structure if and only if $\sigma_*=0$.
\end{itemize}
\end{proposition}

\begin{proof}
We start by proving item (1). Since Lie groups are parallelizable, we have $TG\cong G\times \mathfrak{g}$.
The right action of $H$ on $TG$ satisfies
\[
(g,x)h=(gh,\Ad_G(h^{-1})(x)).
\]
Similarly, $\tau_H\cong G\times\mathfrak{h}$ and the right action of $H$ on $\tau_H$ satisfies
\[
(g,x)h=(gh,\Ad_H(h^{-1})(x)).
\]
Using the $\Ad_H$-invariant decomposition $\mathfrak{g}=\mathfrak{h}\oplus\mathfrak{m}$, we have
\[
\Ad_G(h)(x)=\Ad_H(h)(x_\mathfrak{h})+\sigma(h)(x_\mathfrak{m}),
\]
where $x_\mathfrak{h}\in\mathfrak{h}$, $x_\mathfrak{m}\in\mathfrak{m}$ and $x=x_\mathfrak{h}+x_\mathfrak{m}$.
Consequently, the evaluation map
\[
ev_*^{H\acts F(G)}:
\pi_1(H)\longrightarrow
\pi_1(G)\times\pi_1(\SO(\dim G))
\]
satisfies
\[
ev_*^{H\acts F(G)}
=(i_*,(\Ad_G|_H)_*).
\]
Moreover, the block-diagonal decomposition gives
\[
\delta_{\dim G}\circ(\Ad_G|_H)_*
=\Ad_{H*}+\sigma_*,
\]
where
$\Ad_{H*}:=\delta_{\dim H}\circ(\Ad_H)_*$.

By \cref{spin structures orbit spaces}, $G/H$ admits a spin structure if and only if we have spin structures $f:\pi_1(G)\times\pi_1(\SO(\dim G))\longrightarrow\Z/2$ and
$f_\tau:\pi_1(G)\times\pi_1(\SO(\dim H))\longrightarrow\Z/2$ satisfying that
\[
f\circ(i_*,(\Ad_G|_H)_*)=f_\tau\circ(i_*,(\Ad_H)_*).
\]
Let $\iota$ and $\iota_\tau$ denote the inclusions of $\pi_1(G)$ into the respective first factors, and let
$\phi=f\circ\iota$ and $\phi_\tau=f_\tau\circ\iota_\tau$.

\smallskip
\noindent
By the fibre condition in Lemma \ref{Classification via FG} (encoded via the corresponding $\delta$-maps),
for every $\alpha\in\pi_1(H)$ we have
\[
f\bigl(i_*(\alpha),\,(\Ad_G|_H)_*(\alpha)\bigr)
=\phi(i_*(\alpha))+\Ad_{H*}(\alpha)+\sigma_*(\alpha),
\]
and
\[
f_\tau\bigl(i_*(\alpha),\,(\Ad_H)_*(\alpha)\bigr)
=\phi_\tau(i_*(\alpha))+\Ad_{H*}(\alpha).
\]
Therefore $f\circ(i_*,(\Ad_G|_H)_*)=f_\tau\circ(i_*,(\Ad_H)_*)$ if and only if
\[
\phi\circ i_*+\Ad_{H*}+\sigma_*=\phi_\tau\circ i_*+\Ad_{H*}.
\]
By defining $\psi=\phi+\phi_\tau$, we have that $G/H$ has a spin structure if and only if there exists
$\psi:\pi_1(G)\longrightarrow \Z/2$ satisfying that
\[
\psi\circ i_*=\sigma_*.
\]

Assume that $\psi$ exists and that $\gamma\in \ker i_*$. Then $(\psi\circ i_*)(\gamma)=\sigma_*(\gamma)=0$ and
therefore $\sigma_*(\ker i_*)=0$.

Conversely, assume that $\sigma_*(\ker i_*)=0$. Then we can consider group morphisms
$\overline{i_*}:\pi_1(H)/\ker i_*\longrightarrow\pi_1(G)$ and
$\overline{\sigma_*}:\pi_1(H)/\ker i_*\longrightarrow\Z/2$.
Note that $\overline{i_*}$ is injective. The map $\psi$ exists if and only if there exists a map
$\overline{\psi}:\pi_1(G)\longrightarrow \Z/2$ satisfying that
$\overline{\psi}\circ \overline{i_*}=\overline{\sigma_*}$.

We have a short exact sequence of abelian groups
\[
\begin{tikzcd}
0 \ar{r}{} & \pi_1(H)/\ker i_* \ar{r}{\overline{i_*}} & \pi_1(G) \ar{r}{} & \pi_1(G/H) \ar{r}{} & 0.
\end{tikzcd}
\]
Applying the functor $\Hom(\cdot,\Z/2)$ yields the long exact sequence
\[
\begin{tikzcd}[column sep=small, row sep=large]
0 \ar[r]
& \Hom(\pi_1(G/H),\Z/2) \ar[r]
& \Hom(\pi_1(G),\Z/2)
  \ar[r, "\overline{i}^{\,*}"]
& \Hom(\pi_1(H)/\ker i_*,\Z/2)
  \ar[dll, rounded corners, to path={
    -- ([yshift=-2ex]\tikztostart.south)
    -| (\tikztotarget.north) \tikztonodes
  }]
\\
& \Ext^1_\Z(\pi_1(G/H),\Z/2)
  \ar[r, "\pi^*"]
& \Ext^1_\Z(\pi_1(G),\Z/2) \ar[r]
& \cdots
\end{tikzcd}
\]

Since $\pi^*:\Ext^1_\Z(\pi_1(G/H),\Z/2)\to \Ext^1_\Z(\pi_1(G),\Z/2)$ is injective, exactness implies that
$\overline{i_*}^{\,*}$ is surjective. Hence $\overline{\psi}$ exists, and therefore $G/H$ is spin.

\smallskip

We now prove item (2). The action of $G$ on $G\times \mathfrak{g}$ is given by $g'(g,x)=(g'g,x)$. Similarly, the
action of $G$ on $G\times \mathfrak{h}$ is given by $g'(g,x)=(g'g,x)$. Thus,
$ev_*^{G\acts F(G)}=(\id_{\pi_1(G)},0)$ and $ev_*^{G\acts F(\tau_H)}=(\id_{\pi_1(G)},0)$.

By \cref{invariant spin structures orbit spaces}, $G/H$ has a $G$-invariant spin structure if and only if there
exist spin structures
$f:\pi_1(G)\times\pi_1(\SO(\dim G))\longrightarrow\Z/2$ and
$f_\tau:\pi_1(G)\times\pi_1(\SO(\dim H))\longrightarrow\Z/2$
satisfying that
\[
f\circ(i_*,(\Ad_G|_H)_*)=f_\tau\circ(i_*,(\Ad_H)_*)
\qquad\text{and}\qquad
f\circ(\id_{\pi_1(G)},0)=f_\tau\circ(\id_{\pi_1(G)},0).
\]

The first condition is equivalent (as above) to $\phi\circ i_*+\sigma_*=\phi_{\tau}\circ i_*$, and the second
condition is equivalent to $\phi=\phi_\tau$. Consequently, $G/H$ has a $G$-invariant spin structure if and only if
$\sigma_*=0$. 
\end{proof}

\begin{remark}
The existence of spin structures on homogeneous spaces has been studied in some particular cases (see for example \cite{cahen1988spin,cahen1993spin,gadea2018homogeneous,alekseevsky2019spin}), although it is hard to find general results in the literature. In particular, it is shown in \cite{cahen1993spin} that if $G$ is simply connected, then $G/H$ has a spin structure if and only if there exists a lift of the isotropy representation $\Tilde{\sigma}:H\longrightarrow \Spin(n)$ satisfying $\lambda_n\circ \Tilde{\sigma}=\sigma$. 

The second part of \cref{spin structures homogeneous spaces} was proved in \cite{daura2022g}. Note that it is shown in \cite{daura2022g} that if $G/H$ has a $G$-invariant spin structure, then it is unique.
\end{remark}

\begin{remark}
The hypothesis that $\pi^*:\Ext^1_\Z(\pi_1(G/H),\Z/2)\longrightarrow \Ext^1_\Z(\pi_1(G),\Z/2)$ is injective in
item~(1) is necessary. For example, consider the homogeneous space
\[
M=\operatorname{PSU}(4)/\SO(3),
\]
where the inclusion $\SO(3)\longrightarrow \operatorname{PSU}(4)$ is induced by the inclusion
$\SU(2)\longrightarrow \SU(4)$ defined by $A\mapsto \diag(A,A)$.
The inclusion map induces the injective morphism $i_*:\Z/2\longrightarrow \Z/4$, hence
$\pi_1(G/H)\cong \Z/2$ and the induced map
\[
\pi^*:\Ext^1_\Z(\pi_1(G/H),\Z/2)\longrightarrow \Ext^1_\Z(\pi_1(G),\Z/2)
\]
is trivial. In particular, the condition $\sigma_*(\ker i_*)=0$ is satisfied (since $i_*$ is injective).

We now compute the isotropy representation of $\operatorname{PSU}(4)/\SO(3)$. Consider the inclusion
$j:\SO(3)\longrightarrow \SO(4)\cong (\SU(2)\times \SU(2))/(\Z/2)$. By \cref{hf isotoropy representation}, we have
\[
\sigma=\sigma_{\SO(4)/\SO(3)}+\sigma_{\operatorname{PSU}(4)/\SO(4)}\circ j.
\]
Note that $j_*:\pi_1(\SO(3))\longrightarrow\pi_1(\SO(4))$ is an isomorphism and that
$\SO(4)/\SO(3)\cong S^3$, hence the isotropy representation satisfies
$\sigma_{\SO(4)/\SO(3)*}= \id_{\Z/2}$ (see \cite{daura2022g}).

Thus it remains to compute
\[
\sigma_{\operatorname{PSU}(4)/\SO(4)*}:\pi_1(\SO(4))\longrightarrow \pi_1(\SO(9))
\]
(or equivalently its composition with $\delta_9$, which is the identity since $9\ge 3$).
Using the identification $\SO(4)\cong (\SU(2)\times \SU(2))/\Z/2$, the loop
$\gamma:[0,1]\longrightarrow (\SU(2)\times \SU(2))/\Z/2$ given by
\[
\gamma(t)=\diag(e^{i\pi t}, e^{-i\pi t} ,e^{i\pi t}, e^{-i\pi t})
\]
represents the generator of $\pi_1(\SO(4))$.

The homogeneous space $\operatorname{PSU}(4)/\SO(4)$ is reductive and therefore there exists an
$\Ad(\SO(4))$-invariant subspace $\mathfrak{m}$ satisfying
$\mathfrak{su}(4)=(\mathfrak{su}(2)\oplus \mathfrak{su}(2))\oplus \mathfrak{m}$.
We can decompose $\mathfrak{m}=V\oplus \R$, where $V$ is the subspace of matrices of the form
\[
V=\left\{\left(\begin{array}{@{}c|c@{}}
	0 &   A\\ \hline
	-A^{\dagger}  & 0
\end{array}\right):A\in M_{2\times2}(\C)\right\},
\]
and the second factor is generated by the matrix $\diag(i,i,-i,-i)$.
It is straightforward to see that this decomposition is $\Ad(\SO(4))$-invariant. Thus,
\[
\sigma_{\operatorname{PSU}(4)/\SO(4)*}[\gamma]=[\Ad(\gamma)_{|\mathfrak{m}}]=[\Ad(\gamma)_{|V}].
\]
Write $v\in V$ as
\[
v=\left(\begin{array}{@{}c|c@{}}
0 &   \begin{matrix}
a & b\\
c & d
\end{matrix}
\\ \hline
\begin{matrix}
-\overline{a} & -\overline{c}\\
-\overline{b} & -\overline{d}
\end{matrix} & 0
\end{array}\right),
\qquad a,b,c,d\in\C.
\]
Then
\[
\Ad(\gamma(t))(v)=\left(\begin{array}{@{}c|c@{}}
0 &   \begin{matrix}
a & e^{2\pi it }b\\
e^{-2\pi it }c & d
\end{matrix}
\\ \hline
\begin{matrix}
-\overline{a} & -e^{2\pi it }\overline{c}\\
-e^{-2\pi it }\overline{b} & -\overline{d}
\end{matrix} & 0
\end{array}\right).
\]
Thus $\Ad(\gamma)_{|V}$ is a product of two $2\pi$-rotations in orthogonal $2$-planes, and hence
$\sigma_{\operatorname{PSU}(4)/\SO(4)*}[\gamma]=0\in \pi_1(\SO(9))\cong \Z/2$.
Combining all these observations, we obtain $\sigma_*=\id_{\Z/2}$.
Since $i_*$ and $\sigma_*$ send generators to $2\in\Z/4$ and $1\in\Z/2$, respectively,
\[
\pi_1(FM)\cong (\pi_1(G)\times\Z/2)/\im(i_*,\sigma_*)\cong \Z/4,
\]
and the short exact sequence
\[
\begin{tikzcd}
1\ar{r}{} & \pi_1(\SO(12))\cong\Z/2\ar{r}{} & \pi_1(FM)\cong \Z/4\ar{r}{} &
\pi_1(M)\cong \Z/2\ar{r}{} & 1
\end{tikzcd}
\]
does not split. This implies that $M$ does not have a spin structure.
\end{remark}

\begin{remark}
In item~(1) of \cref{spin structures homogeneous spaces}, the injectivity of
$$\pi^*:\Ext^1_\Z(\pi_1(G/H),\Z/2)\to \Ext^1_\Z(\pi_1(G),\Z/2)$$ is equivalent, by the long exact sequence in the proof,
to the surjectivity of the restriction map
\[
\overline{i}^{\,*}:\Hom(\pi_1(G),\Z/2)\longrightarrow \Hom(\pi_1(H)/\ker i_*,\Z/2).
\]

In particular, we can apply \cref{spin structures homogeneous spaces}.(1) in each of the following situations:
\begin{itemize}
\item If $\pi_1(G/H)$ has no $2$-torsion, then $\Ext^1_\Z(\pi_1(G/H),\Z/2)=0$, so $\pi^*$ is trivially injective.


\item If $\pi_1(G)\cong (\Z/2)^r$ and $H$ is connected, $\pi_1(G/H)\cong (\Z/2)^l$ and $\pi_1(H)/\ker i_*\cong(\Z/2)^s$ for some natural number $l$ and $s$ satisfying $l+s=r$ (this is because $\pi_1(G)\longrightarrow \pi_1(G/H)$ is surjective and $\pi_1(H)/\ker i_*$ is a subgroup of $\pi_1(G)$). If we apply the $\Hom(\cdot ,\Z/2)$ functor to the short exact sequence of abelian groups
\[
\begin{tikzcd}
0 \ar{r}{} & \pi_1(H)/\ker i_*\cong (\Z/2)^s\ar{r}{\overline{i_*}} & \pi_1(G)\cong(\Z/2)^r \ar{r}{} & \pi_1(G/H)\cong (\Z/2)^l \ar{r}{} & 0
\end{tikzcd}
\]
we obtain a surjective morphism $$\overline{i}^{\,*}:\Hom(\pi_1(G),\Z/2)\cong(\Z/2)^r\longrightarrow \Hom(\pi_1(H)/\ker i_*,\Z/2)\cong(\Z/2)^s$$ and therefore $G/H$ is spin if and only if $\sigma_*(\ker i_*)=0$.
\end{itemize}

Moreover, note that the condition $\sigma_*(\ker i_*)=0$ implies that the universal covering
$\widetilde{G/H}$ is spin. Indeed, let $\widetilde{G}\to G$ be the universal covering group and let
$\widetilde{H}\subset\widetilde{G}$ be the full preimage of $H$. Then
\(
\widetilde{G/H}\cong\widetilde{G}/\widetilde{H}^0,
\)
where $\widetilde{H}^0$ is the identity component of $\widetilde{H}$.
The covering $\widetilde{H}^0\to H$ induces an identification
\(
\pi_1(\widetilde{H}^0)\cong\ker i_*.
\)
Since $\sigma_*(\ker i_*)=0$, the isotropy representation of
$\widetilde{H}^0$ lifts to $\Spin(n)$.
Consequently, $\widetilde{G/H}$ is spin.
\end{remark}

\begin{corollary}\label{invariant ss sphere}
Let $G$ be a connected compact Lie group acting by orientation-preserving isometries on the \emph{round} sphere
$S^n$. If $n\ge 2$, then the unique spin structure of $S^n$ is $G$-invariant if and only if
$ev_*^{G\acts F(S^n)}=0$.
\end{corollary}

\begin{proof}
For the round metric one has an identification of the oriented orthonormal frame bundle
$F(S^n)\cong \SO(n+1)$. If $n\ge 2$, then $\pi_1(F(S^n))\cong \pi_1(\SO(n+1))\cong \Z/2$, and the unique spin
structure corresponds (via the Classification via $\pi_1$ lemma) to the nontrivial homomorphism
$f:\pi_1(F(S^n))\to \Z/2$, i.e.\ $f=\id_{\Z/2}$.
Therefore, by \cref{invariant spin structure criteria}, the unique spin structure on $S^n$ is $G$-invariant if and
only if $f\circ ev_*^{G\acts F(S^n)}=0$, which is equivalent to $ev_*^{G\acts F(S^n)}=0$.
\end{proof}

\begin{remark}
The action of $G$ on the round sphere induces a Lie group homomorphism
\[
i:G\longrightarrow \Isom_0(S^n)\cong \SO(n+1).
\]
Note that if the action is effective then $i$ is injective.
Under the identification $F(S^n)\cong \SO(n+1)$, the evaluation map for the $G$-action on $F(S^n)$ equals $i$, hence
\[
ev_*^{G\acts F(S^n)}=i_*:\pi_1(G)\to \pi_1(\SO(n+1)).
\]
\end{remark}

\begin{remark}
The case $n=1$ has some differences. The oriented orthonormal frame bundle is $F(S^1)\cong S^1$.
It has two spin structures, corresponding to the two principal $\Spin(1)\cong \Z/2$-bundles over $S^1$:
the disconnected covering $S^1\sqcup S^1\to S^1$ (trivial bundle) and the connected covering
$S^1\to S^1$ (the degree-$2$ covering).
Equivalently, these correspond to the two homomorphisms $f:\pi_1(S^1)\cong \Z\to \Z/2$, namely $f=0$ and
$f=\text{reduction mod }2$.

The spin structure associated to the disconnected covering (i.e.\ $f=0$) is always $G$-invariant.
The spin structure associated to the connected covering (i.e.\ $f=\text{mod }2$) is $G$-invariant if and only if
\[
ev_*^{G\acts F(S^1)}(\pi_1(G))\subseteq 2\Z\subseteq \pi_1(F(S^1))\cong \Z.
\]
\end{remark}

\subsection{Invariant spin structures and free actions of discrete groups on manifolds}
\label{sec:Invariant spin structures and discrete group actions}

Let $\Gamma$ be a discrete group acting properly (see \cite[\S 1.2]{lee2010seifert} for the precise definition) and by orientation preserving isometries on a connected spin manifold $M$. This action lifts to an action of $\Gamma$ on the total space of the oriented frame bundle $FM\longrightarrow M$, which is proper by \cite[Theorem 1.9.3]{lee2010seifert}. Let $(P,\Lambda)$ be a spin structure on $M$, represented by a class $c\in H^1(FM,\Z/2)$ with the usual fibre restriction. The spin structure is $\Gamma$-invariant if the following two conditions hold:
\begin{itemize}
    \item[1.] The isomorphism class of the spin structure is preserved by $\Gamma$, meaning that $c$ is fixed by the
    induced action of $\Gamma$ on $H^1(FM,\Z/2)$.
    In this case, for each $\gamma\in\Gamma$ there exists a lift of the induced diffeomorphism
    $\gamma:FM\to FM$ to a bundle automorphism $\widetilde\gamma:P\to P$ commuting with the right
    $\Spin(n)$-action and satisfying $\Lambda\circ \widetilde\gamma=\gamma\circ\Lambda$.
    The set of all such lifts forms a group $\widetilde{\Gamma}$ fitting into a short exact sequence
    \[
    \begin{tikzcd}
        1\ar{r}{} &\Z/2\ar{r}{} &\widetilde{\Gamma}\ar{r}{\phi} &\Gamma\ar{r}{}& 1,
    \end{tikzcd}
    \]
    where $\phi(\widetilde\gamma)=\gamma$ and the kernel $\Z/2$ is the group of deck transformations of
    $\Lambda:P\to FM$. The action of $\widetilde{\Gamma}$ on $P$ is proper by \cite[Theorem 1.9.3]{lee2010seifert}, and $\Lambda$ is $\phi$-equivariant:
    \[
    \Lambda(\widetilde\gamma\cdot x)=\phi(\widetilde\gamma)\cdot \Lambda(x)\qquad
    \text{for all }\widetilde\gamma\in\widetilde\Gamma,\ x\in P.
    \]

    \item[2.] There exists a group homomorphism $s:\Gamma\longrightarrow \widetilde{\Gamma}$ such that $\phi\circ s=\id_\Gamma$ and hence the above extension splits. In this case, $\Gamma$ acts properly on $P$ via $\gamma\cdot x := s(\gamma)\cdot x$, and $\Lambda$ is $\Gamma$-equivariant. Note that the above short exact sequence is central, therefore it splits if and only if it is isomorphic to the trivial extension
    \[
    \begin{tikzcd}
        1\ar{r}{} &\Z/2\ar{r}{} &\Z/2\times \Gamma\ar{r}{} &\Gamma\ar{r}{}& 1.
    \end{tikzcd}
    \]
\end{itemize}

\begin{remark}\label{lift Z}
There are cases where one of two conditions automatically hold. For example, assume that we have a connected Lie group $G$ acting properly and by orientation preserving isometries on a spin manifold $M$. Let $\Gamma$ be a discrete subgroup of $G$, then the action of $\Gamma$ on $H^1(FM,\Z/2)$ is trivial since it comes from restricting the trivial action of $G$ on $H^1(FM,\Z/2)$. Consequently, a spin structure is invariant if and only if there exists a section of the group morphism $\phi:\Tilde{\Gamma}\longrightarrow \Gamma$.

Another interesting case is the following. Assume that $\Gamma\cong \Z$ and that the proper action on $M$ is generated by an orientation-preserving isometry $f':M\longrightarrow M$. Let $f:FM\longrightarrow FM$ denote the lifting of $f'$ on the oriented frame bundle $FM$ and let $c\in H^1(FM,\Z/2)$ be a spin structure satisfying $f^*(c)=c$. Since $\Z$ is a free group, the short exact sequence 
    \[
    \begin{tikzcd}
        1\ar{r}{} &\Z/2\ar{r}{} &\widetilde{\Gamma}\ar{r}{\phi} &\Z\ar{r}{}& 1,
    \end{tikzcd}
    \]
always splits and hence the second condition is always satisfied. Thus, a spin structure $c$ is $\Z$-invariant if and only if $f^*(c)=c$.
\end{remark}

We now assume that $\Gamma$ is a discrete group acting freely, properly, and by orientation-preserving isometries
on a connected spin manifold $M$. Then the quotient $M/\Gamma$ is a smooth connected manifold and the projection
$\pi:M\longrightarrow M/\Gamma$ is a regular $\Gamma$-covering (see \cite[Chapter 1]{lee2010seifert}). Endowing $M/\Gamma$ with the unique Riemannian metric for which $\pi$ is a local isometry, the map $\pi$ becomes a Riemannian covering. Since the action preserves the orientation, the quotient manifold $M/\Gamma$ is orientable. This action induces an action of $\Gamma$ on the total space of the oriented frame bundle $FM$ which is free and proper (again by \cite[Theorem 1.9.3]{lee2010seifert}) and consequently, $FM\longrightarrow (FM)/\Gamma$ is a regular $\Gamma$-covering. The frame bundle $q:FM\longrightarrow M$ is $\Gamma$-equivariant and the induced map
\[q':FM/\Gamma\longrightarrow M/\Gamma\]
is a principal $\SO(n)$-bundle.

\begin{lemma}\label{spin structures orbit spaces discrete}
Let $M$ be a spin manifold and let $\Gamma$ be a discrete group acting properly, freely and by orientation-preserving isometries on $M$. Then $M/\Gamma$ is spin if and only if $M$ has a $\Gamma$-invariant spin structure.
\end{lemma}

The proof of \cref{spin structures orbit spaces discrete} is a direct consequence of the following fact.

\begin{lemma}\label{frame bundle discrete group}
The principal $\SO(n)$-bundles $FM/\Gamma\longrightarrow M/\Gamma$ and $F(M/\Gamma)\longrightarrow M/\Gamma$ are
isomorphic.
\end{lemma}

\begin{proof}
Since $\pi:M\longrightarrow M/\Gamma$ is a Riemannian covering, $\pi$ is a local isometry. Hence for each $x\in M$ the differential $d\pi_x:T_xM\to T_{\pi(x)}(M/\Gamma)$ is a linear isometry.

Define
\[
\Theta:FM/\Gamma \longrightarrow F(M/\Gamma)
\]
by
\[
\Theta\bigl([u]\bigr):=d\pi_x\circ u,
\]
where $u:\R^n\to T_xM$ is an oriented orthonormal frame at $x$ and $[u]$ denotes its $\Gamma$-orbit in $FM$.
This is well-defined: if $u'=\gamma\cdot u$ for $\gamma\in\Gamma$, then $\pi(\gamma x)=\pi(x)$ and
$d\pi_{\gamma x}\circ d\gamma_x=d\pi_x$, hence
\[
d\pi_{\gamma x}\circ u' \;=\; d\pi_{\gamma x}\circ d\gamma_x\circ u \;=\; d\pi_x\circ u.
\]
The map $\Theta$ is smooth and $\SO(n)$-equivariant, and it covers $\id_{M/\Gamma}$. It is a diffeomorphism.
Its inverse sends a frame $v$ at $\bar x\in M/\Gamma$ to the class $[u]\in FM/\Gamma$,
where $x\in M$ is any lift of $\bar x$ and $u=(d\pi_x)^{-1}\circ v$; this is independent
of the choice of lift because changing $x$ to $\gamma x$ replaces $u$ by $\gamma\cdot u$. Hence $\Theta$ is an isomorphism of principal $\SO(n)$-bundles.
\end{proof}

\begin{proof}[Proof of \cref{spin structures orbit spaces discrete}]
Assume that $(P,\Lambda)$ is a $\Gamma$-invariant spin structure on $M$. Then $\Gamma$ acts freely and properly on
$P$, the action commutes with the right $\Spin(n)$-action, and $\Lambda$ is $\Gamma$-equivariant. Let
$\overline{\pi}:P\to P/\Gamma$ and $\overline{\pi}_F:FM\to FM/\Gamma$ denote the quotient maps. Since the
$\Spin(n)$-action commutes with $\Gamma$, the quotient $P/\Gamma\to M/\Gamma$ is a principal $\Spin(n)$-bundle.
Moreover, $\Lambda$ descends to a bundle map
\[
\overline{\Lambda}:P/\Gamma \longrightarrow FM/\Gamma,\qquad \overline{\Lambda}([p])=[\Lambda(p)],
\]
which restricts on fibres to a $2$-fold covering. Composing with the bundle isomorphism
$\Theta:FM/\Gamma\to F(M/\Gamma)$ from \cref{frame bundle discrete group}, we obtain
\[
\Lambda':=\Theta\circ \overline{\Lambda}:P/\Gamma \longrightarrow F(M/\Gamma).
\]
Thus $(P/\Gamma,\Lambda')$ is a spin structure on $M/\Gamma$. Hence $M/\Gamma$ is spin.

Conversely, let $(P',\Lambda')$ be a spin structure on $M/\Gamma$. Pulling back along $\pi$ yields a principal
$\Spin(n)$-bundle $\pi^*P'\to M$ together with a bundle map
\[
\pi^*\Lambda':\pi^*P' \longrightarrow \pi^*F(M/\Gamma)
\]
restricting to a $2$-fold covering on fibres. The pullback bundle $\pi^*P'$ carries a natural action of $\Gamma$ which is free, proper, commutes with the $\Spin(n)$-action, and $\pi^*\Lambda'$ is $\Gamma$-equivariant. Hence $\pi^*F(M/\Gamma)\longrightarrow M$ has a $\Gamma$-invariant spin structure. By \cref{frame bundle discrete group}, we have an isomorphism of principal $\SO(n)$-bundles between $FM\longrightarrow M$ and $\pi^*F(M/\Gamma)\longrightarrow M$. Hence transporting the $\Gamma$-invariant spin structure on the pullback bundle via this isomorphism gives a $\Gamma$-invariant spin structure on $FM\to M$. Therefore $M$ has a $\Gamma$-invariant spin structure.
\end{proof}

\Cref{frame bundle discrete group} shows that if $(P,\Lambda)$ is a $\Gamma$-invariant spin structure on $M$, then there is an induced spin structure on $M/\Gamma$. If $f:\pi_1(FM)\to \Z/2$ denotes the morphism induced by $(P,\Lambda)$, then there exists a spin structure $\overline{f}:\pi_1(F(M/\Gamma))\to \Z/2$ such that
\[
\overline{f}\circ \Pi_* = f,
\]
where $\Pi:=\Theta\circ \overline{\pi}_F:FM\to F(M/\Gamma)$ is the canonical $\Gamma$-invariant map induced by
$d\pi$. Conversely, if $\overline{f}:\pi_1(F(M/\Gamma))\to \Z/2$ is a spin structure on $M/\Gamma$, then
$\overline{f}\circ \Pi_*$ is a $\Gamma$-invariant spin structure on $M$.


We will use the following result to count invariant spin structures under finite cyclic coverings.

For a discrete group \(\Gamma\) acting on a spin manifold \(M\), we denote by \(\Spin_\Gamma(M)\) the set of isomorphism classes of spin structures on \(M\) for which the action of \(\Gamma\) admits a lift to the spin bundle.

\begin{lemma}\label{invariant spin structures cyclic groups} Let $\pi:M\longrightarrow M'$ be a regular $\Z/a$-covering of closed connected $n$-dimensional manifolds, where $a$ is even. Assume that $M$ and $M'$ are spin. Then $|\Spin_{\Z/a}(M)|=\tfrac{1}{2}|\Spin(M')|$.
\end{lemma}

\begin{proof}
By \cref{spin structures orbit spaces discrete},
the $\Z/a$-invariant spin structures on $M$ are exactly
the pullbacks of spin structures on $M'$.

The induced covering $\Pi:FM\to FM'$ gives an exact sequence
\[
1\to\pi_1(FM)\xrightarrow{\Pi_*}\pi_1(FM')
\to\Z/a\to1.
\]
By \cref{Classification via FG}, two spin structures
$f_1,f_2:\pi_1(FM')\to\Z/2$ have the same pullback
if and only if $f_1+f_2$ vanishes on $\im\Pi_*$,
or equivalently, factors through $\Z/a$.

Conversely, adding any homomorphism
$\Z/a\to\Z/2$ preserves the fibre condition,
since the inclusion of each \(\SO(n)\)-fibre in \(FM'\) factors through \(\Pi\).
As $\Hom(\Z/a,\Z/2)\cong\Z/2$, every nonempty
fibre of the pullback map has exactly two elements.
Hence
\[
|\Spin_{\Z/a}(M)|=\tfrac12|\Spin(M')|.
\]
\end{proof}

Like in the connected case, \cref{spin structures orbit spaces discrete} has a counterpart for invariant spin structures. Assume that we have a connected Lie group $H$ acting on $M$ and commuting with the action of $\Gamma$. Then $\pi:M\to M/\Gamma$ is $H$-equivariant, and the induced map $\Pi:FM\to F(M/\Gamma)$ is also $H$-equivariant.
In particular, one has
\[
\Pi_*\circ ev_*^{H\acts FM} \;=\; ev_*^{H\acts F(M/\Gamma)}.
\]
Consequently, the same argument as in the proof of \cref{invariant spin structure pullback invariant} yields:

\begin{lemma}\label{spin structures orbit spaces discrete invariant}
Let $M$ be a spin manifold and let $\Gamma$ be a discrete group acting properly, freely and by orientation-preserving isometries on $M$. Assume that a connected Lie group $H$ acts on $M$ and commutes with the action of $\Gamma$. Then $M/\Gamma$ has a $H$-invariant spin structure if and only if $M$ has a $\Gamma$-invariant spin structure which is also $H$-invariant.
\end{lemma}

\begin{proof}
Let $\overline{\pi}_F:FM\to FM/\Gamma$ be the quotient map and let
$\Theta:FM/\Gamma\to F(M/\Gamma)$ be the principal $\SO(n)$-bundle isomorphism of
\cref{frame bundle discrete group}. Set
$\Pi:=\Theta\circ \overline{\pi}_F:FM\longrightarrow F(M/\Gamma).$ 
Since the $H$- and $\Gamma$-actions on $M$ commute, the induced actions on $FM$ commute as well, and hence
$\overline{\pi}_F$ is $H$-equivariant. Moreover, $\Theta$ is $H$-equivariant by construction, since it is induced by
$d\pi$. Therefore $\Pi$ is $H$-equivariant. In particular, we have
$\Pi\circ ev_u=ev_{\Pi(u)}$ for every $u\in FM$, and hence
\[
\Pi_*\circ ev_*^{H\acts FM}
=ev_*^{H\acts F(M/\Gamma)}.
\]

Assume first that $M/\Gamma$ admits a $H$-invariant spin structure, represented by a homomorphism $\overline{f}:\pi_1(F(M/\Gamma))\longrightarrow \Z/2$
satisfying the fibre condition. Pulling back along $\Pi$ defines
\[
f:=\overline{f}\circ \Pi_*:\pi_1(FM)\longrightarrow \Z/2.
\]
Then $f$ is a spin structure on $M$, and it is
$\Gamma$-invariant by construction. Moreover, since $\overline{f}$ is $H$-invariant, we have
\[
f\circ ev_*^{H\acts FM}
=\overline{f}\circ \Pi_*\circ ev_*^{H\acts FM}
=\overline{f}\circ ev_*^{H\acts F(M/\Gamma)}=0,
\]
so $f$ is also $H$-invariant. Thus $M$ admits a $\Gamma$-invariant spin structure which is $H$-invariant.

Conversely, assume that $M$ admits a spin structure
$f:\pi_1(FM)\longrightarrow \Z/2$
which is $\Gamma$-invariant and $H$-invariant. By \cref{spin structures orbit spaces discrete}, the
$\Gamma$-invariance implies that $f$ descends to a spin structure on $M/\Gamma$; equivalently, there exists a spin
structure $\overline{f}:\pi_1(F(M/\Gamma))\longrightarrow \Z/2$ such that $\overline{f}\circ \Pi_*=f$. Then, using again
$\Pi_*\circ ev_*^{H\acts FM}=ev_*^{H\acts F(M/\Gamma)}$, we obtain
\[
\overline{f}\circ ev_*^{H\acts F(M/\Gamma)}
=\overline{f}\circ \Pi_*\circ ev_*^{H\acts FM}
=f\circ ev_*^{H\acts FM}=0,
\]
so $\overline{f}$ is $H$-invariant. Hence $M/\Gamma$ has a $H$-invariant spin structure.
\end{proof}

\section{Spin structures on cohomogeneity one manifolds of type (1) and type (3)}\label{sec: spin structures type (1) and (3)}

The first case we consider is when $M$ admits a cohomogeneity one action of a compact Lie group $G$ of type (3).
Recall that in this case $M$ is equivariantly diffeomorphic to the product $G/H\times \R$ with the trivial action
of $G$ on $\R$, for some closed subgroup $H$ of $G$.

\begin{theorem}\label{main_theorem (3) section}
Assume that $M$ admits a cohomogeneity one action of a compact Lie group $G$ of type (3). Then $M$ is spin if and
only if $G/H$ is spin. Moreover, $M$ has a $G$-invariant spin structure if and only if $G/H$ has a $G$-invariant
spin structure.
\end{theorem}

\begin{proof}
Let $p_1:G/H\times \R\longrightarrow G/H$ and $p_2:G/H\times\R\longrightarrow \R$ denote the usual projection maps.
Note that both of them are $G$-equivariant. We have an isomorphism of $G$-equivariant vector bundles
\[
TM\cong p_1^* T(G/H)\oplus p_2^* T\R.
\] Since $p_1$ is a homotopy equivalence, the pullback along $p_1$ induces a bijection between spin structures on $G/H$ and spin structures on $p_1^*T(G/H)$. Moreover, since $p_1$ is $G$-equivariant, this bijection preserves $G$-invariance. On the other hand, $p_2^*T\R$ is a trivial line bundle with trivial $G$-action and has a unique spin structure, which is $G$-invariant. By \cref{2 out of 3 homotopy},
$M$ is spin if and only if $G/H$ is spin. Furthermore, $M$ has a $G$-invariant spin structure if and only if $G/H$ has a $G$-invariant spin structure.
\end{proof}

We consider now the case where $M$ admits a cohomogeneity one action of a compact Lie group $G$ of type (1).
Recall that we are in the following setting: $M$ is a compact manifold and $G$ is a compact connected Lie group
acting effectively on $M$ such that $M/G\cong S^1$. The orbit map $p:M\longrightarrow S^1$ is a fibre bundle with
fibre $G/H$ and structure group $N_G(H)/H$. Consequently, $M$ is diffeomorphic to the mapping torus
\[
(G/H)_f := (G/H\times[0,1])/(gH,0)\sim(f(gH),1).
\]
The map $f:G/H\longrightarrow G/H$ is a $G$-equivariant isometry and hence there exists $x\in N_G(H)$ such that
$f(gH)=g xH$ for all $gH\in G/H$. Thus, we can restrict the structure group of the bundle $M\longrightarrow S^1$
from $N_G(H)/H$ to the cyclic subgroup generated by $xH$ in $N_G(H)/H$. The main result of this section is the
following:

\begin{theorem}\label{main_theorem (1) section}
Assume that $M$ admits a cohomogeneity one action of a compact Lie group $G$ of type (1). Suppose that $M$ is
orientable and that $H$ is connected. Then:
\begin{itemize}
    \item[1.] $M$ is spin if and only if $G/H$ is spin.
    \item[2.] $M$ has a $G$-invariant spin structure if and only if $G/H$ has a $G$-invariant spin structure.
\end{itemize}
\end{theorem}

To prove \cref{main_theorem (1) section} we need the following lemma. This characterisation of spin structures on mapping tori is well known to experts, although an explicit reference seems difficult to find in the literature.

\begin{lemma}\label{mapping torus}
Let $N$ be a connected smooth manifold and let $f$ be a diffeomorphism. Let $N_f$ denote the associated mapping torus. Then:
\begin{itemize}
	\item[1] $N_f$ is orientable if and only if $N$ is orientable and $f$ preserves the orientation.
	\item[2] $N_f$ is spin if and only if $N$ is spin and $f$ preserves a spin structure of $N$.
\end{itemize}
\end{lemma}

\begin{proof}
By taking the pullback of the map $N_f\longrightarrow S^1$
by the usual principal $\Z$-bundle $\pi:\R\longrightarrow S^1$,
we can describe $N_f$ as the quotient $(N\times\R)/\Z$,
where the action satisfies
\[
v\cdot(x,t)=(f^{\,v}(x),t+v)
\qquad\text{for all }v\in\Z\text{ and }(x,t)\in N\times\R.
\]
The first equivalence follows from the fact that
$(N\times\R)/\Z$ is orientable if and only if the action
preserves the orientation of $N\times\R$, which is equivalent
to $f$ preserving the orientation of $N$.

For the second equivalence, assume that $N$ is oriented
and that $f$ preserves its orientation. Since the $\Z$-action
on $N\times\R$ is free and proper, we may choose a
$\Z$-invariant Riemannian metric. By
\cref{spin structures orbit spaces discrete}, $N_f$ is spin
if and only if $N\times\R$ admits a $\Z$-invariant spin
structure for the above action. By \cref{lift Z}, this is
equivalent to the existence of a spin structure on
$N\times\R$ whose isomorphism class is preserved by the
deck transformation $(x,t)\mapsto(f(x),t+1)$.

The projection onto the first component
$pr_1:N\times\R\longrightarrow N$ induces a correspondence
between spin structures on $N$ and spin structures on
$N\times\R$, by pulling back along $pr_1$ and using the
non-equivariant version of the two-out-of-three lemma
(see the beginning of the proof of
\cref{2 out of 3 homotopy}) with the unique spin structure
on $\R$. Since $\R$ is contractible, this correspondence
is bijective.

Moreover, $pr_1$ is $\Z$-equivariant, where $\Z$ acts on
$N$ through the iterates of $f$. Consequently, the
correspondence is $\Z$-equivariant. Thus $N\times\R$
admits a spin structure preserved by the deck transformation
if and only if $N$ admits a spin structure preserved by $f$.
This completes the proof.
\end{proof}


\begin{remark}\label{mapping torus surfaces}
If $n=1$, then $G/H$ is a closed connected $1$-manifold, hence $G/H\cong S^1$. Given a diffeomorphism $f:S^1\longrightarrow S^1$, then the mapping torus $(S^1)_f$ is the total space of a fiber bundle with fibre and basis $S^1$. In particular, $(S^1)_f$ is either diffeomorphic to the torus $T^2$ or the Klein bottle $K$. If we assume that $f$ is orientation preserving then $(S^1)_f$ is orientable by \cref{mapping torus} and hence it is diffeomorphic to $T^2$. Consequently, $(S^1)_f$ is spin and by \cref{mapping torus} $f$ preserves the two spin structures of $S^1$. 

We also point out that \Cref{mapping torus} can be used to easily prove a result of M.~Atiyah in \cite{atiyah1971riemann} (see also \cite{kallel2010invariant,johnson1980spin}). Namely, every orientation-preserving isometry $f$ on a closed connected Riemannian surface $\Sigma$ preserves a spin structure. Indeed, since $f$ preserves the orientation, the mapping torus $\Sigma_f$ is a closed connected orientable $3$-manifold, and hence it is spin (see, for example,
\cite{benedetti2019framing}). By \cref{mapping torus}, we can conclude that $f$ preserves a spin structure of $\Sigma$.

The result does not extend to $3$-manifolds. For example, consider the manifold $S^2\times S^1$ and let
$f:S^2\times S^1\longrightarrow S^2\times S^1$ be the diffeomorphism defined by $f(x,\theta)=(R_\theta(x),\theta)$,
where $R_\theta$ denotes the rotation around the $z$-axis on $S^2$ with angle $\theta$. This diffeomorphism generates the kernel of the morphism from the mapping class group to the automorphism group of cohomology
\[
\operatorname{MCG}(S^2\times S^1)\longrightarrow \Aut(H^*(S^2\times S^1,\Z)),
\]
see \cite{gluck1962embedding}. A straightforward computation shows that $f$ permutes the two spin structures of
$S^2\times S^1$. Indeed, $S^2\times S^1$ admits a transitive group action of $\SO(3)\times S^1$, where $\SO(3)$ acts on the first
component and $S^1$ acts on the second component. Let $p=((0,0,1),0)\in S^2\times S^1$, then the isotropy subgroup
of $p$ is $\{(R_\theta,0):\theta\in [0,2\pi]\}\cong S^1$. Hence
$S^2\times S^1\cong (\SO(3)\times S^1)/S^1$ and the isotropy representation at $p$,
$\sigma:S^1\longrightarrow \SO(3)$, satisfies $\sigma(R_\theta)=R_\theta$.
Consequently,
\[
F(S^2\times S^1)\cong (\SO(3)\times S^1)\times_\sigma\SO(3).
\]
Then $f:S^2\times S^1\longrightarrow S^2\times S^1$ lifts to a map
$\widetilde{f}:F(S^2\times S^1)\longrightarrow F(S^2\times S^1)$ satisfying
\[
\widetilde{f}([(A,\theta),B])=[(R_\theta A,\theta),B].
\]

Using the long exact sequence of homotopy for the frame bundle of $S^2\times S^1$ and the fact that
$S^2\times S^1$ is spin, we obtain an exact sequence
\[
1\longrightarrow \pi_1(\SO(3))\cong \Z/2 \longrightarrow \pi_1(F(S^2\times S^1))
\longrightarrow \pi_1(S^2\times S^1)\cong \Z \longrightarrow 1.
\]
Since $H^2(\Z,\Z/2)=0$, this extension splits, and hence
\[
\pi_1(F(S^2\times S^1))\cong \Z/2\oplus\Z.
\]
It is generated by the classes of loops $[\alpha]$ and $[\widetilde{\gamma}]$, where
\[
\alpha:I\longrightarrow F(S^2\times S^1), \qquad t\mapsto [(Id,0),R_{2\pi t}]
\]
and
\[
\widetilde{\gamma}:I\longrightarrow F(S^2\times S^1), \qquad t\mapsto [(Id,2\pi t),Id].
\]
Note that $\widetilde{f}(\alpha)=\alpha$ and
\[
\widetilde{f}(\widetilde{\gamma})(t)=[(R_{2\pi t},2\pi t), Id]=[(Id,2\pi t), R_{-2\pi t}].
\]
Therefore, $\widetilde{f}_*([\widetilde{\gamma}])=[\alpha]+[\widetilde{\gamma}]$.

The spin structures of $S^2\times S^1$ are represented by group homomorphisms
$\pi_1(F(S^2\times S^1))\longrightarrow \Z/2$ which send $[\alpha]$ to $1$. Therefore, $S^2\times S^1$ has two spin
structures represented by the group morphisms $f_i\in\Hom(\pi_1(F(S^2\times S^1)),\Z/2)$ satisfying that
$f_i([\alpha])=1$ and $f_i([\widetilde{\gamma}])=i$, where $i=0,1$. The above computation shows that
$f_0\circ \widetilde{f}_*=f_1$, which implies that $f$ permutes the two spin structures. By \cref{mapping torus}, $(S^2\times S^1)_f$ is not spin. 
\end{remark}

We are ready to prove \cref{main_theorem (1) section}.

\begin{proof}[Proof of \cref{main_theorem (1) section}]

We start by proving the first part of \cref{main_theorem (1) section}, that $M$ is spin if and only if $G/H$ is spin. 

Note that if $M$ is spin, then so is $G/H$ by \cref{mapping torus}. Thus it remains to prove that if $G/H$ is spin and $H$ is connected, then $M$ is spin. Let $n=\dim(G/H)$. Since $G$ is compact, we can choose an
$\Ad(N_G(H))$-invariant inner product on
$\mathfrak g$. Setting $\mathfrak m=\mathfrak h^\perp$,
we obtain an $\Ad(N_G(H))$-invariant decomposition
$\mathfrak g=\mathfrak h\oplus\mathfrak m$,
where $\mathfrak m\cong T_{eH}(G/H)$.

By \cref{mapping torus}, it suffices to prove that $G/H$ has a spin structure which is preserved by the element $f\in N_G(H)/H$ such that $(G/H)_f\cong M$. Since $M$ is orientable, the induced diffeomorphism $f$ on $G/H$ preserves the orientation. By \cref{mapping torus surfaces}, the statement holds for $n\le 2$. Thus we may assume
$n\ge 3$, and hence $\pi_1(\SO(n))\cong \Z/2$.

Let $A\le N_G(H)/H$ be the subgroup formed by orientation-preserving isometries of $G/H$, which is either the whole $N_G(H)/H$ or an index 2 subgroup of $N_G(H)/H$. Since $M$ is orientable, $\langle f\rangle\le A$ by \cref{mapping torus}. We will prove that every spin structure on $G/H$ is fixed by the action of $A$. 

Let $A'\le N_G(H)$ denote the preimage of $A$ under the projection $N_G(H)\to N_G(H)/H$. We have a free right action of $H$ on $G\times\SO(n)$ defined by
\[
(g,B)\cdot h=(gh,\Ad(h^{-1})_{|\mathfrak{m}}\,B),
\]
which induces a principal $H$-bundle
\[
q:G\times \SO(n)\longrightarrow G\times_\sigma\SO(n).
\]
Similarly, we define a right action of $A'$ on $G\times \SO(n)$ by
\[
(g,B)\cdot x=(gx,\Ad(x^{-1})_{|\mathfrak{m}}\,B).
\]
This is well defined since $A'$ consists of elements whose adjoint action on $\mathfrak{m}$ is orientation preserving. This action descends to $G\times_\sigma\SO(n)$ and the map $q$ is $A'$-equivariant. Note that the action of $A'$ on $G\times_\sigma\SO(n)$ is not effective, and $H$ is contained in the kernel. The induced action of $A$ on $G\times_\sigma\SO(n)$ coincides with the action obtained by lifting the right action of $A$ on $G/H$ to the frame bundle.

We claim that the action of $A'$ on $H^1(G\times\SO(n),\Z/2)$ is trivial. Indeed,
\[
H^1(G\times\SO(n),\Z/2)\cong H^1(G,\Z/2)\oplus H^1(\SO(n),\Z/2).
\]
Since $G$ is connected and $A'\le G$, the action of $A'$ by right multiplication on $G$ induces the trivial action on $H^1(G,\Z/2)$. For the second summand, let $\gamma:I\longrightarrow \SO(n)$ be a loop whose homotopy class $[\gamma]$ generates $\pi_1(\SO(n))$. For all $x\in A'$, the map $\Ad(x^{-1})_{|\mathfrak{m}}$ is an element of $\SO(n)$, and hence $[\Ad(x^{-1})_{|\mathfrak{m}}\circ\gamma]$ is again a generator of $\pi_1(\SO(n))$. Since $\pi_1(\SO(n))\cong \Z/2$, it follows that $A'$ acts trivially on $\pi_1(\SO(n))$, and therefore trivially on $H^1(\SO(n),\Z/2)$. These two facts imply that $A'$ fixes the generators of $H^1(G\times\SO(n),\Z/2)$, hence acts trivially on $H^1(G\times\SO(n),\Z/2)$.

Since $H$ is connected, the map
\[
q^*:H^1(G\times_\sigma\SO(n),\Z/2)\longrightarrow H^1(G\times\SO(n),\Z/2)
\]
is injective. Moreover, it is $A'$-equivariant. Consequently, the action of $A'$ on $H^1(G\times_\sigma\SO(n),\Z/2)$ is trivial, which implies that the induced action of $A$ on $H^1(G\times_\sigma\SO(n),\Z/2)$ is also trivial. Since the set of
spin structures on $G/H$ is a subset of $H^1(G\times_\sigma\SO(n),\Z/2)$, the action of $A$ fixes all spin structures on $G/H$. In particular, $f$ preserves a spin structure on $G/H$, and therefore $M$ is spin by \cref{mapping torus}.


We now prove the second part of \cref{main_theorem (1) section}, that $M$ has a $G$-invariant spin structure if and only if $G/H$ has a $G$-invariant spin structure. The main observation to prove this result is the fact that the left action of $G$ on $G/H$ and the right action of $N_G(H)/H$ on $G/H$ commute.

Assume that $M$ admits a $G$-invariant spin structure. Let $q:G/H\times \R\longrightarrow M$ be the covering obtained from the pullback of $M\to S^1$ along $\R\to S^1$. It is $G$-equivariant since the left action of $G$ on $G/H\times \R$ and the deck transformation action of $\Z$ on $G/H\times \R$ commute. By \cref{spin structures orbit spaces discrete invariant}, $G/H\times \R$ has a $G$-invariant spin structure. Note that the action of $G$ on the second factor of $G/H\times \R$ is  trivial, and therefore the spin structure on $\R$ is $G$-invariant. Using \cref{2 out of 3 homotopy}, we obtain that $G/H$ has a $G$-invariant spin structure.

Conversely, let $f\in N_G(H)/H$ be such that $(G/H)_f\cong M$ and assume that $G/H$ has a $G$-invariant spin structure, represented by a class $c\in H^1(F(G/H),\Z/2)$. We have that $\langle f\rangle\leq A$ acts trivially on $H^1(F(G/H),\Z/2)$, and in particular, it fixes $c$. Consequently, the product spin structure on $G/H\times \R$ given by $c$ on $G/H$ and the trivial spin structure on $\R$ is $\Z$-invariant under the deck transformation action, as well as $G$-invariant by the action which is the transitive action of $G$ on $G/H$ and trivial on $\R$. By \cref{spin structures orbit spaces discrete invariant}, this $G$-invariant spin structure descends to a $G$-invariant spin structure on $M$.


\end{proof}
We conclude this section by showing that the
connectedness hypothesis on $H$ in
\cref{main_theorem (1) section} is essential.

\begin{proposition}\label{disconnected isotropy counterexample}
There exist a compact connected Lie group $G$,
a closed disconnected subgroup $H\subseteq G$,
and an orientation-preserving $G$-equivariant
diffeomorphism $f:G/H\to G/H$ such that
$G/H$ has exactly two spin structures,
which are exchanged by $f$.

Consequently, the mapping torus $(G/H)_f$
is an orientable cohomogeneity one manifold
of type~(1) which is not spin.
\end{proposition}

To prove the above proposition we need to recall a few facts about the spin group and its realization inside the Clifford algebra (see \cite{friedrich2000dirac, lawson2016spin}).
Let $(\R^n,\langle\cdot,\cdot\rangle)$ be the standard Euclidean vector space. The (real) Clifford algebra $Cl_n$
is the unital associative $\R$-algebra generated by $\R^n$ subject to the relations
\[
v\cdot v = -\langle v,v\rangle\cdot 1 \qquad\text{for all }v\in \R^n.
\]
It is $\Z/2$-graded as $Cl_n=Cl_n^{\mathrm{even}}\oplus Cl_n^{\mathrm{odd}}$. We view $\R^n$ as a subspace of $Cl_n$
via the canonical inclusion.

For a unit vector $u\in S^{n-1}\subset \R^n$, conjugation in $Cl_n$ defines an orthogonal reflection of $\R^n$:
\[
\rho_u:\R^n\to \R^n,\qquad \rho_u(v)= -u\,v\,u^{-1}=v-2\langle u,v\rangle u,
\]
where $u^{-1}=-u$ since $u^2=-1$. Products of an even number of such reflections lie in $\SO(n)$.

The spin group is the subgroup
\[
\Spin(n):=\left\langle u_1u_2\cdots u_{2k}\ \middle|\ k\ge 0,\ u_i\in S^{n-1}\subset \R^n\right\rangle
\subset Cl_n^{\mathrm{even}}.
\]
It acts on $\R^n\subset Cl_n$ by conjugation, and this gives the standard double covering
\[
\lambda_n:\Spin(n)\longrightarrow \SO(n),\qquad \lambda_n(a)(v)=a\,v\,a^{-1}.
\]
The kernel of $\lambda_n$ is $\{\pm 1\}$, so $\lambda_n$ is a connected $2$--fold covering for $n\ge 2$.

We will use the following elementary identities. For an orthonormal basis $(e_1,\dots,e_n)$ of $\R^n$ one has
$e_ie_j=-e_je_i$ for $i\neq j$ and $e_i^2=-1$. Hence, for pairwise distinct indices $i_1,\dots,i_m$,
\[
(e_{i_1}\cdots e_{i_m})^2 = (-1)^{m(m+1)/2}.
\]
In particular, $(e_ie_j)^2=-1$ for $i\neq j$, and if $m$ is even then $(e_{i_1}\cdots e_{i_m})^2=(-1)^{m/2}$.
Moreover, $-1\in Cl_n$ is central, so conjugation by $-1$ acts trivially on $\R^n$. Finally, we use the standard blockwise inclusion $\Spin(2)^r\rightarrow\Spin(2r)$ associated to the orthogonal
decomposition $\R^{2r}=\bigoplus_{j=1}^r Span\{f_{2j-1},f_{2j}\}$.

The first step in the construction of the example is the following lemma, which is
based on \cite[\S 6]{cahen1993spin}.

\begin{lemma}\label{spin structures S2 product}
Consider the free $\Z/2$ action on $(S^2)^r$ defined by the involution
$(x_1,\dots,x_r)\mapsto (-x_1,\dots,-x_r)$. Then $(S^2)^r/(\Z/2)$ is orientable if and only if $r$ is even and
$(S^2)^r/(\Z/2)$ is spin if and only if $r=4k$ for some $k\in\mathbb{N}$.
\end{lemma}

\begin{proof}
Let $\{e_1,e_2,e_3\}$ be the standard orthonormal basis of $\R^3\subset Cl_3$. The action of $\Spin(3)$ on
$S^2\subset \R^3$ given by $x\mapsto axa^{-1}$ is transitive and the isotropy of $e_3$ is $\Spin(2)$ (whose Lie algebra
is spanned by $e_1e_2$), hence $\Spin(3)/\Spin(2)\cong S^2$. Thus $\Spin(3)^r$ acts transitively on $(S^2)^r$ with isotropy
$\Spin(2)^r$ at $(e_3^1,\dots,e_3^r)$, and induces a transitive action on $(S^2)^r/(\Z/2)$. Therefore
\[
(S^2)^r/(\Z/2)\cong \Spin(3)^r/H,
\]
where $H$ fixes the orbit $\{\pm(e_3^1,\dots,e_3^r)\}$. Since $\pi_1((S^2)^r/(\Z/2))\cong \Z/2$ and $\Spin(3)^r$ is
simply connected, $H$ has two connected components.

Set $s:=(e_1^1e_3^1,\dots,e_1^re_3^r)$. A Clifford computation gives
\[
s\,(e_3^1,\dots,e_3^r)\,s^{-1}=(-e_3^1,\dots,-e_3^r),
\]
so $s\in H$, and $H=\langle \Spin(2)^r,s\rangle$.

Let $\sigma:H\longrightarrow \Ogroup(2r)$ be the isotropy representation. One checks that $s$ acts by $-1$ on each $e_1^j$-direction
and by $+1$ on each $e_2^j$-direction, hence
\[
\sigma(s)=\diag(-1,1,-1,1,\dots,-1,1),
\qquad \det(\sigma(s))=(-1)^r.
\]
Thus the quotient is orientable if and only if  $r$ is even.

Assume $r$ is even and write $r=2l$. Then $\sigma(H)\subset SO(2r)$. Since $\Spin(3)^r$ is simply connected,
$\Spin(3)^r/H$ is spin if and only if  $\sigma$ lifts to $\tilde\sigma:H\to \Spin(2r)$. On the identity component
$\Spin(2)^r\le H$, the restriction $\sigma|_{\Spin(2)^r}$ lifts via the standard group homomorphism
$\Spin(2)^r\longrightarrow \Spin(2r)$, so it remains to lift $s$.

Let $f_1,\dots,f_{2r}$ be the standard basis of $\R^{2r}\subset Cl_{2r}$. The preimages of $\sigma(s)$ under
$\lambda_{2r}$ are $\pm f_1f_3\cdots f_{2r-1}$. Conjugation by $s$ inverts each $\Spin(2)$-factor,
and conjugation by $\pm f_1f_3\cdots f_{2r-1}$
induces the same inversion on its image in $\Spin(2r)$. Moreover $s^2=(-1,\dots,-1)\in \Spin(2)^r$, whose image in
$\Spin(2r)$ equals $(-1)^r=1$ (since $r$ is even). Hence a lift exists if and only if 
$(\pm f_1f_3\cdots f_{2r-1})^2=1$. Using $(f_1f_3\cdots f_{2r-1})^2=(-1)^l$, this holds if and only if  $l$ is even, or equivalently, $r=4k$.
\end{proof}

\medskip

We are ready to construct the example.

\begin{proof}[Proof of \cref{disconnected isotropy counterexample}]

Let $M=\Spin(3)^5/K$, where $K$ is generated by $\Spin(2)^5$ and
\[
a_1=(e_1^1e_3^1,\dots,e_1^4e_3^4,1),\qquad
a_2=(1,e_1^2e_3^2,\dots,e_1^5e_3^5).
\]
Then $\Spin(2)^5\triangleleft K$ and $K/\Spin(2)^5\cong \Z/2\oplus\Z/2$, hence there is a regular
$\Z/2\oplus\Z/2$ covering $(S^2)^5\to M$.

Let $H_1=\langle \Spin(2)^5,a_1\rangle$, $H_2=\langle \Spin(2)^5,a_2\rangle$ and
$H_{12}=\langle \Spin(2)^5,a_1a_2\rangle$, where we have $a_1a_2=(e_1^1e_3^1,-1,-1,-1,e_1^5e_3^5)$, and set
$M_i=\Spin(3)^5/H_i$ for $i=1,2,12$. By \cref{spin structures S2 product}, $M_1$ and $M_2$ are spin, whereas
$M_{12}$ is not spin (here $-1\in\Spin(3)$ acts trivially on $S^2$ under $x\mapsto axa^{-1}$, so $a_1a_2$ flips only
two $S^2$-factors). In particular, $M$ is not spin by \cref{spin structures orbit spaces discrete}.

By the long exact sequence of homotopy and \cref{frame bundle discrete group}, the covering $M_1\to M$ induces a commutative diagram
\[
\begin{tikzcd}
& & 1 \ar{d}{}& 1 \ar{d}{}& \\
1 \ar{r}{}& \Z/2 \ar{r}{}\ar{d}{\cong}& \pi_1(FM_1)\cong\Z/2\oplus\Z/2 \ar{r}{}\ar{d}{}&
\pi_1(M_1)\cong \Z/2\ar{r}{}\ar{d}{}& 1\\
1 \ar{r}{}& \Z/2\ar{r}{i_*}& \pi_1(FM)\ar{d}{} \ar{r}{}&
\pi_1(M)\cong\Z/2\oplus \Z/2\ar{r}{}\ar{d}{}& 1\\
& & \Z/2\ar{r}{\cong}\ar{d}{}& \Z/2\ar{d}{} &  \\
& &1& 1 &
\end{tikzcd}
\]
where $\pi_1(FM_1)\cong(\Z/2)^2$ since $M_1$ is spin.

Now $\pi_1(FM)$ has order $8$ and contains $(\Z/2)^2$ as a normal subgroup, hence it is not $\Z/8$ nor the quaternion group.
Since $M$ is not spin, $i_*:\Z/2\to \pi_1(FM)$ has no section, so $\pi_1(FM)\not\cong (\Z/2)^3$. Thus
$\pi_1(FM)\cong D_4$ or $\pi_1(FM)\cong\Z/4\oplus \Z/2$ with $i_*(1)=(2,0)$.
We can rule out $\Z/4\oplus\Z/2$. Indeed, restricting the
extension
\[
1\to \Z/2\xrightarrow{i_*}\pi_1(FM)\to \pi_1(M)\cong(\Z/2)^2\to 1
\]
to the subgroups $\langle(1,0)\rangle$ and $\langle(1,1)\rangle$ of $\pi_1(M)$ yields the short exact sequence
\[
1\to \Z/2\xrightarrow{i_*}\Z/4\to \Z/2\to 1.
\]
Consequently, the two inequivalent $2$-coverings of $M$ corresponding to the subgroups $\langle(1,0)\rangle$ and $\langle(1,1)\rangle$ are not spin.
This contradicts that the only $2$-covering whose total space is not spin is $M_{12}\longrightarrow M$ is not spin. Hence
$\pi_1(FM)\cong D_4$.

Fix $D_4=\langle a,b\mid a^4=b^2=e,\ bab=a^3\rangle$. Then the centre is $\langle a^2\rangle$. Let $\bar a,\bar b$ be the
images of $a,b$ in $\pi_1(M)\cong(\Z/2)^2$, and let $\overline{M}\to M$ be the $2$-cover corresponding to
$\langle \bar b\rangle\le \pi_1(M)$. Then
\[
\begin{tikzcd}
& & 1 \ar{d}{}& 1 \ar{d}{}& \\
1 \ar{r}{}& \Z/2\cong\langle a^2\rangle \ar{r}{}\ar{d}{\cong}&
\pi_1(F\overline{M})\cong \langle a^2,b\rangle \ar{r}{}\ar{d}{}&
\pi_1(\overline{M})\cong \langle \bar b\rangle\ar{r}{}\ar{d}{}& 1\\
1 \ar{r}{}& \Z/2\cong\langle a^2\rangle\ar{r}{i_*}& \pi_1(FM)\cong D_4\ar{r}{}\ar{d}{}&
\pi_1(M)\cong \langle \bar a,\bar b\rangle\ar{r}{}\ar{d}{}& 1\\
& & \Z/2\cong\langle \bar a\rangle\ar{r}{\cong}\ar{d}{}&
\Z/2\cong\langle \bar a\rangle\ar{d}{} &  \\
& &1& 1 &
\end{tikzcd}
\]
and $\pi_1(F\overline{M})=\langle a^2,b\rangle\cong(\Z/2)^2$, so $\overline{M}$ is spin.

Spin structures on $\overline{M}$ correspond to homomorphisms $f\in\Hom(\pi_1(F\overline{M}),\Z/2)\cong H^1(F\overline{M},\Z/2)$ with $f\circ i_*=\mathrm{id}_{\Z/2}$. Under the identification $\pi_1(F\overline{M})=\langle a^2,b\rangle$, the last condition can be written as $f(a^2)=1$. Hence there are exactly two spin structures $f_i:\pi_1(F\overline{M})\longrightarrow\Z/2$, distinguished by $f_i(b)=i$ for $i=0,1$. The involution associated to $\overline{M}\to M$ acts by conjugation by $a$ on $\langle a^2,b\rangle< D_4$. Since $C_a(b)=aba^3=a^2b$, it swaps these two homomorphisms. Thus, this involution permutes the two spin structures, $f_0$ and $f_1$, of $\overline{M}$.

Finally, $\overline{M}\to M$ is equivalent to either $M_1\to M$ or $M_2\to M$; assume it is $M_1\to M$ (if not, $M_1\longrightarrow M$ would be equivalent to the covering induced by the subgroup $\langle\overline{a}\overline{b}\rangle$ and analogous computations would show that the $\Z/2$ action induced by the covering also permutes both spin structures). 
Then $a_2\in K\le N_{\Spin(3)^5}(H_1)$
defines an involution
\[
f_{a_2}:\Spin(3)^5/H_1\to\Spin(3)^5/H_1,
\qquad f_{a_2}(gH_1)=ga_2H_1,
\]
which exchanges the two spin structures of $M_1$.
Moreover, $f_{a_2}$ preserves orientation,
since it is induced by the antipodal map
on four $S^2$-factors.

Consequently, the mapping torus $(M_1)_{f_{a_2}}$
is orientable but not spin by \cref{mapping torus}.
\end{proof}
\section{Spin structures on cohomogeneity one manifolds of type (4)}\label{sec: spin structures on cohomogeneity one manifolds of type (4)}
The aim of this section is to study the existence of
(invariant) spin structures on cohomogeneity one
manifolds of type (4).

When $l=1$, let $q:K\to S^1\cong K/H$ be the
homomorphism defining the action of $K$ on $\R^2$,
and let $j^*$ denote the map on
$\Z/2$-valued characters induced by $j_*:\pi_1(K)\to\pi_1(G)$.
\begin{theorem}\label{main_theorem noncompact section}
Assume that $(H,K,G,i,j,l)$ is a cohomogeneity one tuple
describing a manifold $M$ of type (4), and that $G$ and $K$
are connected Lie groups. Then:
\begin{itemize}
    \item[(1)] If $l\geq 2$, then $M$ is spin if and only if $G/H$ is spin.
    \item[(2)] If $l=1$ and $|\pi_0(H)|$ is even, then $M$ is spin if and only if $G/K$ is spin.
    \item[(3)] If $l=1$ and $|\pi_0(H)|$ is odd, then:
\begin{itemize}
    \item[(a)] If $\delta_2\circ q_*\in\operatorname{Im}(j^*)$,
    then $M$ is spin if and only if $G/K$ is spin.
    \item[(b)] If $\delta_2\circ q_*\notin\operatorname{Im}(j^*)$,
    then at most one of $M$ and $G/K$ is spin.
\end{itemize}
\end{itemize}
\end{theorem}


\begin{remark}
If $l\geq 2$, then $H$ is connected if and only if
$K$ is connected, by the homotopy exact sequence
of $H\to K\to S^l$.
\end{remark}

\begin{proof}
Recall that in case \textup{(4)} we have
\[
M\cong G\times_K\R^{l+1}
\]
by \cref{cohomogenity one structure}. Let
\[
\pi:G\times \R^{l+1}\longrightarrow G\times_K\R^{l+1}\cong M
\]
denote the quotient map. Since the right action of $K$ on $G\times \R^{l+1}$ is free, we can apply
\cref{spin structures orbit spaces} to this action.

We first describe the two bundles appearing in that criterion. Since both $G$ and $\R^{l+1}$ are parallelizable, the product $G\times \R^{l+1}$ is parallelizable. More precisely,
if $n:=\dim G$, then
\[
T(G\times \R^{l+1})\cong (G\times \R^{l+1})\times (\mathfrak g\oplus \R^{l+1}),
\]
and hence its oriented orthonormal frame bundle is canonically
\[
F(G\times \R^{l+1})\cong G\times \R^{l+1}\times \SO(n+l+1).
\]

The vertical bundle $\tau_K:=\ker(d\pi)$ is the bundle tangent to the $K$-orbits. As in the homogeneous-space case,
it is canonically isomorphic to the trivial bundle
\[
(G\times \R^{l+1})\times \mathfrak k \longrightarrow G\times \R^{l+1},
\]
where the right action of $K$ is given by
\[
((g,v),X)\cdot k=((gk,k^{-1}v),\Ad(k^{-1})X).
\]
Consequently,
\[
F(\tau_K)\cong (G\times \R^{l+1})\times \SO(\mathfrak k).
\]
Writing $k=\dim K$, the induced map on fundamental groups
coming from the action of $K$ on the frame bundle of $\tau_K$ is
\[
(ev^{K\acts G\times \R^{l+1}}_*,(\Ad_K)_*):
\pi_1(K)\longrightarrow
\pi_1(G\times \R^{l+1})\times\pi_1(\SO(k))
\cong \pi_1(G)\times\pi_1(\SO(k)).
\]
We now compute the corresponding map for $T(G\times \R^{l+1})$. Choose an $\Ad(K)$-invariant inner product on $\mathfrak g$
and let $\mathfrak m=\mathfrak k^\perp$. Then
$\mathfrak g=\mathfrak k\oplus\mathfrak m$ is a reductive
decomposition, with $\mathfrak m\cong T_{eK}(G/K)$.
Consequently, we have the orthogonal decomposition
\[
\mathfrak g\oplus\R^{l+1}
=\mathfrak k\oplus\mathfrak m\oplus\R^{l+1},
\]
and the action of $K$ decomposes as the direct sum of
the adjoint action on $\mathfrak k$, the isotropy
representation of $G/K$ on $\mathfrak m$, and the given
orthogonal action of $K$ on $\R^{l+1}$.

Therefore, after composing with the appropriate $\delta$-maps as in \cref{Classification via FG}, the induced map on $\pi_1(K)$ is
\[
(ev_*^{K\acts G\times \R^{l+1}},\,\Ad_{K*}+\sigma_{G/K*}+ev_*^{K\acts \R^{l+1}}):
\pi_1(K)\longrightarrow \pi_1(G)\times \Z/2.
\]

Now apply \cref{spin structures orbit spaces}. Since $M\cong (G\times \R^{l+1})/K$, the manifold $M$ is spin if and only if there exist spin structures
\[
f:\pi_1(G)\times \Z/2\longrightarrow \Z/2
\qquad\text{and}\qquad
f_\tau:
\pi_1(G)\times\pi_1(\SO(k))
\longrightarrow\Z/2
\]
on $T(G\times \R^{l+1})$ and $\tau_K$, respectively, such that
\[
f\circ(ev_*^{K\acts G\times \R^{l+1}},\,\Ad_{K*}+\sigma_{G/K*}+ev_*^{K\acts \R^{l+1}})
=
f_\tau\circ
(ev_*^{K\acts G\times\R^{l+1}},(\Ad_K)_*).
\]

Let $\phi$ and $\phi_\tau$ denote the restrictions of
$f$ and $f_\tau$, respectively, to the $\pi_1(G)$-factor.
As in the proof of \cref{spin structures homogeneous spaces}, the above identity is equivalent to
\[
\phi\circ j_*+\Ad_{K*}+\sigma_{G/K*}+ev_*^{K\acts \R^{l+1}}
=
\phi_\tau\circ j_*+\Ad_{K*}.
\]
After canceling the $\Ad_{K*}$ terms and setting $\psi:=\phi+\phi_\tau\in \Hom(\pi_1(G),\Z/2)$,
we obtain:

\[
M \text{ is spin }
\Longleftrightarrow
\exists\,\psi:\pi_1(G)\longrightarrow \Z/2
\text{ such that }
\psi\circ j_*=\sigma_{G/K*}+ev_*^{K\acts \R^{l+1}}.
\]

Thus, it remains to study the term $\sigma_{G/K*}+ev_*^{K\acts \R^{l+1}}$. The next lemma describes the relation between the evaluation map $ev_*^{K\acts \R^{l+1}}$ and the isotropy representation of the transitive action of $K$ on $S^l$ for $l\geq 2$.

\begin{lemma}\label{invariant spin structures on spheres and euclidean spaces}
Assume $l\geq 2$ and that $K$ is a compact connected Lie group acting transitively and by orientation-preserving isometries on
\[
S^l=\{x\in \R^{l+1}:\Vert x\Vert^2=1\}.
\]
Let $H$ be the isotropy subgroup of some point in $S^l$, let $i:H\longrightarrow K$ denote the inclusion and let
$\sigma:H\longrightarrow \SO(l)$ denote the isotropy representation. Then
\[
ev_*^{K\acts \R^{l+1}}\circ i_*=\sigma_*.
\]
\end{lemma}

\begin{proof}
Recall that $T_xS^l=\langle x\rangle^\perp\subset\R^{l+1}$. The action of $K$ on $S^l$ induces a group morphism
\[
\iota:K\longrightarrow \Isom_+(S^l)=\SO(l+1).
\]
Moreover, the injective map
\[
\phi:FS^l\longrightarrow F\R^{l+1},\qquad
\phi(x,\{v_1,\dots,v_l\})=(x,\{x,v_1,\dots,v_l\}),
\]
is $K$-equivariant. Fix a point $x\in S^l$ and let $H=K_x$. Then the elements of $\iota(H)$ are matrices in
$\SO(l+1)$ fixing $x$, and hence they preserve $\langle x\rangle^\perp$. Therefore, for $h\in H$,
\[
h\cdot\phi(x,\{v_1,\dots,v_l\})
=(x,\{x,\iota(h)v_1,\dots,\iota(h)v_l\}).
\]
Thus, under the identification of the fibre of $FS^l$ over $x$ with $\SO(l)$, the action of $H$ on this fibre is precisely given by the isotropy representation
$\sigma:H\to\SO(l)$. Passing to fundamental groups and composing with the corresponding $\delta$-maps gives
\[
ev_*^{K\acts \R^{l+1}}\circ i_*=\sigma_*,
\]
since $\delta_{l+1}$ composed with the map induced by the standard inclusion
$\SO(l)\hookrightarrow\SO(l+1)$ equals $\delta_l$.
\end{proof}

\begin{remark}
For $l\geq 2$, the spin structure on $S^l$ is $K$-invariant if and only if the spin structure on $\R^{l+1}$ is
$K$-invariant. For $l=1$, the connected spin structure of $S^1$ is $K$-invariant if and only if the spin structure
on $\R^2$ is $K$-invariant.
\end{remark}

Therefore, we have
\[
(\sigma_{G/K*}+ev_*^{K\acts \R^{l+1}})\circ i_*=\sigma_{G/K*}\circ i_*+\sigma_{K/H*}.
\]

In the next lemma we describe the relation between the isotropy representations of $G/H$, $G/K$ and $K/H$.

\begin{lemma}\label{hf isotoropy representation}
Let $L\subset H\subset G$ be compact connected Lie groups and let $i:L\longrightarrow H$ denote the inclusion. Then
the isotropy representation of $G/L$ is
\[
\sigma_{G/L}=\sigma_{H/L}\oplus\sigma_{G/H}\circ i.
\]
\end{lemma}

\begin{proof}
Since $L\subset H\subset G$ are compact, the homogeneous spaces are reductive. Thus we have the following
decompositions:
\begin{itemize}
    \item for $G/H$, we have $\mathfrak{g}=\mathfrak{h}\oplus \mathfrak{m}$, where $\mathfrak{m}$ is
    $\Ad(H)$-invariant, and hence also $\Ad(L)$-invariant;
    \item for $H/L$, we have $\mathfrak{h}=\mathfrak{l}\oplus \mathfrak{n}$, where $\mathfrak{n}$ is
    $\Ad(L)$-invariant;
    \item for $G/L$, we have the reductive decomposition
    \[
    \mathfrak{g}=\mathfrak{l}\oplus (\mathfrak{m}\oplus\mathfrak{n}).
    \]
\end{itemize}
Therefore,
\[
\sigma_{G/L}=\Ad(L)_{|\mathfrak{m}\oplus\mathfrak{n}}
=\Ad(L)_{|\mathfrak{m}}\oplus \Ad(L)_{|\mathfrak{n}},
\]
since both subspaces are $\Ad(L)$-invariant. The first component is the restriction of $\sigma_{G/H}$ to the
subgroup $L\subset H$, hence it is $\sigma_{G/H}\circ i$. The second component is precisely the isotropy
representation $\sigma_{H/L}$. This proves the claim.
\end{proof}

Assume now that $l\geq 2$. Then $G$, $H$ and $K$ are connected and
$i_*:\pi_1(H)\longrightarrow \pi_1(K)$ is surjective. Therefore $M$ admits a spin structure if and only if there
exists $\psi:\pi_1(G)\longrightarrow \Z/2$ satisfying
\[
\psi\circ (j\circ i)_*
= \sigma_{G/K*}\circ i_*+ev_*^{K\acts \R^{l+1}}\circ i_*
= \sigma_{G/H*},
\]
by \cref{invariant spin structures on spheres and euclidean spaces} and
\cref{hf isotoropy representation}. This condition holds if and only if $G/H$ is spin, by
\cref{spin structures homogeneous spaces}. This completes the proof in the case $l\geq 2$.

Assume now that $l=1$. As above, $M$ admits a spin structure if and only if there exists a map
$\psi:\pi_1(G)\longrightarrow \Z/2$ satisfying
\[
\psi\circ j_*=\sigma_{G/K*}+ev_*^{K\acts \R^{2}}.
\]
The fibration
\[
\begin{tikzcd}
    H\ar{r}{} & K\ar{r}{q} & S^1
\end{tikzcd}
\]
together with the long exact sequence of homotopy induces an exact sequence
\[
\begin{tikzcd}
1 \ar{r}{} & \pi_1(H)\ar{r}{i_*} & \pi_1(K)\ar{r}{q_*} & \Z\ar{r}{} & \pi_0(H)\ar{r}{} & 1.
\end{tikzcd}
\]
Since the orthogonal action of $K$ on $\R^2$ is given by
$q:K\longrightarrow \SO(2)\cong S^1$, we have
\[
ev_*^{K\acts \R^{2}}=\delta_2\circ q_*.
\]
Thus, the above condition can be written as
\[
\psi\circ j_*=\sigma_{G/K*}+\delta_2\circ q_*.
\]

Recall that $G/K$ is spin if and only if there exists
$\psi':\pi_1(G)\longrightarrow \Z/2$ such that
\[
\psi'\circ j_*=\sigma_{G/K*}.
\]
If $|\pi_0(H)|$ is even, then $\delta_2\circ q_*=0$, and hence $M$ is spin if and only if $G/K$ is spin, by
\cref{spin structures homogeneous spaces}.

Suppose now that $|\pi_0(H)|$ is odd. By the above,
\[
M\text{ is spin}\quad\Longleftrightarrow\quad
\sigma_{G/K*}+\delta_2\circ q_*\in\operatorname{Im}(j^*),
\]
whereas
\[
G/K\text{ is spin}\quad\Longleftrightarrow\quad
\sigma_{G/K*}\in\operatorname{Im}(j^*).
\]
Since $\operatorname{Im}(j^*)$ is a subgroup of
$\Hom(\pi_1(K),\Z/2)$, if
$\delta_2\circ q_*\in\operatorname{Im}(j^*)$,
then these two conditions are equivalent.
Otherwise, they cannot both hold. This completes the proof.

\end{proof}

We have an analogous result for $G$-invariant spin structures.

\begin{theorem}\label{main_theorem noncompact invariant section}
Assume that $(H,K,G,i,j,l)$ is a cohomogeneity one tuple
describing a manifold $M$ of type (4), and that $G$ and $K$
are connected Lie groups. Then:
\begin{itemize}
    \item[(1)] If $l\geq 2$, then $M$ admits a $G$-invariant spin structure if and only if $G/H$ admits a $G$-invariant spin structure.
    \item[(2)] If $l=1$ and $|\pi_0(H)|$ is even, then $M$ admits a $G$-invariant spin structure if and only if $G/K$ admits a $G$-invariant spin structure.
    \item[(3)] If $l=1$ and $|\pi_0(H)|$ is odd, then $M$ and $G/K$ cannot both admit $G$-invariant spin structures. In particular, if one admits a $G$-invariant spin structure, then the other does not.
\end{itemize}
\end{theorem}

\begin{proof}
The proof follows the same arguments used in \cref{main_theorem noncompact section} and \cref{spin structures homogeneous spaces}. By \cref{cohomogenity one structure}, we have a $G$-equivariant diffeomorphism $M\cong G\times_K\R^{l+1}$, where the action of $G$ on $G\times_K \R^{l+1}$ satisfies $g[g',x]=[gg',x]$ for all $g\in G$ and $[g',x]\in G\times_K \R^{l+1}$. This action lifts to a free action of $G$ on $G\times\R^{l+1}$ given by $g(g',x)=(gg',x)$ for all $g\in G$ and $(g',x)\in G\times \R^{l+1}$. Furthermore, this action lifts to the oriented frame bundle $F(G\times\R^{l+1})\cong G\times \R^{l+1}\times\SO(n+l+1)$, where $n=\dim G$. This action satisfies $g(g',x,A)=(gg',x,A)$ for all $g\in G$ and $(g',x,A)\in F(G\times\R^{l+1})$. Thus $ev_*^{G\acts F(G\times \R^{l+1})}:\pi_1(G)\longrightarrow \pi_1(G)\times \Z/2$ is the natural inclusion.

Since $K$ is connected, by \cref{main_theorem noncompact section}, $M$ is spin if and only if there exists $\psi:\pi_1(G)\longrightarrow \Z/2$ such that $\psi\circ j_*=\sigma_{G/K*}+ev_*^{K\acts \R^{l+1}}$. By \cref{invariant spin structure pullback invariant} and \cref{invariant spin structures orbit spaces}, $M$ admits a $G$-invariant spin structure if and only if there exist spin structures on $T(G\times\R^{l+1})$ and on $\tau_K$ satisfying the compatibility condition for the $K$-action and the corresponding invariance condition for the $G$-action. Using the same arguments as in the second part of \cref{spin structures homogeneous spaces}, the latter condition is equivalent to $\phi=\phi_\tau$, in the notation of the proof of \cref{main_theorem noncompact section}. Consequently, $\psi=\phi+\phi_\tau=0$, and the compatibility condition becomes
\[
\sigma_{G/K*}+ev_*^{K\acts \R^{l+1}}=0.
\]
Therefore, $M$ admits a $G$-invariant spin structure if and only if this equality holds.

Assume that $l\geq 2$. Then $G$, $H$ and $K$ are connected and $i_*:\pi_1(H)\longrightarrow \pi_1(K)$ is surjective. Therefore, $\sigma_{G/K*}+ev_*^{K\acts \R^{l+1}}=0$ if and only if $\sigma_{G/K*}\circ i_*+ev_*^{K\acts \R^{l+1}}\circ i_*=0$. As in the non-invariant case, we can use \cref{invariant spin structures on spheres and euclidean spaces} and \cref{hf isotoropy representation} to conclude that
\[
\sigma_{G/K*}\circ i_*+ev_*^{K\acts \R^{l+1}}\circ i_*=\sigma_{G/H*}.
\]
Thus, $M$ has a $G$-invariant spin structure if and only if $G/H$ has a $G$-invariant spin structure, by \cref{spin structures homogeneous spaces}.

Assume now that $l=1$. As above, $M$ admits a $G$-invariant spin structure if and only if $\sigma_{G/K*}+ev_*^{K\acts \R^2}=0$. The fibration \(\begin{tikzcd}
    H\ar{r}{} & K\ar{r}{q} & S^1
\end{tikzcd}
\) together with the long exact sequence of homotopy induces an exact sequence
\[
\begin{tikzcd}
\pi_1(H)\ar{r}{i_*} & \pi_1(K)\ar{r}{q_*} & \Z\ar{r}{} & \pi_0(H)\ar{r}{} & 1.
\end{tikzcd}
\]
Since the orthogonal action of $K$ on $\R^2$ is given by $q:K\longrightarrow \SO(2)\cong S^1$, we have
\[
ev_*^{K\acts \R^2}=\delta_2\circ q_*.
\]
Thus, $M$ has a $G$-invariant spin structure if and only if $\sigma_{G/K*}+\delta_2\circ q_*=0$.

If $|\pi_0(H)|$ is even, then $q_*(\pi_1(K))\subseteq 2\Z$ and therefore $\delta_2\circ q_*=0$. Consequently, $M$ has a $G$-invariant spin structure if and only if $G/K$ has a $G$-invariant spin structure.

If $|\pi_0(H)|$ is odd, then $\delta_2\circ q_*\neq0$. Consequently, $M$ and $G/K$ cannot both have a $G$-invariant spin structure. In particular, if one of them has a $G$-invariant spin structure, then the other does not.
\end{proof}

\section{Spin structures on cohomogeneity one manifolds of type (2)}\label{sec: spin structures on cohomogeneity one manifolds of type (2)}

In this section, we study the existence of (invariant) spin structures
on cohomogeneity one manifolds of type (2). Recall that, with the
notation introduced in \cref{sec: preliminaris cohomogeneity one topology},
if a compact connected Lie group $G$ acts on a manifold $M$ with
cohomogeneity one of type (2), then there exist subgroups
$H$, $K_+$ and $K_-$ fitting into the diagram in
\cref{type (2) commutative diagram} such that $M$ is
$G$-equivariantly diffeomorphic to $M_+\cup_{G/H}M_-$,
where $M_\pm=G\times_{K_\pm}\D^{l_\pm+1}$.

The sets of spin structures and $G$-invariant spin structures
on $M_\pm$ will be denoted by $\Spin(M_\pm)$ and $\Spin_G(M_\pm)$,
respectively. Recall that we denote the inclusions by
$i_\pm:H\longrightarrow K_\pm$ and
$j_\pm:K_\pm\longrightarrow G$.
When $l_\pm=1$ and $K_\pm$ is connected, we have a surjective
group morphism $q_\pm:K_\pm\longrightarrow S^1$.

The manifolds $M_\pm$ have boundary $G$-equivariantly
diffeomorphic to $G/H$. By attaching collars
$\partial M_\pm\times[0,\infty)$ to $M_\pm$, we obtain manifolds
$M_\pm\cup_{\partial M_\pm}(\partial M_\pm\times[0,\infty))$
admitting a cohomogeneity one action of type (4).
Indeed, these manifolds are $G$-equivariantly diffeomorphic
to $G\times_{K_\pm}\R^{l_\pm+1}$, which are also
$G$-equivariantly diffeomorphic to the interiors $M_\pm^o$.
Thus, restriction induces bijections between spin structures
on $M_\pm$ and on $G\times_{K_\pm}\R^{l_\pm+1}$,
and similarly for $G$-invariant spin structures.

Assuming that $K_\pm$ are connected, by
\cref{main_theorem noncompact section}, we have
\[
\Spin(M_\pm)=
\{\psi\in\Hom(\pi_1(G),\Z/2):
\psi\circ j_{\pm*}
=\sigma_{G/K_\pm*}
+ev_*^{K_\pm\acts\R^{l_\pm+1}}\}.
\]
Similarly, by \cref{main_theorem noncompact invariant section},
$M_\pm$ admits a $G$-invariant spin structure if and only if
\[
\sigma_{G/K_\pm*}
+ev_*^{K_\pm\acts\R^{l_\pm+1}}=0.
\]
In this case, the $G$-invariant spin structure is unique.

If $M$ is spin, then $M_\pm$ are also spin. Thus, we need to study when spin structures on $M_+$ and $M_-$ can be glued to obtain a spin structure on $M$. Since $M_+$ and $M_-$ are manifolds with boundary, we will use the following result to glue spin structures.

\begin{lemma}\cite[Lemma 3]{massuyeau2003spin}\label{spin structure surgery}
Let $M_1$ and $M_2$ be spin manifolds of dimension $n$ with connected boundaries $S_1$ and $S_2$ and let $f:S_2\longrightarrow S_1$ be an orientation-reversing diffeomorphism. Define the set $J=\{(x,y)\in \Spin(M_1)\times \Spin(M_2):f^*(x_{|S_1})=y_{|S_2}\}$ and assume that it is not empty. Then $M=M_1\cup_f M_2$ is spin and the map $r:\Spin(M)\longrightarrow \Spin(M_1)\times \Spin(M_2)$ such that $r(x)=(x_{|M_1},x_{|M_2})$ is injective and its image is $J$.
\end{lemma}

\begin{proposition}\label{gluing spin structures}
Assume that $M$ is a closed connected manifold with a cohomogeneity one action of type (2) of a compact connected  Lie group $G$ where $K_+$ and $K_-$ are connected. Assume that:
\begin{itemize}
    \item[1.] either $l_+$ or $l_-$ are equal or greater than $2$ or
    \item[2.] $l_+=l_-=1$ and for either $+$ or $-$ we have $\delta_2\circ q_{\pm*}\notin \im j^*_\pm$. 
\end{itemize}
Then, if $M_+$ and $M_-$ are spin then $M$ is spin.
\end{proposition}

We remark that there exist manifolds $M$ with a cohomogeneity one action of type (2) such that $M_+$ and $M_-$ have a spin structure, but that $M$ does not have a spin structure. We construct these examples at the end of this section. 

\begin{proof}[Proof of \cref{gluing spin structures}]
Note that each spin structure on $M_\pm$ induces a spin structure on the boundary $G/H$, thus we have maps
$\phi_\pm:\Spin(M_\pm)\longrightarrow \Spin(G/H)$.
By \cref{spin structure surgery}, $M$ is spin if and only if there exist
$x_\pm\in\Spin(M_\pm)$ such that
$\phi_+(x_+)=\phi_-(x_-)\in\Spin(G/H)$.
Our next goal is to explicitly describe the maps $\phi_\pm$.

\textbf{Part 1: $\phi_\pm$ are injective.}
Recall that in \cref{sec: preliminaris cohomogeneity one topology}
we defined the $G$-equivariant inclusions
\[
\iota_\pm:M_\pm'=G\times_{K_\pm}(K_\pm/H)
\longrightarrow M_\pm,
\]
and the $G$-equivariant diffeomorphisms
\[
\zeta_\pm:G/H\longrightarrow
G\times_{K_\pm}(K_\pm/H).
\]
The inclusions of the boundary $G/H$ into $M_\pm$ are precisely
$\iota_\pm\circ\zeta_\pm:G/H\longrightarrow M_\pm$.
Thus, the maps $\phi_\pm$ are the compositions
\[
\begin{tikzcd}
\Spin(M_\pm)\ar{r}{\iota_\pm^*}
& \Spin(M_\pm')\ar{r}{\zeta_\pm^*}
& \Spin(G/H).
\end{tikzcd}
\]

Since $\zeta_\pm$ are diffeomorphisms, the maps
$\zeta_\pm^*$ are bijective. Hence, to prove that
$\phi_\pm$ are injective, it is enough to prove that
$\iota_\pm^*$ are injective.

Recall that $\pi_1(M_\pm)\cong\pi_1(G/K_\pm)$,
$\pi_1(M_\pm')\cong\pi_1(G/H)$ and
$\iota_{\pm*}=\overline{i}_{\pm*}$, which is
surjective by \cref{iota=i}.
Consequently, the induced maps
\[
\iota_\pm^*:H^1(M_\pm;\Z/2)
\longrightarrow H^1(M_\pm';\Z/2)
\]
are injective. Since spin structures form an affine
space over $H^1(-;\Z/2)$, the restriction maps
$\iota_\pm^*:\Spin(M_\pm)\longrightarrow\Spin(M_\pm')$
are also injective. Therefore, $\phi_\pm$ are injective.

\textbf{Part 2: Study of $\phi_\pm$ when $l_+\geq2$ or $l_-\geq2$.}
Assume that $l_+\geq2$ (the case where $l_-\geq2$ is analogous).
Since $K_+$ is connected, so is $H$. Moreover, since
$K_+/H\cong S^{l_+}$ is simply connected, the long exact
sequence of homotopy implies that
$\iota_{+*}:\pi_1(M_+')\longrightarrow\pi_1(M_+)$
is an isomorphism. Consequently,
\[
\iota_+^*:H^1(M_+;\Z/2)\longrightarrow H^1(M_+';\Z/2)
\]
is an isomorphism. Since both $M_+$ and $M_+'$ are spin
and spin structures form an affine space over $H^1(-;\Z/2)$,
the map $\iota_+^*:\Spin(M_+)\longrightarrow\Spin(M_+')$
is a bijection, and hence so is $\phi_+$.
Thus, if $M_-$ is spin, then
\[
\phi_+(\Spin(M_+))\cap\phi_-(\Spin(M_-))\neq\emptyset.
\]
This completes the proof of \cref{gluing spin structures}
under assumption 1.

\textbf{Part 3: Study of $\iota_\pm^*$ when $l_\pm=1$.}
To study the case where $l_\pm=1$, we need to describe
$\iota_\pm^*$ and $\zeta_\pm^*$ explicitly. Recall that
\[
\Spin(M_\pm)=
\{\psi\in\Hom(\pi_1(G),\Z/2):
\psi\circ j_{\pm*}
=\sigma_{G/K_\pm*}+\delta_2\circ q_{\pm*}\}.
\]

Since $l_\pm=1$, we have
$M_\pm'\cong(G\times S^1)/K_\pm$, defined using the
inclusions
\[
k\longmapsto(j_\pm(k),q_\pm(k)^{-1}).
\]
The isotropy representations $\sigma'$ satisfy
$\sigma'_*=\sigma_{G/K_\pm*}$, by
\cref{l=1 homogeneous space}. Since $K_\pm$ are
connected, we have
\[
\Spin(M'_\pm)=
\{\psi\in\Hom(\pi_1(G)\times\Z,\Z/2):
\psi\circ(j_{\pm*},-q_{\pm*})
=\sigma_{G/K_\pm*}\}.
\]

Given $\psi\in\Hom(\pi_1(G)\times\Z,\Z/2)$,
we will write it as
$\psi=\psi_{\pi_1(G)}+\psi_\Z$,
where the subscripts denote the restrictions to each
component of the product $\pi_1(G)\times\Z$.

The inclusion maps
$\iota_\pm:M_\pm'\longrightarrow M_\pm$
lift to the inclusions
$G\times S^1\longrightarrow G\times\D^2$,
which are $K_\pm$-equivariant. Thus, we have
commutative diagrams
\[
\begin{tikzcd}
G\times S^1\ar{r}{}\ar{d}{}
& G\times\D^2\ar{d}{}\\
M_\pm'\ar{r}{\iota_\pm} & M_\pm.
\end{tikzcd}
\]

Using the outward normal to stabilise the frames along
the boundary, and applying the long exact sequence of
homotopy, we obtain commutative diagrams
\[
\begin{tikzcd}
\pi_1(K_\pm)
\ar{rr}{id_{K_\pm}}
\ar{d}{(j_{\pm*},-q_{\pm*},\sigma_{G/K_\pm*})}
&&
\pi_1(K_\pm)
\ar{d}{(j_{\pm*},
\delta_2\circ q_{\pm*}+\sigma_{G/K_\pm*})}\\
\pi_1(G)\times\Z\times\Z/2
\ar{rr}{(id_{\pi_1(G)},\delta_2+id_{\Z/2})}
\ar{d}{}
&&
\pi_1(G)\times\Z/2
\ar{d}{}\\
\pi_1(F((G\times S^1)/K_\pm))
\ar{rr}{\widetilde{\iota_{\pm*}}}
&&
\pi_1(F(G\times_{K_\pm}\D^2)).
\end{tikzcd}
\]

Given spin structures $\psi_\pm\in\Spin(M_\pm)$,
the maps
$\iota_\pm^*:\Spin(M_\pm)\longrightarrow\Spin(M_\pm')$
satisfy
\[
(\iota_\pm^*\psi_\pm)_{\pi_1(G)}=\psi_\pm
\qquad\text{and}\qquad
(\iota_\pm^*\psi_\pm)_\Z=\delta_2
\]
for any $\psi_\pm\in\Spin(M_\pm)$.

We are ready to prove that $M$ is spin under
assumption 2. The result will follow from the next
claim: $\iota_+^*:\Spin(M_+)\longrightarrow\Spin(M_+')$
is bijective if and only if
$\delta_2\circ q_{+*}\notin\im j_+^*$
(the case for $-$ is analogous).

We will prove the contrapositive statement, that
$\iota_+^*:\Spin(M_+)\longrightarrow\Spin(M_+')$
is not bijective if and only if
$\delta_2\circ q_{+*}\in\im j_+^*$.
Indeed, assume that $\iota_+^*$ is not bijective.
Since $\iota_+^*$ is injective and its image consists
of spin structures with $\Z$-component $\delta_2$,
there exists a spin structure with $\Z$-component zero.
Thus, there exists
$\psi'\in\Hom(\pi_1(G),\Z/2)$ satisfying
$\psi'\circ j_{+*}=\sigma_{G/K_+*}$.
Moreover, since $M_+$ is spin, there exists
$\psi\in\Hom(\pi_1(G),\Z/2)$ satisfying
\[
\psi\circ j_{+*}+\delta_2\circ q_{+*}
=\sigma_{G/K_+*}.
\]
Consequently,
$(\psi'-\psi)\circ j_{+*}=\delta_2\circ q_{+*}$,
and therefore
$\delta_2\circ q_{+*}\in\im j_+^*$.

Conversely, assume that
$\delta_2\circ q_{+*}\in\im j_+^*$,
and therefore there exists
$\eta\in\Hom(\pi_1(G),\Z/2)$ such that
$\eta\circ j_{+*}=\delta_2\circ q_{+*}$.
Since $M_+$ is spin, there exists
$\psi\in\Hom(\pi_1(G),\Z/2)$ satisfying
\[
\psi\circ j_{+*}+\delta_2\circ q_{+*}
=\sigma_{G/K_+*}.
\]
Then
$(\psi+\eta)_{\pi_1(G)}+0_\Z\in\Spin(M_+')$
is not in the image of $\iota_+^*$.

Thus, without loss of generality, assume that
the condition in assumption 2 holds for $+$.
We have $\delta_2\circ q_{+*}\notin\im j_+^*$
and hence $\phi_+$ is bijective.
Since $M_-$ is spin, we have
\[
\phi_+(\Spin(M_+))
\cap\phi_-(\Spin(M_-))\neq\emptyset,
\]
and thus $M$ is spin.
This completes the proof of \cref{gluing spin structures}
under assumption 2.
\end{proof}

\begin{corollary}
Assume that either $b_+$ or $b_-$ are odd. If $M_+$ and $M_-$ are spin then $M$ is spin.
\end{corollary}
\begin{proof}
Assume that $b_+$ is odd (the other case is analogous).
Since $\iota_{+*}:\pi_1(M_+')\to\pi_1(M_+)$
is surjective and its kernel has odd order,
every homomorphism $\pi_1(M_+')\to\Z/2$
factors uniquely through $\pi_1(M_+)$.
Consequently,
\[
\iota_+^*:H^1(M_+;\Z/2)
\longrightarrow H^1(M_+';\Z/2)
\]
is an isomorphism. Since both manifolds are spin,
$\iota_+^*:\Spin(M_+)\to\Spin(M_+')$
is bijective, and hence so is $\phi_+$.
The result follows from \cref{spin structure surgery}.
\end{proof}

Our next goal is to find criteria when $l_+=l_-=1$
and $\delta_2\circ q_{\pm*}\in\im j_\pm^*$.
We assume throughout that $M_+$ and $M_-$ are spin.
Let $\eta_\pm\in H^1(G,\Z/2)$ satisfy
$\delta_2\circ q_{\pm*}=\eta_\pm\circ j_{\pm*}$.
We claim that $\Spin(M_\pm')$ is the disjoint union of
\[
S_{\pm1}
=\{(\psi,\delta_2)\in\Spin(M_\pm'):
\psi\in\Spin(M_\pm)\}
=\iota_\pm^*(\Spin(M_\pm))
\]
and
\[
S_{\pm0}
=\{(\psi+\eta_\pm,0)\in\Spin(M_\pm'):
\psi\in\Spin(M_\pm)\}.
\]

Indeed, the same arguments as in the last part of
\cref{gluing spin structures} show that $S_{\pm0}$
and $S_{\pm1}$ are well-defined subsets, and they
are clearly disjoint. To prove that
$S_{\pm0}\sqcup S_{\pm1}=\Spin(M_\pm')$,
we consider
$(\psi_{\pi_1(G)},\psi_\Z)\in\Spin(M_\pm')$,
and consequently we have
$\psi_{\pi_1(G)}\circ j_{\pm*}
-\psi_\Z\circ q_{\pm*}=\sigma_{G/K_\pm*}$.
On the one hand, if $\psi_\Z=\delta_2$ then
$\psi_{\pi_1(G)}\circ j_{\pm*}
=\sigma_{G/K_\pm*}+\delta_2\circ q_{\pm*}$
and consequently
$\psi_{\pi_1(G)}\in\Spin(M_\pm)$.
Thus, $(\psi_{\pi_1(G)},\psi_\Z)\in S_{\pm1}$.
On the other hand, if $\psi_\Z=0$ then
$\psi_{\pi_1(G)}\circ j_{\pm*}=\sigma_{G/K_\pm*}$.
Then
$(\psi_{\pi_1(G)}-\eta_\pm)\circ j_{\pm*}
=\sigma_{G/K_\pm*}+\delta_2\circ q_{\pm*}$,
which implies that
$\psi_{\pi_1(G)}-\eta_\pm\in\Spin(M_\pm)$.
Thus, $(\psi_{\pi_1(G)},\psi_\Z)\in S_{\pm0}$.
Since
$|S_{\pm0}|=|S_{\pm1}|=|\Spin(M_\pm)|$,
we can conclude that
$|\Spin(M_\pm')|=2|\Spin(M_\pm)|$.

We note that the set $S_{\pm0}$ is independent
of the choice of $\eta_\pm$. Indeed, two spin
structures in $\Spin(M_\pm)$ differ by an element
of $\ker j_\pm^*$. Similarly, two elements in
$(j_\pm^*)^{-1}(\delta_2\circ q_{\pm*})$
also differ by an element of $\ker j_\pm^*$.
Therefore, given another
$\eta_\pm'\in(j_\pm^*)^{-1}(\delta_2\circ q_{\pm*})$,
we have
\[
(\psi+\eta_\pm',0)
=((\psi+\eta_\pm'-\eta_\pm)+\eta_\pm,0)
\in S_{\pm0},
\]
since
$\psi+\eta_\pm'-\eta_\pm\in\Spin(M_\pm)$.
In particular, since $K_\pm$ are connected,
we have
$H^1(G/K_\pm,\Z/2)\cong\ker j_\pm^*$
and hence
$|\Spin(M_\pm)|=|H^1(G/K_\pm,\Z/2)|$.

Our next goal is to describe the maps
$\zeta_\pm^*:\Spin(M_\pm')\longrightarrow\Spin(G/H)$.
The main difficulty is that $H$ may not be connected,
and hence we do not have a description of the set
of spin structures of $G/H$ using the isotropy
representation. To circumvent this problem, we
consider the regular covering
$\pi:G/H_0\longrightarrow G/H$, where $H_0$ is
the connected component of the identity of $H$.
Recall that $H/H_0\cong\Z/a$ and set
$I=j_+\circ i_+=j_-\circ i_-$.
Then
\[
\Spin(G/H_0)
=\{\psi\in\Hom(\pi_1(G),\Z/2):
\psi\circ I_*=\sigma_{G/H*}\},
\]
where we use that
$\sigma_{G/H_0}=(\sigma_{G/H})_{|H_0}$.

The $G$-equivariant diffeomorphisms
$\zeta_\pm:G/H\longrightarrow(G\times S^1)/K_\pm$
satisfy $\zeta_\pm(gH)=(g,e)K_\pm$.
We have a commutative diagram
\[
\begin{tikzcd}
H \ar{r}{i_\pm}\ar{d}{j_\pm\circ i_\pm}
& K_\pm\ar{d}{(j_\pm,q_\pm^{-1})}\\
G \ar{r}{(\id_G,e)}\ar{d}{}
& G\times S^1\ar{d}{}\\
G/H\ar{r}{\zeta_\pm}
& (G\times S^1)/K_\pm
\end{tikzcd}
\]
which induces a commutative diagram between the
homotopy groups of the frame bundles of $G/H$
and $(G\times S^1)/K_\pm$,
\[
\begin{tikzcd}
\pi_1(H)
\ar{r}{i_{\pm*}}
\ar{d}{((j_\pm\circ i_\pm)_*,\sigma_{G/H*})}
&
\pi_1(K_\pm)
\ar{d}{(j_{\pm*},-q_{\pm*},\sigma_{G/K_\pm*})}
\\
\pi_1(G)\times\Z/2
\ar{r}{(\id_{\pi_1(G)},0,\id_{\Z/2})}
\ar{d}{\rho_*}
&
\pi_1(G)\times\Z\times\Z/2
\ar{d}{\rho'_{\pm*}}
\\
\pi_1(F(G/H))
\ar{r}{\overline{\zeta_\pm}_*}
\ar{d}{\partial}
&
\pi_1(F((G\times S^1)/K_\pm))
\ar{d}{}
\\
\pi_0(H)\ar{r}{}\ar{d}{}
& 0
\\
0 &
\end{tikzcd}
\]
where
$\rho:G\times\SO(n-1)\longrightarrow F(G/H)$
and
$\rho'_\pm:G\times S^1\times\SO(n-1)
\longrightarrow F((G\times S^1)/K_\pm)$
are the orbit maps.

There are two important observations.
Firstly, since $M_\pm$ are spin, there exists
$\psi_\pm\in\Hom(\pi_1(G),\Z/2)$ such that
\[
\psi_\pm\circ j_{\pm*}
=\sigma_{G/K_\pm*}+\delta_2\circ q_{\pm*}.
\]
Consequently,
\[
(0,1,0)\notin
\im(j_{\pm*},-q_{\pm*},\sigma_{G/K_\pm*}),
\]
since $\psi_\pm+\delta_2+\id_{\Z/2}$ vanishes
on this image but takes the value $1$ on $(0,1,0)$.

Secondly, since $I_*=j_{\pm*}\circ i_{\pm*}$,
the homotopy exact sequence gives
\[
\pi_1(F(G/H_0))
\cong
\frac{\pi_1(G)\times\Z/2}
{\im(I_*,\sigma_{G/H*})}.
\]
Let
$\rho_0:G\times\SO(n-1)\longrightarrow F(G/H_0)$
denote the corresponding orbit map.
Then $\rho_*=\pi_*\circ\rho_{0*}$, and
$\rho_{0*}$ is surjective.
Since
\[
\overline{\zeta_\pm}_*\circ\rho_*
=\rho'_{\pm*}\circ
(\id_{\pi_1(G)},0,\id_{\Z/2}),
\]
we conclude that
\[
\pi^*\circ\zeta_\pm^*:
\Spin(M_\pm')\longrightarrow\Spin(G/H_0)
\]
is given by
\[
\psi=\psi_{\pi_1(G)}+\psi_\Z
\longmapsto\psi_{\pi_1(G)}.
\]
Consequently, $\pi^*\circ\phi_\pm$ is given by
\[
\psi\longmapsto\psi+\delta_2\longmapsto\psi.
\]
In particular, $\pi^*\circ\phi_\pm$ is injective.
Note that since
$M_\pm'\cong(G\times S^1)/K_\pm$
and $G/H_0$ are homogeneous spaces with connected
isotropy subgroups, the groups
$\pi_1(M_\pm')$ and $\pi_1(G/H_0)$ are abelian.

To study spin structures on $M$ in the remaining
cases, we need to distinguish between the case
where $a$ is odd and the case where $a$ is even,
since the obstruction to the existence of spin
structures will be different.

We first assume that $a$ is odd. We note that,
as in the proof of \cref{main_theorem (1) section},
the action of $H/H_0\cong\Z/a$ on $G/H_0$
fixes all the spin structures. Since $a$ is odd,
any central exact sequence
\[
\begin{tikzcd}
0\ar{r}{} & \Z/2\ar{r}{}
& \widetilde\Gamma\ar{r}{}
& \Z/a\ar{r}{} & 0
\end{tikzcd}
\]
splits. Consequently, all the spin structures
are $\Z/a$-invariant and
$\pi^*:\Spin(G/H)\longrightarrow
\Spin_{\Z/a}(G/H_0)=\Spin(G/H_0)$
is a bijection. Therefore,
$|\Spin(G/H_0)|=2|\Spin(M_\pm)|$
and
\[
\phi_+(\Spin(M_+))\cap\phi_-(\Spin(M_-))
\neq\emptyset
\]
if and only if
\[
\pi^*\phi_+(\Spin(M_+))
\cap\pi^*\phi_-(\Spin(M_-))
\neq\emptyset.
\]

We consider the group morphisms
\[
j'_{\pm*}:\pi_1(K_\pm)\longrightarrow
\pi_1(G)/2\pi_1(G),
\]
and the characters
\[
\chi_\pm:\im j'_{\pm*}\longrightarrow\Z/2
\]
defined by
\[
\chi_\pm(j'_{\pm*}x)
=\sigma_{G/K_\pm*}(x)+\delta_2\circ q_{\pm*}(x).
\]
These maps are well-defined. Indeed, consider
another $x'\in\pi_1(K_\pm)$ satisfying
$j'_{\pm*}(x)=j'_{\pm*}(x')$.
Since $M_\pm$ are spin, we can choose a spin
structure $\psi_\pm$ such that
\[
\psi_\pm\circ j_{\pm*}(x)
=\chi_\pm(j'_{\pm*}(x)).
\]
The equality
$j'_{\pm*}(x)=j'_{\pm*}(x')$
implies that
$\psi_\pm\circ j_{\pm*}(x)
=\psi_\pm\circ j_{\pm*}(x')$
and consequently,
$\chi_\pm(j'_{\pm*}x)
=\chi_\pm(j'_{\pm*}x')$.

\begin{proposition}\label{gluing spin structures a odd}
Assume that $M$ is a closed connected manifold with a cohomogeneity one action of type (2) of a compact connected  Lie group $G$ where $K_+$ and $K_-$ are connected. Assume that $l_+=l_-=1$, that $\delta_2\circ q_{\pm*}\in \im j^*_\pm$ and that $a$ is odd. If $M_+$ and $M_-$ are spin, then:
\begin{itemize}
    \item[1.] If $\ker j^*_+\neq \ker j^*_-$ then $M$ is spin.
    \item[2.] If $\ker j^*_+= \ker j^*_-$, then $\im j'_{+*}=\im j'_{-*}$, and $M$ is spin if and only if $\chi_+=\chi_-$.
\end{itemize}
\end{proposition}

\begin{proof}
By the above description of the maps $\pi^*\phi_\pm$, $M$ is spin if and only if there exist $\psi\in \Hom(\pi_1(G),\Z/2)$ satisfying $\psi\circ j_{+*}=\sigma_{G/K_+*}+\delta_2\circ q_{+*}$ and $\psi\circ j_{-*}=\sigma_{G/K_-*}+\delta_2\circ q_{-*}$.

Since $M_\pm$ are spin, $\Spin(M_\pm)$ are nonempty affine spaces with translation spaces $\ker j_\pm^*$. Moreover, the assumption $\delta_2\circ q_{\pm*}\in\im j_\pm^*$ implies that the restriction maps $\Spin(M_\pm)\to\Spin(M'_\pm)$ are not bijective. By the description of $\Spin(M_\pm')$, their images have half the cardinality of $\Spin(G/H)$. Since $a$ is odd, $\pi^*:\Spin(G/H)\longrightarrow \Spin(G/H_0)$ is a bijection. Consequently, $\Spin(M_\pm)$ are affine hyperplanes in the same affine space $\Spin(G/H_0)$. The description of $\pi^*\phi_\pm$ implies that the translation spaces of the hyperplanes $\pi^*\phi_\pm(\Spin(M_\pm))$ are precisely $\ker j_\pm^*$. Consequently, if $\ker j^*_+\neq \ker j^*_-$ then these hyperplanes have different translation spaces and therefore intersect. Hence $M$ is spin, proving the first part of the statement.

We now prove the second part. The first claim follows from the fact that 
$$\ker j_\pm^*=\operatorname{Ann}(\im j_{\pm*}')=\{f\in \Hom(\pi_1(G)/2\pi_1(G),\Z/2):f(j'_{\pm*}x)=0 \text{ for all $j'_{\pm*}x\in \im j_{\pm*}'$}\}.$$
Since $\pi_1(G)/2\pi_1(G)$ is a finite dimensional $\Z/2$-vector space, the equality $\ker j^*_+= \ker j^*_-$ implies $\im j'_{+*}=\im j'_{-*}$. This equality allows us to compare both characters $\chi_\pm$.

Set $U=\im j'_{+*}=\im j'_{-*}$. Since the two affine hyperplanes have the same translation space, they are either equal or disjoint. Every character $\psi:\pi_1(G)\to\Z/2$ factors uniquely through $\pi_1(G)/2\pi_1(G)$. Under this identification, the descriptions
of $\pi^*\phi_\pm$ give \[\pi^*\phi_\pm(\Spin(M_\pm))=\left\{\overline{\psi}\in\Hom(\pi_1(G)/2\pi_1(G),\Z/2):\overline{\psi}|_U=\chi_\pm\right\}.\]
If $\chi_+=\chi_-$, these two nonempty sets coincide, and hence $M$ is spin. If $\chi_+\neq\chi_-$, no character can restrict to both $\chi_+$ and $\chi_-$ on $U$. The two images are therefore disjoint, and $M$ is not spin.
\end{proof}

We focus now on the case where $a$ is even. In this setting, we have $\delta_2\circ q_{\pm*}=0$. Thus, $\delta_2\circ q_{\pm*}\in\im j^*_\pm$. Moreover, we can choose $\eta_\pm=0$ and therefore
$$ S_{\pm1}=\{(\psi,\delta_2)\in \Spin(M_\pm'):\psi\in \Spin(M_\pm)\}=\iota_\pm^*(\Spin(M_\pm))$$
and $$ S_{\pm0}=\{(\psi,0)\in \Spin(M_\pm'):\psi\in \Spin(M_\pm)\}.$$

Recall that $\pi^*\circ\zeta_\pm^*:\Spin(M_\pm')\longrightarrow\Spin(G/H_0)$ sends $(\psi,\psi_\Z)$ to $\psi$. By the above description of $S_{\pm0}$ and $S_{\pm1}$, we have
$$\Spin(M_\pm')=\{(\psi,t\delta_2):\psi\in\Spin(M_\pm),\ t\in\{0,1\}\}.$$
Consequently, identifying $\Spin(G/H_0)$ and $\Spin(M_\pm)$ with their descriptions as subsets of $\Hom(\pi_1(G),\Z/2)$, we obtain
$$\pi^*(\Spin(G/H))=\Spin(M_\pm).$$
In particular, we have
$$S=\Spin(M_+)=\Spin(M_-)$$
as subsets of $\Hom(\pi_1(G),\Z/2)$, where the equality refers to their corresponding characters.

We use the description $$\Spin(G/H)=\{f\in\Hom(\pi_1(F(G/H)),\Z/2):f\circ \alpha_{F(G/H)}=\id_{\Z/2}\},$$ where $\alpha_{F(G/H)}:\Z/2\longrightarrow\pi_1(F(G/H))$ is induced by the inclusion of a fiber. Given $\psi\in S$, we define $$f_\pm=\zeta_\pm^*(\psi,\delta_2).$$ We also define $\eta=\overline{\delta}_2\circ\partial\in\Hom(\pi_1(F(G/H)),\Z/2)$, where $\partial:\pi_1(F(G/H))\longrightarrow\pi_0(H)\cong\Z/a$ is the connecting homomorphism and $\overline{\delta}_2:\Z/a\longrightarrow\Z/2$ is reduction modulo $2$. Then $M$ is spin if and only if there exists $\psi\in S$ satisfying that $f_+=f_-$, or alternatively $M$ is not spin if and only if for all $\psi\in S$, we have $f_-=f_++\eta$.

To be able to detect this obstruction, we define elements
$$\gamma_\pm=(\overline{\zeta_\pm}_*)^{-1}\bigl(\rho'_{\pm*}(0,1,0)\bigr).$$
Since $q_{\pm*}(\pi_1(K_\pm))=a\Z$, the elements $\partial\gamma_+$ and $\partial\gamma_-$ are generators of $\Z/a$. Thus, we pick $r\in\Z$ satisfying $r\partial\gamma_+=\partial\gamma_-$ and an element $(\gamma,\epsilon)\in \pi_1(G)\times\Z/2$ satisfying $\rho_*(\gamma,\epsilon)=r\gamma_+-\gamma_-$.

\begin{lemma}
Let $f\in\Spin(G/H)$ and $t\in\Z/2$.
Then $f=\zeta_\pm^*(\psi,t\delta_2)$ for some
$\psi\in\Hom(\pi_1(G),\Z/2)$ if and only if
$f(\gamma_\pm)=t$.
\end{lemma}

\begin{proof}
Since $\zeta_\pm^*$ is a bijection, we can write
$f=\zeta_\pm^*(\psi,s\delta_2)$ for some
$\psi\in\Hom(\pi_1(G),\Z/2)$ and $s\in\Z/2$.
By the definition of $\gamma_\pm$,
\[
f(\gamma_\pm)
=(\psi,s\delta_2,\id_{\Z/2})(0,1,0)=s.
\]
Thus, $f(\gamma_\pm)=t$ if and only if $s=t$,
which proves the claim.
\end{proof}

\begin{proposition}\label{gluing spin structures a even}
Assume that $l_+=l_-=1$ and that $a$ is even. Let $(\gamma,\epsilon)$ be as above. If $M_+$ and $M_-$ are spin then $M$ is not spin if and only if there exists $x\in\pi_1(K_+)$ such that
\[j_{+*}(x)\equiv\gamma\pmod{2\pi_1(G)}\qquad\text{and}\qquad\sigma_{G/K_+*}(x)+\epsilon=1.\]
\end{proposition}

\begin{proof}
For $\psi\in S$, let $f_+=\zeta_+^*(\psi,\delta_2)$. The preceding lemma gives $f_+(\gamma_+)=1$ and shows that $f_+$ extends over $M_-$ precisely when $f_+(\gamma_-)=1$. Since $r$ is odd and $r\gamma_+-\gamma_-=\rho_*(\gamma,\epsilon)$, we have
\[
f_+(\gamma_-)
=r f_+(\gamma_+)-f_+(\rho_*(\gamma,\epsilon))
=1+\psi(\gamma)+\epsilon.
\]
As $\psi$ ranges over $S$, the structures $f_+$ exhaust the image of $\phi_+$. Hence $M$ is spin if and only if $\psi(\gamma)+\epsilon=0$ for some $\psi\in S$.
Suppose first that there exists $x\in\pi_1(K_+)$ such that $j_{+*}(x)\equiv\gamma\pmod{2\pi_1(G)}$. Every $\psi\in S$ vanishes on $2\pi_1(G)$ and satisfies $\psi\circ j_{+*}=\sigma_{G/K_+*}$, since $a$ is even. Therefore
\[
\psi(\gamma)=\psi(j_{+*}(x))=\sigma_{G/K_+*}(x)
\qquad\text{for every }\psi\in S.
\]
In this case, $M$ is not spin if and only if $\sigma_{G/K_+*}(x)+\epsilon=1$.
Suppose now that no $x\in\pi_1(K_+)$ satisfies the congruence. Then the class of $\gamma$ in $\pi_1(G)/2\pi_1(G)$ lies outside $\im j'_{+*}$. We can extend a basis of $\im j'_{+*}$ by first adjoining the class of $\gamma$ and then completing it to a basis of $\pi_1(G)/2\pi_1(G)$. Thus, we can define an element of $\Hom(\pi_1(G)/2\pi_1(G),\Z/2)$ which sends $\im j'_{+*}$ to 0 and $\gamma$ to 1. This morphism induces an element $\lambda\in\Hom(\pi_1(G),\Z/2)$ satisfying
\[
\lambda\circ j_{+*}=0,\qquad \lambda(\gamma)=1.
\]
For any $\psi\in S$, both $\psi$ and $\psi+\lambda$
belong to $S$, and their values on $\gamma$ differ.
Consequently, one of them satisfies
$\psi(\gamma)+\epsilon=0$, and hence $M$ is spin.
\end{proof}


We now study $G$-invariant spin structures. Note that, since the inclusions of $M_\pm^o\longrightarrow M$ are $G$-equivariant, if $M$ has a $G$-invariant spin structure, then $M_\pm^o$ and hence $M_\pm$ also have a $G$-invariant spin structure. Thus, as in the non-invariant case, we need to study when the existence of $G$-invariant spin structures on $M_\pm$ imply the existence of a $G$-invariant spin structure on $M$. 

In the $G$-invariant case, \cref{gluing spin structures}
and \cref{gluing spin structures a odd} give the following analogue:

\begin{proposition}\label{gluing G-inv spin structures}
Assume that $M$ is a closed connected manifold with a cohomogeneity one action of type (2) of a compact connected  Lie group $G$ where $K_+$ and $K_-$ are connected. Moreover, assume that one of the following conditions holds:
\begin{itemize}
    \item[1.] $l_+\geq 2$ or $l_-\geq 2$, or
    \item[2.] $l_+=l_-=1$ and $a$ is odd. 
\end{itemize}

If $M_+$ and $M_-$ have a $G$-invariant spin structure then $M$ has a $G$-invariant spin structure.
\end{proposition}

To prove \cref{gluing G-inv spin structures} we need an analogous statement to \cref{spin structure surgery} for invariant spin structures.

\begin{lemma}\label{spin structure surgery G-inv}
Let $M_1$ and $M_2$ be spin manifolds of dimension $n$ with connected boundaries $S_1$ and $S_2$ and let $G$ be a compact connected Lie group acting by orientation preserving isometries on $M_1$ and $M_2$. Let $h:S_2\longrightarrow S_1$ be a $G$-equivariant orientation-reversing diffeomorphism. If $J_G=\{(x,y)\in \Spin_G(M_1)\times \Spin_G(M_2):h^*(x_{|S_1})=y_{|S_2}\}$ is non-empty, then $M=M_1\cup_h M_2$ has a $G$-invariant spin structure.
\end{lemma}

\begin{proof}
Since $J_G\neq\emptyset$, there exist
$x\in\Spin_G(M_1)$ and $y\in\Spin_G(M_2)$
such that $h^*(x_{|S_1})=y_{|S_2}$.
Choose $G$-invariant spin structures
$\Lambda_i:P_i\longrightarrow FM_i$, $i=1,2$,
representing $x$ and $y$, respectively.

The outward unit normals along $S_i$ determine
reductions of the oriented frame bundles
$FM_i|_{S_i}$ to $FS_i$, and hence induce spin
structures on $S_i$. Since $G$ preserves the
outward normals, these induced spin structures
are also $G$-invariant.

The equality $h^*(x_{|S_1})=y_{|S_2}$
implies that, using compatible collars and
reversing the normal direction, there exists
an isomorphism of principal $\Spin(n)$-bundles
\[
\widetilde{h}:P_2|_{S_2}
\longrightarrow P_1|_{S_1}
\]
covering the corresponding identification of
the oriented frame bundles.

Since $h$ is $G$-equivariant, the induced map
on the oriented frame bundles is also
$G$-equivariant. For every $g\in G$, the maps
$\widetilde{h}$ and
$g\circ\widetilde{h}\circ g^{-1}$
therefore cover the same map and differ by
a deck transformation in $\Z/2$.
Since $G$ is connected, this deck transformation
must be trivial. Thus $\widetilde{h}$ is
$G$-equivariant.

Using compatible $G$-invariant collars, we can
glue $P_1$ and $P_2$ along $\widetilde{h}$ to
obtain a spin structure
$\Lambda:P\longrightarrow FM$ on $M$.
Since $\widetilde{h}$ is $G$-equivariant,
the lifted actions of $G$ on $P_1$ and $P_2$
glue to an action of $G$ on $P$ covering
the given action on $FM$.
Hence $P$ is a $G$-invariant spin structure
on $M$.
\end{proof}

\begin{proof}[Proof of \cref{gluing G-inv spin structures}]

Since the inclusion $\iota_\pm:M_\pm'\longrightarrow M_\pm$ and the diffeomorphism $\zeta_\pm:G/H\longrightarrow M_\pm'$ are $G$-equivariant, the restriction maps $\phi_\pm:\Spin(M_\pm)\longrightarrow\Spin(G/H)$ send $G$-invariant spin structures to $G$-invariant spin structures. We denote these maps by $\phi_\pm^G:\Spin_G(M_\pm)\longrightarrow\Spin_G(G/H)$. Hence, by \cref{spin structure surgery G-inv}, $M$ has a $G$-invariant spin structure if and only if there exist $x_\pm\in\Spin_G(M_\pm)$ such that $\phi_+^G(x_+)=\phi_-^G(x_-)$.

Recall that $G$-invariant spin structures on $G/H$ correspond to lifts of the isotropy representation of $H$ to the spin group. Since any two such lifts differ by a homomorphism $H\longrightarrow\Z/2$, there is at most one $G$-invariant spin structure on $G/H$ whenever $\Hom(H,\Z/2)=0$.

\textbf{Part 1: The case where $l_+\geq2$ or $l_-\geq2$.}
Assume that $l_+\geq2$ (the other case is analogous). Since $K_+$ is connected and $K_+/H\cong S^{l_+}$ is simply connected, the long exact sequence of homotopy implies that $H$ is connected. Hence $\Hom(H,\Z/2)=0$, and therefore $G/H$ admits at most one $G$-invariant spin structure. Since $M_+$ and $M_-$ have $G$-invariant spin structures, their restrictions to $G/H$ must agree. By \cref{spin structure surgery G-inv}, $M$ admits a $G$-invariant spin structure.

\textbf{Part 2: The case where $l_+=l_-=1$ and $a$ is odd.}
Recall that $H/H_0\cong\Z/a$. Since $a$ is odd, every homomorphism $H\longrightarrow\Z/2$ is trivial. Thus, also in this case, $G/H$ admits at most one $G$-invariant spin structure. Consequently, the restrictions of the $G$-invariant spin structures on $M_+$ and $M_-$ agree, and $M$ admits a $G$-invariant spin structure by \cref{spin structure surgery G-inv}.
\end{proof}

\begin{remark}
Assume that $l_\pm=1$, that $a$ is odd and that
$M_\pm$ admits a $G$-invariant spin structure.
Then
$\sigma_{G/K_\pm*}=\delta_2\circ q_{\pm*}\neq0$.
Consequently, the unique $G$-invariant spin structure
on $M_\pm'$ is not $G\times S^1$-invariant.
In particular, when
$\delta_2\circ q_{\pm*}\in\im j_\pm^*$,
the spin structures in $S_{\pm0}$ cannot be
$G$-invariant.
\end{remark}
\begin{proposition}
\label{gluing G-invariant spin structures a even}
Assume that $M$ is a closed connected manifold with a cohomogeneity one action of type (2) of a compact connected Lie group $G$ where $K_+$ and $K_-$ are connected. Assume that $l_+=l_-=1$ and $a$ is even. If $M_+$ and $M_-$ have a $G$-invariant spin structure then $M$ has a $G$-invariant spin structure if and only if $\epsilon=0$.
\end{proposition}

\begin{proof}
Since $a$ is even, we have $\delta_2\circ q_{\pm*}=0$, and hence $M_\pm$ admit $G$-invariant spin structures precisely when $\sigma_{G/K_\pm*}=0$. Since $K_\pm$ are connected, the unique $G$-invariant spin structures on $M_\pm$ are represented by the zero character, and their restrictions to $M_\pm'$ are represented by $(0,\delta_2)$. Set $f_\pm=\zeta_\pm^*(0,\delta_2)$. Both $f_+$ and $f_-$ pull back to the spin structure represented by the zero character on $G/H_0$. They therefore agree if and only if they take the same value on $\gamma_-$, since the difference of two spin structures with the same pullback is either $0$ or $\eta$, and $\eta(\gamma_-)=1$.

By the definition of $\gamma_\pm$, we have $f_+(\gamma_+)=f_-(\gamma_-)=1$. Since $a$ is even and $\partial\gamma_\pm$ are generators of $\Z/a$, the integer $r$ is odd. Thus, using $r\gamma_+-\gamma_-=\rho_*(\gamma,\epsilon)$, we obtain
$$f_+(\gamma_-)=r f_+(\gamma_+)-f_+(\rho_*(\gamma,\epsilon))=1+\epsilon.$$
Thus $f_+=f_-$ if and only if $\epsilon=0$. By \cref{spin structure surgery G-inv}, $M$ admits a $G$-invariant spin structure precisely when $\epsilon=0$.
\end{proof}


We finish this section by providing an example of an oriented non-spin manifold $M$ admitting a cohomogeneity one action of a compact Lie group $G$ such that $M_+$ and $M_-$ are spin, showing that the additional hypotheses in \cref{gluing spin structures} cannot be omitted in general.

For $n\geq 2$, consider the group $G=\SO(2n+2)$ and the subgroups $K_+=\U(n)$, $K_-=\SU(n)\times S^1$ and $H=\SU(n)$. We consider natural inclusions $$i_+:\SU(n)\longrightarrow\U(n)\qquad\text{and}\qquad i_-:\SU(n)\longrightarrow \SU(n)\times S^1.$$ We also define inclusions
$$j_+:\U(n)\longrightarrow\SO(2n)\longrightarrow\SO(2n+2)\qquad\text{and}\qquad j_-:\SU(n)\times S^1\longrightarrow \SO(2n)\times \SO(2)\longrightarrow \SO(2n+2).$$
Using the above inclusions we can define $M_+=G\times_{K_+}\D^2$ and $M_-=G\times_{K_-}\D^2$, and $M=M_+\cup_{G/H}M_-$ is a closed connected manifold admitting a cohomogeneity one action of $G$. Since all Lie groups $G$, $K_\pm$ and $H$ are connected, $M$ is orientable.

\begin{lemma}
The manifolds $M_+$ and $M_-$ are spin. The manifold $M$ is spin if and only if $n$ is odd.
\end{lemma}

\begin{proof}
By \cref{main_theorem noncompact section}, $M_+$ is spin if and only if there exists $\psi_+\in \Hom(\pi_1(G),\Z/2)=\Hom(\Z/2,\Z/2)=\{0,\id\}$ satisfying $\psi_+\circ j_{+*}=\sigma_{G/K_+*}+\delta_2\circ q_{+*}$. The manifold $\SO(2n)/\U(n)$ is a symmetric space, and we know by \cite{cahen1988spin} that $\sigma_{\SO(2n)/\U(n)*}=0$ if $n$ is odd and $\sigma_{\SO(2n)/\U(n)*}=\delta_2$ if $n$ is even. Similarly, we also have $\sigma_{\SO(2n+2)/\SO(2n)*}=0$. Using \cref{hf isotoropy representation}, we obtain that $\sigma_{G/K_+*}=0$ if $n$ is odd and $\sigma_{G/K_+*}=\delta_2$ if $n$ is even.

On the other hand, $q_+:\U(n)\longrightarrow S^1$ induces the identity map on the fundamental group and $j_{+*}=\delta_2$. Thus, if $n$ is odd then $M_+$ is spin if and only if there exists $\psi_+\in\Hom(\pi_1(G),\Z/2)$ satisfying $\psi_+\circ\delta_2=\delta_2$. We can choose $\psi_+=\id$ satisfying the previous equation. Similarly, if $n$ is even then $M_+$ is spin if and only if there exists $\psi_+\in\Hom(\pi_1(G),\Z/2)$ satisfying $\psi_+\circ\delta_2=0$. We can choose $\psi_+=0$ satisfying the previous equation. Hence, in both cases $M_+$ is spin.

In the same way as above, $M_-$ is spin if and only if there exists $\psi_-\in \Hom(\pi_1(G),\Z/2)$ satisfying $\psi_-\circ j_{-*}=\sigma_{G/K_-*}+\delta_2\circ q_{-*}$. Since $\sigma_{\SO(2n+2)/(\SO(2n)\times\SO(2))*}=0$ (see \cite{cahen1988spin}) and $\SU(n)$ is simply connected, we can use \cref{hf isotoropy representation} to conclude that $\sigma_{G/K_-*}=0$. We also have that $q_-:\SU(n)\times S^1\longrightarrow S^1$ induces the identity morphism on the fundamental groups and $j_{-*}=\delta_2$. Thus, $M_-$ is spin if and only if there exists $\psi_-\in\Hom(\pi_1(G),\Z/2)$ satisfying $\psi_-\circ\delta_2=\delta_2$. We can choose $\psi_-=\id$ and hence $M_-$ is spin.

Note that the spin structures on $M_+$ and $M_-$ are unique, since $j_{\pm*}$ are surjective and hence $G/K_\pm$ are simply connected. Moreover, $\SO(2n+2)/\SU(n)$ has two spin structures. Since $H$ is connected, we can use the description of $\pi^*\circ\phi_\pm$ established above. In this case, $H=H_0$, so $\Spin(G/H)=\{0,\id\}$ and $\phi_\pm(\psi_\pm)=\psi_\pm$. Assume that $n$ is odd. Then $\psi_+=\id$ and therefore $\phi_+(\psi_+)=\id$. Moreover, $\phi_-(\psi_-)=\id$, and therefore $M$ is spin by \cref{spin structure surgery}.

On the other hand, if $n$ is even then $\psi_+=0$ and $\psi_-=\id$, which implies that $\phi_+(\psi_+)\neq \phi_-(\psi_-)$ and hence $M$ does not admit a spin structure.
\end{proof}

\begin{remark}
The same example also illustrates the necessity of the
hypothesis that $a$ is odd in
\cref{gluing G-inv spin structures}.
Indeed, replace $G=\SO(2n+2)$ by
$\widetilde G=\Spin(2n+2)$, acting on the same
manifold through the double covering, and let
$\widetilde K_\pm$ and $\widetilde H$ be the
preimages of $K_\pm$ and $H$.
Then $\widetilde K_\pm$ are connected, while
$\widetilde H$ has two connected components,
so $a=2$.
Since $\widetilde G$ is simply connected,
the spin structures on $M_\pm$ are
$\widetilde G$-invariant.
However, for even $n$, the manifold $M$
is not spin, and hence cannot admit a
$\widetilde G$-invariant spin structure.
\end{remark}

\section{Examples}\label{sec: examples}

\begin{example}
Let $n\geq 3$. We have a cohomogeneity one action of
$\SO(n)$ on $S^n$, obtained by restricting the transitive
action of $\SO(n+1)$ on $S^n$. The two singular orbits are
fixed points, and we can decompose $S^n$ into two disks
\[
M_\pm=\SO(n)\times_{\SO(n)}\D^n\cong\D^n.
\]
Since $M_\pm$ are spin and $l_\pm=n-1\geq2$,
by \cref{gluing spin structures}, $S^n$ is spin.
Together with the fact that $S^2$ is spin, this shows
that all spheres of dimension at least $2$ are spin.
\end{example}

\begin{example}
It is well known that $\CP^n$ is spin if and only if $n$ is odd.
Let us provide another proof using \cref{gluing spin structures}.
Let $n\geq 2$ and consider the action of $\SU(n)$ on $\CP^n$
defined by
\[
A[z_0:z_1:\dots:z_n]
=[z_0:A(z_1,\dots,z_n)^t].
\]
It is straightforward to see that this is a cohomogeneity one
action with $\CP^n/\SU(n)\cong[-1,1]$.
One singular orbit is the fixed point $[1:0:\dots:0]$,
hence $K_-=\SU(n)$. The other singular orbit is
$\CP^{n-1}\cong\{[0:z_1:\dots:z_n]\in\CP^n\}$,
and hence $K_+=\U(n-1)$.
The principal orbits are spheres
$\SU(n)/\SU(n-1)\cong S^{2n-1}$.

Since $\SU(n)\times_{\SU(n)}\D^{2n}\cong\D^{2n}$
is spin and $l_-=2n-1\geq2$, by
\cref{gluing spin structures}, $\CP^n$ is spin if and only if
$\SU(n)\times_{\U(n-1)}\D^2$ is spin.

We claim that $\SU(n)\times_{\U(n-1)}\D^2$ is spin
if and only if $\CP^{n-1}$ is not spin.
Indeed, the map $q:\U(n-1)\longrightarrow S^1$
induces an isomorphism on fundamental groups, while
$j_*$ is trivial since $\SU(n)$ is simply connected.
Consequently, there does not exist
$\eta\in\Hom(\pi_1(\SU(n)),\Z/2)$ satisfying
$\eta\circ j_*=\delta_2\circ q_*$.
By item (3)(b) of \cref{main_theorem noncompact section},
if $\SU(n)\times_{\U(n-1)}\D^2$ is spin,
then $\CP^{n-1}$ is not spin.

Conversely, assume that $\CP^{n-1}$ is not spin.
Since $\pi_1(\U(n-1))\cong\Z$, the map
$\sigma_{\SU(n)/\U(n-1)*}$ is either trivial or
$\delta_2$. If it were trivial, then $\CP^{n-1}$
would be spin by \cref{spin structures homogeneous spaces}.
Therefore,
\[
\sigma_{\SU(n)/\U(n-1)*}=\delta_2\circ q_*,
\]
and hence
\[
\sigma_{\SU(n)/\U(n-1)*}+\delta_2\circ q_*=0.
\]
By the proof of \cref{main_theorem noncompact section},
$\SU(n)\times_{\U(n-1)}\D^2$ is spin,
which completes the proof of the claim.

In conclusion, $\CP^n$ is spin if and only if
$\CP^{n-1}$ is not spin. Since $\CP^1\cong S^2$
is spin, we conclude that $\CP^n$ is spin if and only if
$n$ is odd.
\end{example}

\begin{example}\label{cohomogeneity 1 Grassmannians}
Let $n_1,n_2\geq2$ and set $G=\SO(n_1+n_2+1)$.
Consider the chains of inclusions
\[
H=\SO(n_1)\times\SO(n_2)
\subseteq K_1=\SO(n_1+1)\times\SO(n_2)
\subseteq G
\]
and
\[
H\subseteq K_2=\SO(n_1)\times\SO(n_2+1)
\subseteq G.
\]
They define a cohomogeneity one manifold
\[
M=X_1\cup_{G/H}X_2,
\qquad
X_1=G\times_{K_1}\D^{n_1+1},
\quad
X_2=G\times_{K_2}\D^{n_2+1}.
\]
With the natural gluing, $M$ is the oriented real
Grassmannian
\[
M\cong
\frac{\SO(n_1+n_2+2)}
{\SO(n_1+1)\times\SO(n_2+1)},
\]
with the action of $G$ induced by its standard
inclusion in $\SO(n_1+n_2+2)$.

Since $n_1,n_2\geq2$, both singular orbits have
codimension at least three, so $M$ is spin if and
only if its principal orbit $G/H$ is spin.
The isotropy representation of $G/H$ is
\[
\R^{n_1}\otimes\R^{n_2}
\oplus\R^{n_1}\oplus\R^{n_2}.
\]
Let $e_i$ denote a generator of $\pi_1(\SO(n_i))$,
which is $\Z$ when $n_i=2$ and $\Z/2$ otherwise.
The inclusion $j:H\longrightarrow G$ satisfies
$j_*(e_1,0)=j_*(0,e_2)=1\in\Z/2$, while
\[
\sigma_{G/H*}(e_1,0)=n_2+1\mod2,
\qquad
\sigma_{G/H*}(0,e_2)=n_1+1\mod2.
\]
By \cref{spin structures homogeneous spaces},
the isotropy character factors through $j_*$
if and only if $n_1+n_2$ is even. Consequently,
\[
M\text{ is spin}
\quad\Longleftrightarrow\quad
n_1+n_2\text{ is even}.
\]

We can also determine when the spin structure
is $G$-invariant. Since $K_i/H$ is simply connected,
the maps $\pi_1(H)\longrightarrow\pi_1(K_i)$
are surjective. Thus, the invariant spin
obstructions for the disk bundles vanish if
and only if $\sigma_{G/H*}=0$.
By \cref{gluing G-inv spin structures},
$M$ therefore admits a $G$-invariant spin
structure precisely when both $n_1$ and $n_2$
are odd. 
\Cref{spin structures homogeneous spaces} together with standard computations imply that $M$ is $\SO(n_1+n_2+2)$-invariant if and only if $n_1+1$ and $n_2+1$ are even. Thus, $M$ has a $\SO(n_1+n_2+2)$-invariant spin structure if and only if it has a $\SO(n_1+n_2+1)$-invariant spin structure.

We now consider the case $n_1=n\geq3$ and $n_2=1$.
Here the principal orbit is the Stiefel manifold
$\SO(n+2)/\SO(n)$, which is spin, but one singular
orbit has codimension two. The disk bundles are
\[
X_1=\SO(n+2)\times_{\SO(n+1)}\D^{n+1},
\qquad
X_2=\SO(n+2)\times_{\SO(n)\times\SO(2)}\D^2.
\]
The first is spin. For the second, the character
$\sigma_{G/K_2*}+\delta_2\circ q_{2*}$ takes the
values $0$ and $n+1\mod2$ on the generators of
$\pi_1(K_2)=\Z/2\times\Z$, whereas $j_{2*}$ sends
both generators to $1\in\Z/2$.
By \cref{main_theorem noncompact section},
$X_2$ is spin if and only if $n$ is odd.
Consequently, $M$ is spin precisely when $n$
is odd, and in this case the spin structure
is also $G$-invariant.
The case $n=2$ is considered explicitly below,
using the double cover
$\Spin(4)\cong\SU(2)\times\SU(2)$.
\end{example}

\begin{example}
We illustrate the preceding codimension-two case
in the smallest dimension, using the action of
$\SU(2)\times\SU(2)\cong\Spin(4)$.

The oriented real Grassmannian
$Gr_2^+(\R^5)=\SO(5)/(\SO(2)\times\SO(3))$
admits a cohomogeneity one action of
$\SU(2)\times\SU(2)$ which induces the chains
of inclusions
\[
\Delta S^1\subseteq T^2
\subseteq\SU(2)\times\SU(2)
\]
and
\[
\Delta S^1\subseteq\Delta\SU(2)
\subseteq\SU(2)\times\SU(2).
\]

Since $\Delta\SU(2)/\Delta S^1\cong S^2$ and
$(\SU(2)\times\SU(2))/\Delta S^1$ is spin,
we know by \cref{main_theorem noncompact section}
that
$(\SU(2)\times\SU(2))
\times_{\Delta\SU(2)}\D^3$
is spin.

Therefore, to study the spin structure of
$Gr_2^+(\R^5)$, we need to study the first
chain of inclusions.
Since $T^2/\Delta S^1\cong S^1$, we compute
$q_*$ and $\sigma_{\SU(2)\times\SU(2)/T^2*}$,
as in the proof of
\cref{main_theorem noncompact section}.
Indeed, $q_*:\Z^2\longrightarrow\Z$ satisfies
\[
q_*(1,0)=1,\qquad q_*(0,1)=-1.
\]
Moreover,
$(\SU(2)\times\SU(2))/T^2\cong S^2\times S^2$,
and each $S^1$ factor of $T^2$ acts on the
corresponding tangent plane with weight $2$.
Therefore,
$\sigma_{\SU(2)\times\SU(2)/T^2*}=0$.

Hence,
\[
\delta_2\circ q_*+
\sigma_{\SU(2)\times\SU(2)/T^2*}\neq0.
\]
On the other hand, $j_*=0$ since
$\SU(2)\times\SU(2)$ is simply connected.
Therefore,
$(\SU(2)\times\SU(2))\times_{T^2}\D^2$
does not admit a spin structure, and hence
$Gr_2^+(\R^5)$ is not spin, as also shown
in \cite{cahen1988spin}.
\end{example}

\begin{example}
The construction in \cref{cohomogeneity 1 Grassmannians} has complex and
quaternionic analogues.

For the complex case, let $n_1,n_2\geq1$ and
$G=\U(n_1+n_2+1)$. Consider
\[
H=\U(n_1)\times\U(n_2)
\subseteq K_1=\U(n_1+1)\times\U(n_2)
\subseteq G
\]
and
\[
H\subseteq K_2=\U(n_1)\times\U(n_2+1)
\subseteq G.
\]
The corresponding cohomogeneity one manifold is
\[
M=X_1\cup_{G/H}X_2,
\qquad
X_1=G\times_{K_1}\D^{2n_1+2},
\quad
X_2=G\times_{K_2}\D^{2n_2+2}.
\]
This is the complex Grassmannian
\[
M\cong
\frac{\U(n_1+n_2+2)}
{\U(n_1+1)\times\U(n_2+1)}.
\]

Since both singular orbits have codimension
at least four, $M$ is spin if and only if
$G/H$ is spin. The inclusion
$j:H\longrightarrow G$ induces
$j_*:\Z^2\longrightarrow\Z$ with
$j_*(1,0)=j_*(0,1)=1$, while
\[
\sigma_{G/H*}(1,0)=n_2+1\mod2,
\qquad
\sigma_{G/H*}(0,1)=n_1+1\mod2.
\]
As in the real case, the isotropy character
factors through $j_*$ if and only if
$n_1+n_2$ is even. Thus $M$ is spin precisely
when $n_1+n_2$ is even, and admits a
$G$-invariant spin structure precisely when
$n_1$ and $n_2$ are both odd.

For the quaternionic case, let $n_1,n_2\geq1$,
set $G=\Sp(n_1+n_2+1)$, and consider
\[
H=\Sp(n_1)\times\Sp(n_2)
\subseteq K_1=\Sp(n_1+1)\times\Sp(n_2)
\subseteq G
\]
and
\[
H\subseteq K_2=\Sp(n_1)\times\Sp(n_2+1)
\subseteq G.
\]
The corresponding manifold is
\[
M=Y_1\cup_{G/H}Y_2,
\qquad
Y_1=G\times_{K_1}\D^{4n_1+4},
\quad
Y_2=G\times_{K_2}\D^{4n_2+4},
\]
which is the quaternionic Grassmannian
\[
M\cong
\frac{\Sp(n_1+n_2+2)}
{\Sp(n_1+1)\times\Sp(n_2+1)}.
\]
Since both $G$ and $H$ are simply connected,
$G/H$ admits a $G$-invariant spin structure.
Consequently, by \cref{gluing G-inv spin structures},
$M$ admits a $G$-invariant spin structure
for all $n_1,n_2\geq1$.
\end{example}

\begin{example}
\Cref{main_theorem noncompact section} can also
be used to study spin structures on homogeneous
spaces. The oriented real Grassmannian
$Gr_2^+(\R^9)$ admits a cohomogeneity one action
of $\Spin(7)$ with associated chains of inclusions
\[
\SU(3)\subseteq G_2\subseteq\Spin(7)
\]
and
\[
\SU(3)\subseteq S^1\cdot\SU(3)\subseteq\Spin(7).
\]

Since $Gr_2^+(\R^9)$ is not spin and
$\Spin(7)\times_{G_2}\D^7$ is spin
(recall that $G_2/\SU(3)\cong S^6$),
we conclude that
$\Spin(7)\times_{S^1\cdot\SU(3)}\D^2$
is not spin.

Recall that
$S^1\cdot\SU(3)\cong(S^1\times\SU(3))/(\Z/3)$.
Since $\SU(3)$ is simply connected,
$q_*:\Z\longrightarrow\Z$ is the identity map
and $\rho\circ q_*=\rho$.
Moreover, $j_*$ is trivial, and since
$\pi_1(S^1\cdot\SU(3))\cong\Z$,
the character
$\sigma_{\Spin(7)/(S^1\cdot\SU(3))*}$
is either trivial or $\rho$.
If
$\sigma_{\Spin(7)/(S^1\cdot\SU(3))*}=\rho$,
then
\[
\sigma_{\Spin(7)/(S^1\cdot\SU(3))*}
+\rho\circ q_*=0
\]
and $\Spin(7)\times_{S^1\cdot\SU(3)}\D^2$
would be spin, which is a contradiction.
Therefore,
$\sigma_{\Spin(7)/(S^1\cdot\SU(3))*}=0$,
and $\Spin(7)/(S^1\cdot\SU(3))$ is spin.
\end{example}

\begin{example}
For a non-negative integer $p$, the Bazaikin space
\[
B_p^{13}=
\diag(z,z,z,z,z^{2p-1})\setminus \SU(5)/
\Sp(2)\diag(1,1,1,1,\overline{z}^{2p+3})
\]
admits a cohomogeneity one action by $\SU(4)$
(see \cite{grove2008positively}). The isotropy groups are
$K_+=\Sp(2)\cdot\Z/2$,
$K_-=\SU(2)\cdot S^1_{p,p+1}$ and
$H=\SU(2)\cdot\Z/2$.
In particular, $K_+/H\cong S^7$,
$K_-/H\cong S^1$ and $G/K_+\cong\RP^5$.
Thus, we have the homogeneous fibration
\[
S^7\longrightarrow G/H\longrightarrow\RP^5.
\]

We claim that $G/H$ is not spin.
Indeed, the involution
$\iota=\diag(1,-1,-1,1)\in H$
acts on $\mathfrak{su}(4)/\mathfrak{su}(2)$
with a six-dimensional $(-1)$-eigenspace.
Therefore, any lift of its isotropy action to
$\Spin(12)$ has square $-1$, and hence the
isotropy representation of $H$ does not lift
to $\Spin(12)$.
Since $G=\SU(4)$ is simply connected, this
implies that $G/H$ is not spin.

Consequently, $B_p^{13}$ is not spin for any $p$,
since a principal orbit $G/H$ is an oriented
hypersurface with trivial normal bundle.
In particular, this includes the homogeneous
Berger space $B_1^{13}$.
\end{example}
\begin{example}
Consider the Brieskorn manifolds
\[
W_d^{2n-1}=
\left\{(z_0,Z)\in\C\oplus\C^n:
z_0^d+\sum_{i=1}^n z_i^2=0,\ 
|z_0|^2+|Z|^2=1\right\},
\]
where $n\geq3$ and $d\geq1$. If $n$ and $d$ are odd then $W^{2n-1}_d$ is homeomorphic to a sphere. Moreover, it is exotic if $d\equiv\pm3\pmod8$ and $2n-1\equiv1\pmod8$, see see \cite[p.~343]{grove2000curvature}. Thus, when $n$ and $d$ are odd they admit a unique spin structure.

From now on we assume that $n$ and $d$ are odd. There exists a cohomogeneity one action of $G=S^1\times\SO(n)$ on $W_d^{2n-1}$ given by
\[
(e^{i\theta},A)(z_0,Z)
=(e^{2i\theta}z_0,e^{id\theta}AZ).
\]
The group diagram, described in \cite{grove2006nonnegative}, has principal isotropy
\[
H=\{(\varepsilon,\diag(\varepsilon I_2,A)):\varepsilon=\pm1,\ A\in\SO(n-2)\}
\]
and singular isotropy groups
\[
K_-=\{(e^{-i\theta},\diag(R(d\theta),A)):
A\in\SO(n-2)\},
\]
and
\[
K_+=
\{(\det B,\diag(\det B,B)):B\in O(n-1)\}.
\]
Here $K_-/H\cong S^1$, $K_+/H\cong S^{n-2}$, and $H$ has two connected components. We note that whenever $W_d^{2n-1}$ is exotic then it is diffeomorphic to an exotic Kervaire sphere of dimension $2n-1$. By \cite{straume1994compact}, if a compact connected Lie group $G$ acts effectively on $W^{2n-1}_d$ then $\dim G\leq1+\tfrac{n(n-1)}{2}$. Since the
cohomogeneity one action above is effective and $\dim S^1\times\SO(n)=1+\frac{n(n-1)}{2}$, it realizes this bound.

We claim that the unique spin structure on of $W_d^{2n-1}$ is not invariant by the above cohomogeneity one action. We first consider the disk bundle $M_-=G\times_{K_-}\D^2$. Since $K_-$ is connected, we can apply \cref{main_theorem noncompact invariant section}. Using the decomposition $\mathfrak{so}(n)=\mathfrak{so}(2)\oplus\mathfrak{so}(n-2)\oplus(\R^2\otimes\R^{n-2})$ and taking the quotient by $\mathfrak k_-=\R(-1,d)\oplus\mathfrak{so}(n-2)$, we obtain
\[
\mathfrak g/\mathfrak k_-
\cong\R\oplus(\R^2\otimes\R^{n-2})
\]
as $K_-$-representations.

The circle factor of $K_-$ acts on
$\R^2\otimes\R^{n-2}$ as the direct sum of
$n-2$ rotations of weight $d$. Hence
$\sigma_{G/K_-*}$ takes the value
$d(n-2)\mod2$ on its generator.
On the other hand, the $\SO(n-2)$ factor acts
through two copies of its standard representation,
so its contribution to $\sigma_{G/K_-*}$ vanishes.
Finally, the normal disk is parametrised by $z_0$,
on which the circle acts by
$z_0\mapsto e^{-2i\theta}z_0$,
while $\SO(n-2)$ acts trivially. Thus $ev_*=0$.
Consequently, by
\cref{main_theorem noncompact invariant section},
$M_-$ admits a $G$-invariant spin structure
if and only if
\[
d(n-2)\equiv nd\equiv0\pmod2.
\]
Since a $G$-invariant spin structure on $W_d^{2n-1}$ restricts to one on $M_-$,
we conclude that $W_d^{2n-1}$ cannot admit a $G$-invariant spin structure, since $d$ and $n$ are odd. In particular, the nine-dimensional exotic Kervaire sphere $W_3^9$ is spin, but its natural $S^1\times\SO(5)$-action does not lift to its spin structure.

This example can be understood as a counterpart that the unique spin structure of $S^m$ is not invariant by the compact connected Lie group acting effectively on $S^m$ of largest dimension, namely $\SO(m+1)$. 
\end{example}

\bibliographystyle{alpha}
\bibliography{refs.bib}

\end{document}